\documentclass[10pt]{article}

\usepackage{epic,eepic,epsfig,amssymb,amsmath,amsthm,graphics,stmaryrd}

\usepackage{float}
\usepackage{authblk}
\usepackage{multimedia}
\usepackage{xcolor}
\usepackage{tikz}
\usepackage{algorithm2e}
\usepackage{pgfplots}
\usepackage{makecell}            
\usepackage[margin=1in]{geometry}
\usepackage{longtable}

\AtEndDocument{\bigskip{\footnotesize%
  \textsc{Laboratoire J.A. Dieudonné, Université Côte d'Azur, Parc Valrose, 06108 Nice Cedex 2, France} \par  
  \textit{E-mail address}: \texttt{ludovic.rifford@math.cnrs.fr} 
}}
  
\def\R{\mathbb R}
\def\N{\mathbb N}
\def\Z{\mathbb Z}
\def\Q{\mathbb Q}
\def\I{\mathbb I}
\def\b{\mbox{b}}
\def\k{\mbox{k}}
\def\c{\mbox{c}}
\def\m{\hspace{-0.075cm}\shortmid \hspace{-0.075cm}}

\newtheorem{theorem}{Theorem}
\newtheorem{lemma}{Lemma}
\newtheorem{proposition}{Proposition}

\newtheorem{definition}{Definition\rm}
\newtheorem{remark}{Remark}

\makeindex

\begin{document}

\title{On the time for a runner to get lonely}


\author{Ludovic Rifford}



\maketitle

\begin{abstract}
The Lonely Runner Conjecture asserts that if $n$ runners with distinct constant speeds run on the unit circle $\R/\Z$ starting from $0$ at time $0$, then each runner will at some time $t>0$ be lonely in the sense that she/he will be separated by a distance at least $1/n$ from all the others at time $t$. In investigating the size of $t$, we show that an upper bound for $t$ in terms of a certain number of rounds (which, in the case where the lonely runner is static, corresponds to the number of rounds of the slowest non-static runner) is equivalent to a covering problem in dimension $n-2$. We formulate a conjecture regarding this covering problem and prove it to be true for $n=3,4,5,6$. Then, we use our method of proof to demonstrate that the Lonely Runner Conjecture with Free Starting Points is satisfied for $n=3,4$. Finally, we show that the so-called gap of loneliness in one round (with respect to the Lonely Runner Conjecture), where we have $m+1$ runners including one static runner,  is bounded from below by $1/(2m-1)$  for all integer $m\geq 2$.
\end{abstract}

\section{Introduction}\label{SECintroduction}

The Lonely Runner Conjecture, which was stated for the first time in this form  by Bienia {\it et al.} in \cite{bggst98}, formulates in terms of runners running on the circular track $\R / \Z$  a conjecture introduced originally by Wills \cite{wills67} in the context of diophantine approximation and by Cusick \cite{cusick73} in obstruction theory  ($\{x\}$  denotes the fractional part of $x$, that is,  $\{x\}=x-\lfloor x\rfloor$): \\

\noindent {\bf Lonely Runner Conjecture.}
{\em For every integer $n\geq 2$, if $x_1, \ldots, x_n:[0,\infty) \rightarrow \R /\Z$ are $n$ runners with different constant speeds $v_1, \ldots, v_n \in \R$ starting from $0$ at time $0$, that is such that $x_i(t)=\left\{ tv_i\right\}$ for all $t\geq 0$, then for each runner there is a time where she/he is at distance at least $1/n$ from the others, which in other words means that for every $i\in \{1, \ldots,n\}$, there is  $t>0$ such that 
$$
\left|x_j(t)-x_i(t) \right|  \in \left[ \frac{1}{n}, 1-\frac{1}{n}\right] \qquad \forall j\in \bigl\{1,\ldots,n\bigr\} \setminus \bigl\{i\bigr\}.
$$}

\begin{figure}[H]
\centering
\begin{tikzpicture}[scale=0.4]
\draw (0,0) grid(13,5);
\foreach \x in {0,...,13}
{
	\draw (\x,-0.25) node[below]{\tiny \color{black}{$\x$}};
}
\foreach \y in {0,...,5}
{
	\draw (-0.25,\y) node[left]{\tiny \color{black}{$\y$}};
}
\foreach \x in {0,...,13}
{
	\foreach \y in {0,...,5}
	{
		\filldraw[black] (\x,\y)circle(1pt);
	}
}
\foreach \x in {0,...,12}
{
	\foreach \y in {0,...,4}
	{
		\fill[color=red] (\x+0.33,\y+0.33) -- (\x+0.66,\y + 0.33) -- (\x+0.66,\y+0.66) -- (\x+0.33,\y+0.66) -- cycle;  
	}
}
\draw [color=blue] (0,0) -- (11.5,5);
\draw [color=blue] (0,0) -- (13,3.345);
\draw [color=blue] (0,0) -- (2.123,5);
\end{tikzpicture}
\caption{The Lonely Runner Conjecture in obstruction form asserts that every blue ray starting from the origin must intersect one of the red sets $\mathcal{K}^m_k(1/(m+1))$ (with $k\in \N^m$)\label{FigLight}}
\end{figure}
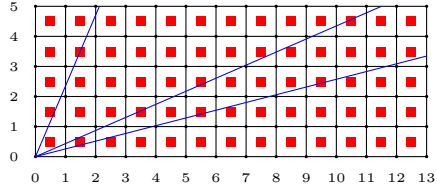

One can observe that for any $v, v' \in \R$, any $t>0$ and any $\delta \in (0,1/2]$, one has that $|\{tv\}-\{tv'\}|$ belongs to $[\delta,1-\delta]$ if and only if   $\{t|v'-v|\}$ belongs to   $[\delta,1-\delta]$ (see Lemma \ref{LEMLRCequivLRCOF}). This result allows to show that the Lonely Runner Conjecture holds true for $n\geq 2$ if and only if for every $w_1, \ldots, w_{n-1} >0$ there is $t>0$ such that $\{tw_i\}$ belongs to $[1/n,1-1/n]$ for all $i=1, \ldots, n-1$ (see Proposition \ref{PROPLRCequivLRCOF}). For every $\delta \in(0,1/2]$, we define the set $\mathcal{K}(\delta) \subset (0,\infty)$ by (we denote by $\N =\{0,1,2, \ldots\}$ the set of non-negative integers and by $\N^*=\N \setminus \{0\}$ the set of positive integers)
$$
\mathcal{K}(\delta) := \bigcup_{k\in \N} \mathcal{K}_k(\delta)  \quad \mbox{with} \quad \mathcal{K}_k(\delta):=  k+ \left[ \delta,1-\delta\right] \quad \forall k \in \N,
$$
and denote its $m$-th power (with $m\in \N^*$) by
$$
\mathcal{K}^m(\delta) =  \bigcup_{k\in \N^m}\mathcal{K}^m_k(\delta)    \quad \mbox{with} \quad\mathcal{K}^m_k(\delta):=  k+ \left[ \delta,1-\delta\right]^m \quad \forall k=(k_1, \ldots, k_m) \in \N^m.
$$ 
Then, given  $m\in \N^*$ and $m$ real numbers $w_1, \ldots, w_m>0$, $\{tw_1\}, \ldots, \{tw_m\}$ belong to $[1/(m+1),1-1/(m+1)]$ for some $t>0$ if and only if the ray $\{t(w_1, \ldots, w_m) \, \vert \, t>0\}$ intersects the set $\mathcal{K}^m(1/(m+1))$. Thus, we infer that the Lonely Runner Conjecture is equivalent to the following conjecture in obstruction theory (see Figure \ref{FigLight}):\\

\noindent {\bf Lonely Runner Conjecture - Obstruction Form.}
{\em Let $m \geq 1$ be an integer, then for any real numbers $w_1, \ldots, w_{m} >0$, there exists $t>0$ such that:
\begin{eqnarray*}
t  \left(w_1, \ldots, w_m\right) \in \mathcal{K}^m(1/(m+1)).
\end{eqnarray*}}

The above conjecture is obvious for $m=1$. For $m\geq 2$, by symmetry and dilation the conjecture holds true if and only if for any $w_1, \ldots, w_{m} >0$ with 
$$
1=w_1 \leq w_2,  \ldots , w_m,
$$
there is $\lambda >0$ such that 
$$
\lambda = \lambda w_1 \in \mathcal{K}(1/(m+1)) \quad \mbox{and} \quad \lambda \left(w_2, \ldots, w_m\right)  \in \mathcal{K}^{m-1}(1/(m+1)),
$$
which means that each coordinate of $(w_2,\ldots, w_m)\in [1,\infty)^{m-1}$ can be written as the quotient of an element of $ \mathcal{K}^{m-1}(1/(m+1))$ with the same element of $\mathcal{K}(1/(m+1))$.
For every integer $d\geq 1$, given two sets $A \subset \R^{d}$ and $B\subset (0,\infty)$, we use the notation $A/B$ to denote the subset of $\R^{d}$ given by
$$
A/B := \left\{ \frac{a}{b} =  \left(\frac{a_1}{b}, \ldots, \frac{a_d}{b}\right) \, \vert \, a=\left(a_1, \ldots, a_d\right) \in A, \, b \in B\right\}.
$$
If $B$ is a singleton, that is, of the form $\{b\}$, then we also write $A/b$ for $A/\{b\}$. Furthermore, from now on, for every integer $d\geq 1$, we set 
$$
\delta_d := \frac{1}{d+2}.
$$
The above discussion shows that the Lonely Runner Conjecture  is equivalent to the following conjecture:\\

\noindent {\bf Lonely Runner Conjecture - Covering Form.}
{\em For every integer $d\geq 1$, there holds 
$$
[1,\infty)^d \subset \mathcal{K}^{d}\left( \delta_d \right) / \mathcal{K} \left( \delta_d \right).
$$}

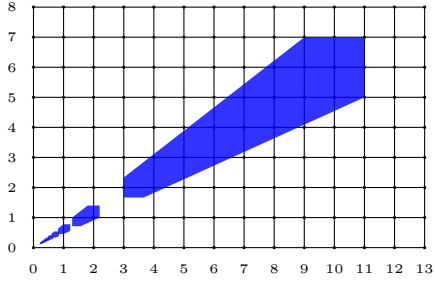
\begin{figure}[H]
\centering
\begin{tikzpicture}[scale=0.4]
\draw (0,0) grid(13,8);
\foreach \x in {0,...,13}
{
	\draw (\x,-0.25) node[below]{\tiny \color{black}{$\x$}};
}
\foreach \y in {0,...,8}
{
	\draw (-0.25,\y) node[left]{\tiny \color{black}{$\y$}};
}
\foreach \x in {0,...,13}
{
	\foreach \y in {0,...,8}
	{
		\filldraw[black] (\x,\y)circle(1pt);
	}
}
\foreach \k in {0,...,10}
{
	\pgfmathsetmacro\a{\k + 0.75};
	\pgfmathsetmacro\b{\k + 0.25};
	\fill[color=blue, opacity=.8] ( 2.25 / \a,1.25/ \a ) -- ( 2.75 / \a,1.25 / \a) -- ( 2.75 / \b,1.25 / \b ) -- ( 2.75 / \b,1.75 / \b ) -- ( 2.25 / \b,1.75 / \b) -- ( 2.25 / \a,1.75 / \a ) -- cycle; 
}
\end{tikzpicture}
\caption{The set $\mathcal{K}_{(2,1)}^2/\mathcal{K}(\delta_2)$ \label{FigLeafDen}}
\end{figure}

In other words, this form of the Lonely Runner Conjecture states that for every integer $d\geq 1$, the set $[1,\infty)^d$ can be covered by the union of the sets $ \mathcal{K}_k^{d}\left( \delta_d \right) / \mathcal{K} \left( \delta_d \right)$ with $k\in \N^d$, each of which is the countable union of convex polytopes (see Figure \ref{FigLeafDen} and Section \ref{SECPrelim}). As we will see, this last formulation of the conjecture is convenient to address the problem of the time required for a runner to get lonely.\\

The Lonely Runner Conjecture is known to be true for $n$ runners with $n\leq 7$ (the integer $n$ refers here to the first statement of the conjecture given in this paper). As seen above, the case  $n=2$ is trivial. The case of three runners, also very simple, was proved by Wills \cite{wills67} in 1967. Several proofs in the context of diophantine approximation, by Betke and Wills \cite{bw72} and Cusick \cite{cusick74}, respectively in 1972 and 1974, established the conjecture for $n=4$. The case $n=5$ was first settled with a computer-assisted proof by Cusick and Pomerance \cite{cp84} in 1984 and later demonstrated with a simpler proof by Biena {\it et al.} \cite{bggst98} in 1998. It is also worth noting that Chen \cite{chen90,chen91a} proved a conjecture which happens to be equivalent to the Lonely Runner Conjecture for the cases $n=3,4,5$. The case $n=6$ was proved with a long and complicated proof by Bohmann, Holzman and Kleitman \cite{bhk01} in 2001, which was later simplified by Renault \cite{renault04} in 2004.  Finally, the conjecture for seven runners was established by Barajas and Serra \cite{bs08} in 2008 and it remains open for all integers $n\geq 8$. Several other problems related to the Lonely Runner Conjecture have also been profusely studied such as the gap of loneliness \cite{ps16,chen94,cc99,tao18} or the validity of the conjecture under various hypotheses on the velocities \cite{pandey09,rtv04,bs09,dubickas11,tao18}. However, to our knowledge the question of the size of the time required for a runner to get lonely has surprisingly not been really adressed before, this is the purpose of the present paper.  \\

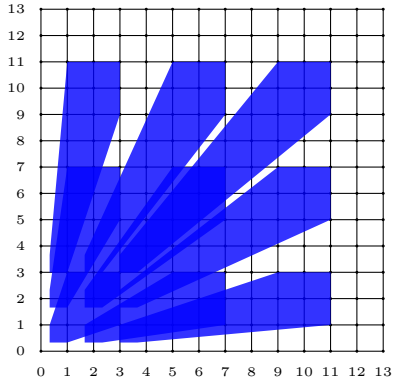
\begin{figure}[H]
\centering
\begin{tikzpicture}[scale=0.35]
\draw (0,0) grid(13,13);
\foreach \x in {0,...,13}
{
	\draw (\x,-0.25) node[below]{\color{black}{\tiny $\x$}};
}
\foreach \y in {0,...,13}
{
	\draw (-0.25,\y) node[left]{\tiny \color{black}{$\y$}};
}
\foreach \x in {0,...,13}
{
	\foreach \y in {0,...,13}
	{
		\filldraw[black] (\x,\y)circle(1pt);
	}
}
\foreach \k in {0,...,2}
{
	\foreach \l in {0,...,2}
	{
		\pgfmathsetmacro\c{4*\k + 3};
		\pgfmathsetmacro\d{4*\l + 3};
		\pgfmathsetmacro\e{4*\k + 1};
		\pgfmathsetmacro\f{4*\l + 1};
		\fill[color=blue, opacity=.8] ( 0.33+ 1.33*\k,0.33+ 1.33*\l) -- (1+ 1.33*\k,0.33+ 1.33*\l) -- (\c,\f) -- (\c,\d) -- (\e,\d) -- ( 0.33+ 1.33*\k,1+ 1.33*\l) -- cycle; 
	}
}
\end{tikzpicture}
\caption{The square  $[1,6]^2$ can be covered by the unions of the sets $\mathcal{K}_{(k,l)}^2 (\delta_2)/\mathcal{K}_0(\delta_2)$ with $(k,l) \in \{0,1,2\}^2$ \label{FigSimu}}
\end{figure}

The starting point of the present paper is the result of computer experiments. As shown in Figure \ref{FigSimu} (for $d=2$ and a not too large parameter $C$), we can check easily for $d=2$ and $d=3$ (the integer $d$ refers here to the Lonely Runner Conjecture in covering form) that large sets of the form $[1,C]^d$ can indeed be covered by the union of the sets $\mathcal{K}_k^d (\delta_d)/\mathcal{K}_0(\delta_d)$ for $k\in \N^d$. This observation suggests that one can maybe expect a strong Lonely Runner Conjecture in covering form where the denominator could be taken to be a finite union of $N$ sets $\mathcal{K}_k(\delta_d)$  of the form 
$$
\mathcal{K}_{\llbracket 0,N-1\rrbracket}(\delta_d) :=\bigcup_{k=0}^{N-1} \mathcal{K}_k(\delta_d),
$$
for some integer $N \in \N^*$ ($\llbracket 0,N-1\rrbracket$ stands for $\N \cap [0,N-1]$). Moreover, as suggested by the title of the paper, this type of result would provide an upper bound on the time required for a runner to get lonely. As a matter of fact, assume that for some integer $d\geq 1$, there exists an integer $N\geq 1$ such that 
$$
[1,\infty)^d \subset \mathcal{K}^{d}\left( \delta_d \right) / \mathcal{K}_{\llbracket 0,N-1\rrbracket} \left( \delta_d \right).
$$
Then, setting $m=d+1$, for any $w_1, \ldots, w_m >0$ with $w_1 \leq w_2, \ldots, w_m$, we have
\[
\left(\frac{w_2}{w_1}, \cdots, \frac{w_m}{w_1}\right) \in [1,\infty)^d \subset  \mathcal{K}^{d}\left( \delta_d \right) / \mathcal{K}_{\llbracket 0,N-1 \rrbracket} \left( \delta_d \right),
\]
so there is $\lambda \in \mathcal{K}_{\llbracket 0,N-1 \rrbracket}(\delta_d)$ such that  
\[
\frac{\lambda}{w_1} w_1 =\lambda   \in \mathcal{K}_{\llbracket 0,N-1 \rrbracket}(\delta_d) \quad \mbox{and} \quad \frac{\lambda}{w_1} \left(w_2, \ldots, w_m \right)  = \lambda \left( \frac{w_2}{w_1}, \ldots, \frac{w_m}{w_1} \right)  \in \mathcal{K}^{d}(\delta_d)
\]
which implies that 
\[
t:=\frac{\lambda}{w_1} \in \left[ 0, N/w_1\right] \mbox{ satisfies } 
t \left( w_1, \ldots, w_m \right) \in \mathcal{K}^{d+1}(\delta_d) =  \mathcal{K}^{m}(1/(m+1)).
\]
On the other hand, if we know that for every $w_1, \ldots, w_m>0$ (with $m=d+1$), there is $t\in [0, N/w_1]$ such that
\[
t \left( w_1, \ldots, w_m \right) \in \mathcal{K}^{m}(1/(m+1)),
\]
then for every tuple $(z_1, \ldots, z_d) \in [1,\infty)^d$, there is $t \in [0, N]$ such that
\[
t \left( 1,z_1, \ldots, z_d \right) \in \mathcal{K}^{m}(1/(m+1)) = \mathcal{K}^{d+1}(\delta_d)
\]
which gives 
\[
\left(z_1,\ldots, z_d\right) \in \mathcal{K}^{d}\left( \delta_d \right) / \mathcal{K}_{\llbracket 0,N-1\rrbracket} \left( \delta_d \right),
\]
because $t\in [0,N]\cap \mathcal{K}(\delta_d)=\mathcal{K}_{\llbracket 0,N-1 \rrbracket}(\delta_d)$. Therefore, we investigate in the present paper the following conjecture:\\

\noindent {\bf Strong Lonely Runner Conjecture.}
{\em For every integer $d\geq 1$, there is an integer $N\geq 1$ such that the following equivalent properties hold:\\

\noindent (Covering from) $[1,\infty)^d \subset \mathcal{K}^{d}\left( \delta_d \right) / \mathcal{K}_{\llbracket 0,N-1 \rrbracket} \left( \delta_d \right)$.\\

\noindent (Obstruction form) For any $w_1, \ldots, w_{d+1}>0$, there exists $t\geq 0$  such that
\[
t \leq \frac{N}{\min\{w_1, \ldots, w_{d+1} \}} \quad \mbox{and} \quad t  \left(w_1, \ldots, w_{d+1}\right) \in \mathcal{K}^{d+1}(\delta_d).
\]
\noindent (Runner form) If $x_1, \ldots, x_{n}:[0,\infty) \rightarrow \R /\Z$ are $n=d+2$ runners with different constant speeds $v_1, \ldots, v_{n} \in \R$ starting from $0$ at time $0$, then for each $i\in \{1, \ldots,n\}$, there exists $t\geq 0$ such that
\[
t \leq \frac{N}{\min \left\{ |v_i -v_j| \, \vert \,  j\in \{1,\ldots, n\}\setminus \{i\} \right\}}
\]
and
\[
 \left|x_j(t)-x_i(t) \right|  \in \left[ \frac{1}{n}, 1-\frac{1}{n}\right] \qquad \forall  j\in \{1,\ldots, n\}\setminus \{i\}.
\]
}\\

It is clear that if the above conjecture is satisfied for some $N \in \N^*$, then there is a minimal integer for which it is satisfied. This remark justifies the following definition: \\

\begin{definition}
For every integer $d\geq 1$, we denote by $N_d$ the least integer $N\geq 1$ such that 
$$
[1,\infty)^d \subset \mathcal{K}^{d}\left( \delta_d \right) / \mathcal{K}_{\llbracket 0,N -1\rrbracket} \left( \delta_d \right),
$$
with the convention that $N_d=\infty$ if the set of $N \in \N^*$ satisfying the above inclusion is empty and 
$$
[1,\infty)^d \subset \mathcal{K}^{d}\left( \delta_d \right) / \mathcal{K} \left( \delta_d \right).
$$
\end{definition}

In other words, for every $d\in \N^*$, $N_d$ is the minimal number of rounds $N$ such that if there are $d+1$ runners with positive speeds and one static runner (so we have $d+2$ runners), starting all from $0$ at time $0$ on the circular track $\R / \Z$, then there is a moment within the first $N$ rounds of the slowest non-static runner such that the static runner is at distance at least $1/(d+2)$ from the others. To summarize the status of the Lonely Runner Conjecture in terms of $N_d$, it is currently known that $N_d=\infty$ for $d=1, 2, 3, 4, 5$ and $N_1=1$ (see {\it e.g.} \cite[Section 6.]{bhs18}). Our first result is the following:\\

\begin{theorem}\label{THM1}
We have $N_1=N_2=N_3=1$. \\
\end{theorem}

In fact, we found out, after the release of a first version of the present paper, that the above result has already been proved, although stated in a different way (see Remark \ref{REMChen}), by Chen in \cite{chen90,chen91a} (see also \cite{chen91b,chen93}). Theorem \ref{THM1} seemingly provides an answer in low dimension to the Question 9 raised by Beck, Hosten and Schymura  in \cite{bhs18}. The proofs for the cases $d=1, 2$ are easy and the one for $d=3$, rather tedious, proceeds in two steps. We first show that the covering property $[1,\infty)^3 \subset \mathcal{K}^{3}\left( \delta_3 \right) / \mathcal{K}_{0} \left( \delta_3 \right)$ is satisfied at infinity, that is, outside of a certain compact set, and then we check by hand that it holds also on that compact set. Let us explain the idea of the asymptotic part. The covering property at infinity is in fact equivalent to showing that for every $(z_1,z_2,z_3) \in [1,\infty)^3$ with $1<z_1<z_2<z_3$ and $z_1$ large enough, we cannot have (see Proposition \ref{PROPcaracx}) 
\[
\left[ \frac{1}{5}, \frac{4}{5} \right] = \mathcal{K}_{0}(\delta_3) \subset  \bigcup_{i=1}^3 \mathcal{B}(\delta_3)/z_i,
\]
where each $\mathcal{B}(\delta_3)/z_i$ for $i=1,2,3$ is defined as 
\[
\mathcal{B}(\delta_3)/z_i := \bigcup_{k\in \N}  \left( \mathcal{B}_k^-/z_i,  \mathcal{B}_k^+/z_i\right) \quad \mbox{with} \quad \mathcal{B}_k^- := \frac{5k-1}{5}, \, \mathcal{B}_k^+ := \frac{5k+1}{5} \, \, \, \forall k\in \N.
\]
 We can show that if the above inclusion is satisfied with $z_1$ large, then $\mathcal{K}_0(\delta_3)$ contains many sets, called $z_1$-kwais, of the form $\mathcal{K}_k/z_1$ (with $k\in \N$) each of which is covered by the union $\mathcal{B}(\delta_3)/z_2 \cup \mathcal{B}(\delta_3)/z_3$. 
 
\vspace{0.2cm}
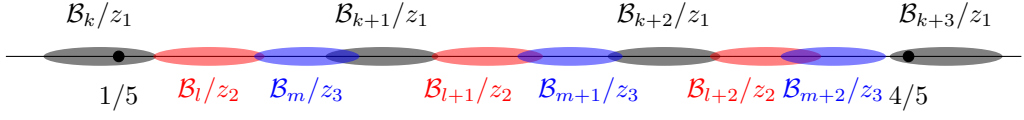
\begin{figure}[H]
\centering
\begin{tikzpicture}
\draw (0,2) -- (13.5,2);
\fill [color=black, opacity=0.5]  (1.25,2) ellipse  (0.75cm and 0.125 cm);
\draw (1.25,2.25) node[above]{$\mathcal{B}_{k}/z_1$};
\fill [color=black, opacity=0.5]  (5,2) ellipse  (0.75cm and 0.125 cm);
\draw (5,2.25) node[above]{$\mathcal{B}_{k +1}/z_1$};
\fill [color=black, opacity=0.5]  (8.75,2) ellipse  (0.75cm and 0.125 cm);
\draw (8.75,2.25) node[above]{$\mathcal{B}_{k+2}/z_1$};
\fill [color=black, opacity=0.5]  (12.5,2) ellipse  (0.75cm and 0.125 cm);
\draw (12.5,2.25) node[above]{$\mathcal{B}_{k+3}/z_1$};
\filldraw[black] (1.5,2)circle(2pt);
\draw (1.5,1.75) node[below]{\color{black}{$1/5$}};
\filldraw[black] (12,2)circle(2pt);
\draw (12,1.75) node[below]{\color{black}{$4/5$}};
\fill [color=red, opacity=0.5] (2.7,2) ellipse (0.74cm and 0.125 cm);
\draw (2.7,1.25) node[above]{\color{red}{$\mathcal{B}_{l}/z_2$}};
\fill [color=red, opacity=0.5] (6.4,2) ellipse (0.74cm and 0.125 cm);
\draw (6.15,1.25) node[above]{\color{red}{$\mathcal{B}_{l+1}/z_2$}};
\fill [color=red, opacity=0.5] (10.1,2) ellipse (0.74cm and 0.125 cm);
\draw (9.65,1.25) node[above]{\color{red}{$\mathcal{B}_{l+2}/z_2$}};
\fill [color=blue, opacity=0.5] (4,2) ellipse (0.7cm and 0.125 cm);
\draw (4,1.25) node[above]{\color{blue}{$\mathcal{B}_{m}/z_3$}};
\fill [color=blue, opacity=0.5] (7.5,2) ellipse (0.7cm and 0.125 cm);
\draw (7.75,1.25) node[above]{\color{blue}{$\mathcal{B}_{m+1}/z_3$}};
\fill [color=blue, opacity=0.5] (11,2) ellipse (0.7cm and 0.125 cm);
\draw (11,1.25) node[above]{\color{blue}{$\mathcal{B}_{m+2}/z_3$}};
\end{tikzpicture}
\caption{We have $( \mathcal{B}_{l+j+1}^-/z_2 -  \mathcal{B}_{k+j+1}^-/z_1)- ( \mathcal{B}_{l+j}^-/z_2 -  \mathcal{B}_{k+j}^-/z_1) = 1/z_2-1/z_1<0 $ and $( \mathcal{B}_{m+j+1}^-/z_3 -  \mathcal{B}_{k+j+1}^-/z_1)- ( \mathcal{B}_{m+j}^-/z_3 -  \mathcal{B}_{k+j}^-/z_1) = 1/z_3-1/z_1<0 $ for $j=1,2$, which means that  $\mathcal{B}_{l+j}^-/z_2$ and  $\mathcal{B}_{m+j}^-/z_3$ are shifted to the left (with respect to $ \mathcal{B}_{k+j}^-/z_1$) when $j$ increases\label{FIGk}}
\end{figure}
\vspace{0.2cm}

Then, we can show that each $z_1$-kwai has to be covered by $\mathcal{B}(\delta_3)/z_2 \cup \mathcal{B}(\delta_3)/z_3$ in the same manner (see Lemma \ref{LEM3Subcoverings}) and obtain a contradiction, provided that $z_1$ is sufficiently large, from the fact that, since $z_3>z_2>z_1$, the sets $\mathcal{B}_{l+j}/z_2$ and $\mathcal{B}_{m+j}/z_3$ are shifted to the left along consecutive $z_1$-kwais (see Figure \ref{FIGk}). In fact, as we shall see, the argument of the first step that we have just described amounts to checking that some convex polytopes have empty interior and the second step consists in verifying that a given family of convex polytopes covers a compact subset of $\R^d$, two properties that can be checked automatically by a computer.  Since the combinatorics of the $d=4$ case are not too involved, our method also allows us to give a computer-assisted proof of the following result: \\

\begin{theorem}\label{THM2}
We have $N_4=2$. \\
\end{theorem}

This theorem says that the inclusion $[1,\infty)^4 \subset \mathcal{K}^{4}\left( \delta_4 \right) / \mathcal{K}_{\llbracket 0,1\rrbracket} \left( \delta_4 \right)$ holds true while the inclusion $[1,\infty)^4 \subset \mathcal{K}^{4}\left( \delta_4 \right) / \mathcal{K}_{0} \left( \delta_4 \right)$ does not. The first inclusion follows principally from computer checks while the second one is a consequence of the fact that (see Section \ref{SEC4Dnot})
\[
 \left(\frac{13(1-\epsilon_1)}{5},  \frac{3211(1-\epsilon_2)}{935}, \frac{247(1-\epsilon_3)}{55}, \frac{61009(1-\epsilon_4)}{10285} \right) \notin \mathcal{K}^{4}\left( \delta_4 \right) / \mathcal{K}_{0}(\delta_4),
\]
for every $\epsilon_1,\epsilon_2,\epsilon_3,\epsilon_4>0$ satisfying $\epsilon_1<\epsilon_3<\epsilon_2<\epsilon_4<32/3211$. If we come back to the initial statement of the Lonely Runner Conjecture, the latter property implies for example that if we consider six runners $x_1, \ldots, x_6:[0,\infty) \rightarrow \R /\Z$ starting from $0$ at time $0$ with speeds $v_1, \ldots, v_6$ given by
\begin{multline*}
v_1=0, \quad v_2=5\cdot 11^2 \cdot 17 =10285, \quad v_3= 11^2 \cdot 13 \cdot 17-1=26740,\\
 v_4 = 11 \cdot 13^2  \cdot 19-2= 35319, \quad v_5 = 11 \cdot 13 \cdot 17 \cdot 19-2=46187, \quad v_6=  13^2 \cdot 19^2-4=61005,
\end{multline*}
 then the static runner $x_1$ is never separated by a distance at least $1/6$ from the others during the first round of $x_2$  but it is the case at least once during $x_2$'s second round, for instance at $t=5/(4v_2)$.  Our method does not extend easily to higher dimension, so we do not know for example if $N_5<\infty$ ($N_5=\infty$ has been proved by Barajas and Serra \cite{bs08}). We can just mention that we checked by computer that $N_5\geq 2$.\\

Our proof of Theorem \ref{THM2} is based on the notion of weak chains that allows us to deal with coverings by sets of the form 
\[
\bigcup_{i=1}^4 \left( h_i + \mathcal{B}(\delta_4)/z_i\right) \quad \mbox{with} \quad \mathcal{B}(\delta_4)/z_i := \bigcup_{k\in \N}  \left( \bigl(k-\delta_4\bigr)/z_i, \bigl(k+\delta_4\bigr)/z_i\right),
\]
where $h_1,\ldots, h_4$ belong to $\R$. This approach happens to be well adapted to treat the following form of the Lonely Runner Conjecture where the runners have are not supposed to start necessarily from the origin ($\sharp A$ stands for the cardinal of a set $A$):\\

\noindent {\bf Lonely Runner Conjecture with Free Starting Positions.}
{\em For every integer $n\geq 2$, if $x_1, \ldots, x_n:[0,\infty) \rightarrow \R /\Z$ are $n$ runners with constant speeds $v_1, \ldots, v_n \in \R$, that is such that $x_i(t)=\left\{ \bar{x}_i+tv_i\right\}$ for all $t\geq 0$ for some tuple $(\bar{x}_1,\ldots,\bar{x}_n) \in (\R /\Z)^n$, then the following property holds for each runner $x_i$ with $i\in \{1, \ldots,n\}$: If 
\begin{eqnarray}\label{CardSpeeds}
\sharp \Bigl\{ |v_j-v_i| \, \vert \, j\in\bigl\{1,\ldots,n\bigr\} \setminus \bigl\{i\bigr\}\Bigr\} = n-1,
\end{eqnarray}
 then there is  $t> 0$ such that 
$$
\left| x_j(t)- x_i(t) \right|  \in \left[ \frac{1}{n}, 1-\frac{1}{n}\right] \qquad \forall j\in \bigl\{1,\ldots,n\bigr\} \setminus \bigl\{i\bigr\}.
$$}\\

Note that the assumption (\ref{CardSpeeds}) is important. As a matter of fact, for $n=3$, we can consider for every parameter $v\in \R$, three runners $x_1, x_2, x_3:[0,\infty) \rightarrow \R /\Z$ given by 
\[
x_1(t)=v\, t, \quad x_2(t) = \left\{\frac{1}{4}+ (v+1) \,t\right\}, \quad x_3(t) = \left\{\frac{1}{4}+(v-1) \, t\right\} \qquad \forall t \geq 0.
\]
Then the runner $x_1$ will be lonely at time $t>0$ if and only if there are two integers $k,l$ such that $1/4+t\in [k+1/3,k+2/3]$ and $1/4-t\in [l+1/3,l+2/3]$ which is equivalent to
\[
k+\frac{1}{12} \leq t \leq k+\frac{5}{12} \quad \mbox{and} \quad (-l-1) + \frac{7}{12} \leq t \leq (-l-1) + \frac{11}{12},
\]
which has no solution. The Lonely Runner Conjecture with Free Starting Positions is equivalent to the following conjecture (see Proposition \ref{PROPLRCFSPequivLRCFSPOF} and compare with Proposition \ref{PROPLRCequivLRCOF}):\\

\noindent {\bf Lonely Runner Conjecture with Free Starting Positions - Obstruction Form.}
{\em Let $m \geq 1$ be an integer, then  any affine line $\mathcal{L} \subset \R^m$ having a direction vector with distinct positive coordinates intersects infinitely many hypercubes of the form 
\[
\mathcal{K}_k^m(1/(m+1)) := k+ \left[ \frac{1}{m+1},\frac{m}{m+1}\right]^m,
\] 
with $k\in \N^m$.}

\begin{figure}[H]
\centering
\begin{tikzpicture}[scale=0.35]
\draw (0,0) grid(11,11);
\foreach \x in {0,...,10}
{
	\foreach \y in {0,...,10}
	{
		\fill[color=red] (\x+0.33,\y+0.33) -- (\x+0.66,\y + 0.33) -- (\x+0.66,\y+0.66) -- (\x+0.33,\y+0.66) -- cycle;  
	}
}
\draw [color=blue] (0,0) -- (8.5,11);
\draw [color=blue] (3.27,0) -- (5.54,11);
\draw [color=blue] (0,3.1) -- (11,7.3);
\end{tikzpicture}
\caption{In dimension $2$, the Lonely Runner Conjecture with Free Starting Points in Obstruction Form asserts that every blue line, which is not parallel to the horizontal axis, the vertical axis and the diagonal $y=x$, intersects infinitely many red squares $\mathcal{K}^2_k(1/3)$ (with $k\in \N^2$) \label{FigLightFree}}
\end{figure}
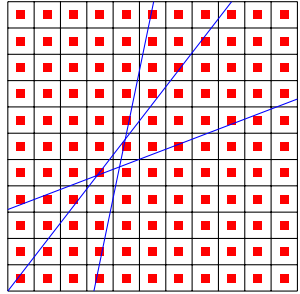

Note that the condition on the coordinates of the direction vector of the affine line is the counterpart of (\ref{CardSpeeds}) in the previous statement of the conjecture. The method used to prove Theorem \ref{THM2} allows us to show the following result:

\begin{theorem}\label{THM3}
The Lonely Runner Conjecture with Free Starting Positions  is satisfied for $n=2,3,4$.
\end{theorem}

Note that the above result for $2$ runners is obvious and the proof for $n=3$ appeared in \cite{bhs18}. The proof of Theorem \ref{THM3} is given in Section \ref{SECTHM3}.\\

Our last result is concerned with the gap of loneliless that we mentioned above. According to \cite{ps16}, for every integer $m\geq 2$, the gap of loneliness $D_m$ is defined as the supremum  of $\delta \in (0,1/2]$ such that for every $w_1, \ldots, w_m>0$ there exists $t>0$ such that $t(w_1,\ldots,w_m)\in \mathcal{K}^m(\delta)$. Following the discussion we had before the statement of the Strong Lonely Runner Conjecture, this definition is equivalent to saying that 
\[
D_m := \sup \Bigl\{ \delta \in (0,1/2] \, \vert \,  [1,\infty)^{m-1} \subset \mathcal{K}^{m-1}\left( \delta \right) / \mathcal{K} (\delta ) \Bigr\}.
\]
 Let us introduce the following definition: \\

\begin{definition}
For every integer $d \geq 1$ and every $N\in \N^*\cup \{\infty\}$, we denote by $\Delta_d^N$, called gap of loneliness in $N$ rounds, the supremum of all real numbers $\delta \in (0,1/2]$ such that 
\[
[1,\infty)^d \subset  \mathcal{K}^{d}\left( \delta \right) / \mathcal{K}_{\llbracket 0,N -1\rrbracket} \left( \delta \right).
\]
\end{definition}
\bigskip

We note that if $N$ in the above definition is finite then the supremum is actually a maximum and we check easily that we have for every integer $d\geq 1$,
\[
\Delta^N_{d+1} \leq \Delta^N_{d} \quad \forall N\in \N^*\cup \{\infty\} \quad \mbox{and}  \quad \Delta^{N}_d \leq \Delta^{N+1}_d \leq \Delta_d^{\infty}=D_{d+1} \quad \forall N\in \N^*.
\]
Our Theorem \ref{THM1} yields $\Delta_d^1 =\delta_d=1/(d+2)$ for $d=1,2,3$, Theorem \ref{THM2}  implies $\Delta_4^1 <1/6=\delta_4$ and $\Delta_4^2=\delta_4$, and the result by Barajas and Serra \cite{bs08} gives $D_6=\Delta_5^{\infty} =\delta_5$. Moreover, as it is well-known the Dirichlet approximation Theorem \cite{cassels72} implies that 
\[
D_{d+1}=\Delta^{\infty}_{d} \leq  \frac{1}{d+2} \qquad \forall d \in \N^*.
\]
 In the spirit of \cite{chen94,cc99,ps16}, we are able to give a lower bound for $\Delta_d^1$ of the form, $1/(2(d+1))+c/(d+1)^2$ with $c>0$, for all $d\in \N^*$. For every $\delta \in (0,1/2]$ and every $z\geq 1$, we call $\mathcal{K}_0(\delta)$-measure of $z$ the nonnegative real number $\mathcal{M}_0^{\delta}(z)$ defined by ($\mathcal{L}^1$ denotes the Lebesgue measure in $\R$)
\[
\mathcal{M}_0^{\delta}(z) := \mathcal{L}^1 \left( \mathcal{B}(\delta)/z \cap \mathcal{K}_0(\delta) \right).
\]
We can check easily (see Section \ref{SECK0Measure}) that 
\[
\mathcal{M}_0^{\delta}(z) \leq \frac{2\delta(1-2\delta)}{1-\delta} < 2\delta \qquad \forall z \geq 1, \, \forall \delta \in (0,1/3].
\]
This result implies that if for some $d\in \N^*$ and $\delta \in (0,1/3]$, we have 
\[
[1,\infty)^d \subsetneq  \mathcal{K}^{d}\left( \delta \right) / \mathcal{K}_{0} \left( \delta \right),
\]
then any $(z_1,\ldots,z_d) \in [1,\infty)^d$ (see Proposition \ref{PROPcaracx}) such that 
\[
\mathcal{K}_{0}(\delta) \subset \left( \cup_{i=1}^d \mathcal{B}(\delta)/z_i\right)
\]
satisfies 
\begin{eqnarray*}
1-2 \delta = \mathcal{L}^1 \left( \mathcal{K}_0(\delta)\right) & = & \mathcal{L}^1 \left(  \cup_{i=1}^d (\mathcal{B}(\delta)/z_i \cap \mathcal{K}_0(\delta)) \right) \\
& \leq & \sum_{i=1}^d  \mathcal{M}_0^{\delta}(z_i)\leq \frac{2\delta(1-2\delta)d}{1-\delta},
\end{eqnarray*}
which implies that $\delta\geq 1/(2d+1)$. In conclusion, the gap of loneliness in one round satisfies the following inequalities: 

\begin{proposition}\label{PROPGap}
For every integer $d\geq 1$, we have 
\[
 \frac{1}{2d+1} \leq \Delta_d^1 \leq \frac{1}{d+2}.
\]
\end{proposition}

This result implies that for every integer $m\geq 2$, we have 
\[
D_{m} = \Delta_{m-1}^{\infty} \geq \frac{1}{2m-1}.
\]
This bound is very far from the lower bound obtained by Tao in \cite{tao18}. In fact, the reader familiar with the Lonely Runner Conjecture has probably observed that we do not restrict our attention to integral speeds. As a matter of fact, it can be shown that the conjecture holds true if and only if it is satisfied for rational speeds (see \cite{bhk01}). This fact allows Tao to work with integral speeds and to show, thanks to a fine study of the size of unions of the so-called Bohr sets, that there is an absolute constant $c>0$ (see \cite[Theorem 1.2]{tao18}) such that
\[
D_m \geq \frac{1}{2m} + \frac{c\log m}{m^2(\log(\log m))^2},
\]
for all sufficiently large $m$. We do not know if $ \Delta_d^1$ may verify a lower bound of the form $1/(2(d+1))+1/((d+1)^2o(1))$ for large integers $d$. Moreover, noting that  the Strong Lonely Runner Conjecture is verified for some pair $d,N\in \N^*$ if and only if  
\[
[1,\infty)^d \cap \Q^d \subset \mathcal{K}^{d}\left( \delta_d \right) / \mathcal{K}_{\llbracket 0,N-1 \rrbracket} \left( \delta_d \right),
\]
we don't know either if Tao's method can be applied to improve the results presented in this paper. \\

The paper is organized as follows: Several definitions and a preliminary result are introduced in Section \ref{SECPrelim}, the proofs of Theorems \ref{THM1}, \ref{THM2} and \ref{THM3} are given respectively in Sections \ref{SECProofTHM1}, \ref{SECProofTHM2} and \ref{SECTHM3}, and Section \ref{SECK0Measure} deals with the upper bound on $\mathcal{K}_0(\delta)$-measures.

\section{Preliminary results}\label{SECPrelim}

We introduce several definitions and state in Proposition \ref{PROPcaracx} the characterization of the covering property that will be our main tool in the proofs of Theorems \ref{THM1} and \ref{THM2}.

\subsection{Feathers, bridges, kwais and beams}\label{SECNotions}

Given $d\in \N^*$, $N\in \N^* \cup\{\infty\}$, and $\delta \in (0,1/2)$, we introduce several definitions. 

\begin{definition}
$\quad$
\begin{itemize}
\item[(i)] We call $(d,\delta)$-feather, or simply feather, any set of the form 
$$
F_{k,l}^d(\delta) :=  \mathcal{K}_k^d(\delta)/\mathcal{K}_l(\delta),
$$
with $k=(k_1,\ldots, k_d)  \in \N^d$ and $l\in \N$.
\item[(ii)] For every real number $z\geq 1$, we call $z$-bridge any set of the form 
$$
\mathcal{B}_k(\delta)/z := \left( \mathcal{B}_k^-(\delta)/z,  \mathcal{B}_k^+(\delta)/z  \right) = \left(\frac{k-\delta}{z},\frac{k+\delta}{z}\right),
$$
with $k\in \N$ and we denote its length (which does not depend on $k$) by $\b^{\delta}(z) = 2\delta/z$. 
\item[(iii)] For every real number $z\geq 1$, we call $z$-kwai any set of the form 
$$
\mathcal{K}_k(\delta)/z := \left[ \mathcal{K}_k^-(\delta)/z,  \mathcal{K}_k^+(\delta)/z  \right] = \left[\frac{k+\delta}{z},\frac{k+1-\delta}{z}\right],
$$
with $k\in \N$ and we denote its length (which does not depend on $k$) by $\k^{\delta}(z) = (1-2\delta)/z$.
\item[(iv)] For every $(z_1, \ldots, z_d) \in [1,\infty)^d$, we call $z$-beam the set 
\[
\mathcal{P}^d_{\delta , N} (z_1, \ldots, z_d) := \bigcup_{l \in \llbracket 0,N-1\rrbracket} P^d_{\delta,l} (z_1, \ldots, z_d)  \subset \R^d,
\]
where each $P^d_{\delta,l} (z_1, \ldots, z_d)$  (with $l \in \llbracket 0,N-1\rrbracket$) is the set of $(\lambda_1, \ldots, \lambda_d)\in \R^d$ verifying
  \[
(l+\delta)z_i+\delta -1 \leq  \lambda_i \leq (l+1-\delta) z_i - \delta \qquad \forall i=1, \ldots,d
\]
and
\[
\lambda_j z_i -  \lambda_i z_j \geq (\delta-1)z_i + \delta z_j \qquad \forall i,j=1, \ldots,d \mbox{ with } i\neq j.
\]
\end{itemize}
\end{definition}

We recall that a convex polytope is a convex compact set with a finite number of extreme points, or equivalently a compact set given by the intersection of a finite number of half-spaces. We refer for example the reader to Ziegler's monograph \cite{ziegler95} for more details on convex polytopes. 

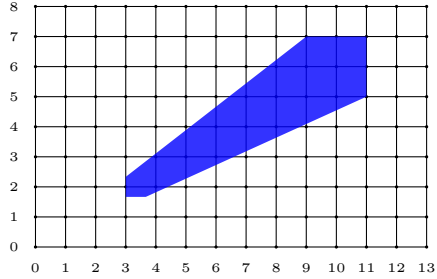
\begin{figure}[H]
\centering
\begin{tikzpicture}[scale=0.4]
\draw (0,0) grid(13,8);
\foreach \x in {0,...,13}
{
	\draw (\x,-0.25) node[below]{\tiny \color{black}{$\x$}};
}
\foreach \y in {0,...,8}
{
	\draw (-0.25,\y) node[left]{\tiny \color{black}{$\y$}};
}
\foreach \x in {0,...,13}
{
	\foreach \y in {0,...,8}
	{
		\filldraw[black] (\x,\y)circle(1pt);
	}
}
\foreach \k in {0,...,0}
{
	\pgfmathsetmacro\a{\k + 0.75};
	\pgfmathsetmacro\b{\k + 0.25};
	\fill[color=blue, opacity=.8] ( 2.25 / \a,1.25/ \a ) -- ( 2.75 / \a,1.25 / \a) -- ( 2.75 / \b,1.25 / \b ) -- ( 2.75 / \b,1.75 / \b ) -- ( 2.25 / \b,1.75 / \b) -- ( 2.25 / \a,1.75 / \a ) -- cycle; 
}
\end{tikzpicture}
\caption{The $(2,\delta_2)$-feather $\mathcal{K}_{(2,1)}^2(\delta_2)/\mathcal{K}_0(\delta_2)$ is the intersection of the rectangle $[3,11] \times [5/3,7]$ with the positive cone bounded by the two lines of equations $7x_1-9x_2=0$ and $11x_2-5x_1=0$\label{FigLeaf}}
\end{figure}

The proof of the following result is postponed to Appendix \ref{APPPROPFeather}. 

\begin{proposition}\label{PROPFeather} The following properties hold:
\begin{itemize}
\item[(i)]  For every $k=(k_1,\ldots, k_d)  \in \N^d$ and every $l\in \N$, $F_{k,l}^d(\delta)$ is the convex polytope equal to the intersection of the hyperrectangle $R_{k,l}^d(\delta)$ corresponding to the set of  $(z_1,\ldots,z_d) \in \R^d$ satisfying 
\begin{eqnarray}\label{28septeq1}
 \frac{k_i+\delta}{l+1-\delta} \leq z_i \leq \frac{k_i+1-\delta}{l +\delta} \qquad \forall i=1, \ldots,d
\end{eqnarray}
with the convex (positive) cone $C_{k}^d(\delta)$  given by 
$$
C_k^d (\delta):=(0,\infty) \cdot  \mathcal{K}_k^d(\delta), 
$$ 
which correspond to the set of  $(z_1,\ldots,z_d) \in (0,\infty)^d$ satisfying
\begin{eqnarray}\label{28septeq2}
\frac{z_i}{k_i+\delta}  - \frac{z_j}{k_j+1-\delta} \geq 0 \qquad \forall i,j =1, \ldots,d \mbox{ with } i\neq j,
\end{eqnarray}
and in addition  the lower face of of $F_{k,l}^d(\delta)$ defined as
\[
 \partial^- F_{k,l}^d(\delta) =  \Bigr\{ z \in F_{k,l}^d(\delta) \, \vert \, tz\notin F_{k,l}^d(\delta) \, \forall t \in [0,1)\Bigr\} 
\]
verifies
\[
 \partial^- F_{k,l}^d(\delta) = F_{k,l}^d(\delta) \, \cap \, \bigcup_{i=1}^d \left\{ z \in \R^d \, \vert \, z_i = \frac{k_i+\delta_d}{l+1-\delta_d}\right\}.
\]
\item[(ii)] For every $z\in [1,\infty)$ and every $k\in \N$, we have  
\[
\mathcal{K}_0(\delta) \setminus (\mathcal{B}_k (\delta)/z) \neq \emptyset.
\] 
\item[(iii)] If $\delta \leq 1/4$, then for every $z\in [(2-\delta)/(1-\delta),\infty)$, the set $\mathcal{K}_0(\delta)\setminus (\mathcal{B}(\delta)/z)$ contains a $z$-kwai, that is, a closed interval of length $(1-2\delta)/z$, and for every $z\in [1,(2-\delta)/(1-\delta))$ the set $\mathcal{K}_0(\delta)\setminus (\mathcal{B}(\delta)/z)$ contains at least one closed interval of length $\geq (1-\delta)/z-\delta$.
\item[(iv)] For every $(z_1, \ldots, z_d) \in [1,\infty)^d$ and every $l \in \llbracket 0,N-1\rrbracket$, $P^d_{\delta,l} (z_1, \ldots, z_d)$ is a convex polytope.
\end{itemize}
\end{proposition}

We now state the result whose assertion (ii) will be the main tool in the proofs of Theorems \ref{THM1} and \ref{THM2}.

\begin{proposition}\label{PROPcaracx}
For every $(z_1, \ldots, z_d) \in [1,\infty)^d$, the following properties are equivalent:
 \begin{itemize}
 \item[(i)] $(z_1, \ldots, z_d) \in \mathcal{K}^{d}(\delta) / \mathcal{K}_{\llbracket 0,N-1 \rrbracket}(\delta) = \cup_{k\in \N^d}  \cup_{l\in \llbracket 0,N-1 \rrbracket} F_{k,l}^d(\delta)$.
 \item[(ii)] $\mathcal{K}_{\llbracket 0,N-1\rrbracket}(\delta) \cap \left( \cap_{i=1}^d \mathcal{K}(\delta)/z_i \right) = \mathcal{K}_{\llbracket 0,N-1\rrbracket}(\delta) \setminus \left( \cup_{i=1}^d \mathcal{B}(\delta)/z_i\right)   \neq \emptyset$.
 \item[(iii)] $\mathcal{P}^d_{\delta,N} (z_1, \ldots, z_d) \cap \N^d \neq \emptyset$.
  \end{itemize}
\end{proposition}

\begin{proof}[Proof of Proposition \ref{PROPcaracx}]
Let  $(z_1,\ldots, z_d)\in [1,\infty)^d$ be fixed. If $(z_1,\ldots, z_d) \in \mathcal{K}^{d}(\delta) / \mathcal{K}_{\llbracket 0,N-1 \rrbracket}(\delta)$, then there are an integer $l\in [0,N-1]$ and $\lambda \in \mathcal{K}_l(\delta)$ such that $\lambda z_i \in \mathcal{K}(\delta)$ for all $i=1,\ldots, d$, which implies that
\[
\lambda \in \mathcal{K}_l(\delta) \cap \left(\cap_{i=1}^d \mathcal{K}(\delta)/z_i \right) = \mathcal{K}_{l}(\delta) \setminus \left( \cup_{i=1}^d \mathcal{B}(\delta)/z_i\right).
\]
So we have (i)$\Rightarrow$(ii). If the set $\mathcal{K}_{\llbracket 0,N-1\rrbracket}(\delta) \cap \left( \cap_{i=1}^d \mathcal{K}(\delta)/z_i \right)$ is not empty, then there are $\lambda \in \R$, $l\in  \llbracket 0,N-1\rrbracket$ and $k_1, \ldots, k_d\in \N$ such that
\[
\lambda \in \mathcal{K}_{l}(\delta) \quad \mbox{and} \quad \lambda \in \mathcal{K}_{k_i}(\delta)/z_i \quad \forall i=1, \ldots,d, 
\]
which, by setting $k:=(k_1,\ldots, k_d)$, implies that
\[
(z_1, \ldots, z_d) \in \mathcal{K}_{k}(\delta) /\mathcal{K}_l(\delta) = F_{k,l}^d(\delta).
\]
By Proposition \ref{PROPFeather} (i), we infer that the inequalities (\ref{28septeq1})-(\ref{28septeq2}) are satisfied, then we have 
 \begin{eqnarray}\label{17sept3}
(l+\delta)z_i+\delta -1 \leq  k_i \leq (l+1-\delta) z_i - \delta \qquad \forall i=1, \ldots,d
\end{eqnarray}
and
\begin{eqnarray}\label{17sept4}
k_j z_i -  z_i x_j \geq (\delta-1)z_i + \delta z_j \qquad \forall i,j=1, \ldots,d \mbox{ with } i\neq j,
\end{eqnarray}
which shows that $\mathcal{P}^d_{\delta,N} (z_1, \ldots, z_d) \cap \N^d \neq \emptyset$. Thus, we have (ii)$\Rightarrow$(iii). Finally to prove (iii)$\Rightarrow$(i), we note that if (iii) holds then there are $k=(k_1,\ldots,k_d)\in \N^d$ and $l\in \llbracket 0,N-1\rrbracket$ such that the above inequalities (\ref{17sept3})-(\ref{17sept4}) hold. Then the equivalent inequalities (\ref{28septeq1})-(\ref{28septeq2}) hold, which means that $(z_1,\ldots, z_d)$ belongs to the feather $F_{k,l}^d(\delta)$ and in other words that it is in the set $\mathcal{K}^{d}(\delta) / \mathcal{K}_{\llbracket 0,N-1 \rrbracket}(\delta)$.
 \end{proof}

\begin{remark}\label{REMChen}
Proposition \ref{PROPFeather} (i) together with Proposition \ref{PROPcaracx} (i) show that the Lonely Runner Conjecture in covering form holds true for an integer $d\geq 1$ if and only if for every $(z_1, \ldots, z_d) \in [1,\infty)^d$, there are $(k_1,\ldots, k_d)  \in \N^d$ and $l\in \N$ such that 
\[
 \frac{k_i+\delta_d}{l+1-\delta_d} \leq z_i \leq \frac{k_i+1-\delta_d}{l +\delta_d} \qquad \forall i=1, \ldots,d
\]
and
\[
\frac{z_i}{k_i+\delta_d}  - \frac{z_j}{k_j+1-\delta_d} \geq 0 \qquad \forall i,j =1, \ldots,d \mbox{ with } i\neq j,
\]
where we recall that $\delta_d=1/(d+2)$. Moreover, we can observe that if $(z_1, \ldots, z_d) \in [1,\infty)^d$, $(k_1,\ldots, k_d)  \in \N^d$ and $l\in \N$ are given then setting $n=d+1$, $a_1=1$, $l_1=l$ and $a_i = z_{i-1}, l_i=k_{i-1}$ for all $i=2,\ldots,n$ allows to write the above inequalities in the equivalent form
\[
a_il_j-a_jl_i \leq \frac{n}{n+1}\, a_j - \frac{1}{n+1} \, a_i \qquad \forall i,j = 1, \ldots, n,
\]
which correspond exactly to the set of inequalities of a conjecture stated by Chen in \cite{chen90}. Chen has proven his conjecture for $n=2,3, 4$, that is $d=1,2,3$, in \cite{chen91a} and that in this case $l_1$ can be taken to be equal to $0$ see \cite[Remark 3 p. 197]{chen91a} which is exactly the statement of Theorem \ref{THM1}. 
\end{remark}

\begin{remark}
Our approach is reminiscent to those  Beck, Hosten and Schymura \cite{bhs18} (see for example Proposition 1 p. 3) who work in dimension $d+1$ and with integral speeds (see also \cite{hm17}).
\end{remark}

\begin{remark}
Characterization (iii) allows to see the conjecture as a problem of geometry of numbers {\it à la} Minkowski (see \cite{cassels72}), it won't be used in this paper. 
\end{remark} 

\subsection{$\mathcal{K}_0$-coverings, chains and $1$-subchains}\label{SECCovering}

Let $d\in \N^*$ and $N\in \N^*\cup \{\infty\}$ be fixed. Proposition \ref{PROPcaracx} shows that a way to prove that the covering property 
$$
[1,\infty)^d \subset \mathcal{K}^{d}\left( \delta_d \right) / \mathcal{K}_{\llbracket 0,N-1\rrbracket} \left( \delta_d \right)
$$
is satisfied is to demonstrate that for every $(z_1,\ldots, z_d)\in [1,\infty)^d$, the set  
 \[
 \mathcal{K}_{\llbracket 0,N-1\rrbracket}(\delta_d) \setminus \left( \cup_{i=1}^d \mathcal{B}(\delta_d)/z_i\right) = \mathcal{K}_{\llbracket 0,N-1\rrbracket}(\delta_d) \cap \left( \cap_{i=1}^d \mathcal{K}(\delta_d)/z_i \right) 
 \]
 is not empty. Our proof of the covering property at infinity will follow from the study of that set, so from now on, in order to simplify the notations, we write $\mathcal{B}$, $\mathcal{K}$, $\b(z),\k(z)$ instead of $\mathcal{B}(\delta_d), \mathcal{K}(\delta_d), \b^{\delta_d}(z),\k^{\delta_d}(z)$, and if $(z_1,\ldots, z_d)\in [1,\infty)^d$ is fixed then we set $\b_i:=\b^{\delta_d}(z_i)$ and $\k_i:=\k^{\delta_d}(z_i)$.

\begin{definition}
We say that a tuple $(z_1, \ldots, z_d) \in [1,\infty)^d$ is strictly well-ordered if  
\[
1< z_1< z_2 \cdots < z_d
\]
and we say that it is a $\mathcal{K}_0$-covering,  if  it is strictly well-ordered and satisfies
$$
\mathcal{K}_0 \subset \bigcup_{i=1}^d \mathcal{B}/z_i.
$$
\end{definition}

\begin{figure}[H]
\centering
\begin{tikzpicture}[scale=0.95]
\draw (-1,2) -- (14,2);
\fill [color=black, opacity=0.5]  (0.64,2) ellipse  (1cm and 0.125 cm);
\draw (0.64,2.25) node[above]{$\mathcal{B}_{k_1}/z_1$};
\fill [color=black, opacity=0.5]  (6.64,2) ellipse  (1cm and 0.125 cm);
\draw (6.64,2.25) node[above]{$\mathcal{B}_{k_5}/z_1$};
\fill [color=black, opacity=0.5]  (12.64,2) ellipse  (1cm and 0.125 cm);
\draw (13,2.25) node[above]{$\mathcal{B}_{k_{10}}/z_1$};
\filldraw[black] (1.5,2)circle(2pt);
\draw (1.5,1.75) node[below]{\color{black}{$1/6$}};
\filldraw[black] (12,2)circle(2pt);
\draw (12,1.75) node[below]{\color{black}{$5/6$}};
\fill [color=red, opacity=0.5] (2.55,2) ellipse (0.95cm and 0.125 cm);
\draw (2.55,1.25) node[above]{\color{red}{$\mathcal{B}_{k_2}/z_2$}};
\fill [color=red, opacity=0.5] (8.25,2) ellipse (0.95cm and 0.125 cm);
\draw (8.25,1.25) node[above]{\color{red}{$\mathcal{B}_{k_6}/z_2$}};
\fill [color=blue, opacity=0.5] (4.71,2) ellipse (0.94cm and 0.125 cm);
\draw (4.71,1.25) node[above]{\color{blue}{$\mathcal{B}_{k_4}/z_3$}};
\fill [color=blue, opacity=0.5] (10.35,2) ellipse (0.94cm and 0.125 cm);
\draw (10.35,1.25) node[above]{\color{blue}{$\mathcal{B}_{k_8}/z_3$}};
\fill [color=brown, opacity=0.5] (3.78,2) ellipse (0.44cm and 0.125 cm);
\draw (3.78,2.25) node[above]{\color{brown}{$\mathcal{B}_{k_3}/z_4$}};
\fill [color=brown, opacity=0.5] (9.06,2) ellipse (0.44cm and 0.125 cm);
\draw (9.06,2.25) node[above]{\color{brown}{$\mathcal{B}_{k_7}/z_4$}};
\fill [color=brown, opacity=0.5] (11.7,2) ellipse (0.44cm and 0.125 cm);
\draw (11.7,2.25) node[above]{\color{brown}{$\mathcal{B}_{k_9}/z_4$}};
\end{tikzpicture}
\caption{A chain $(1243124341\vert k_1k_2k_3k_4k_5k_6k_7k_8k_9k_{10})$ associated with a $\mathcal{K}_0$-covering $(z_1, z_2, z_3,z_4) \in [1,\infty)^4$ \label{FIGcovering}}
\end{figure}
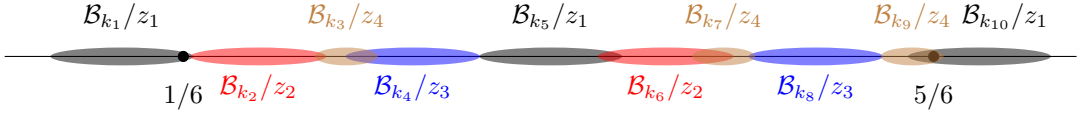

By Proposition \ref{PROPFeather} (ii), the set $\mathcal{K}_0$ cannot be covered by only one of the $\mathcal{B}/z_i$ (with $i\in \{1,\ldots,d\}$), this observation justifies the following definition (see Figure \ref{FIGcovering}):  

\begin{definition}
We call chain associated with a $\mathcal{K}_0$-covering $(z_1, \ldots, z_d) \in [1,\infty)^d$, that we denote by $(i_1\cdots i_{\ell}\vert k_1\cdots k_{\ell})$, any pair of ordered families of indices $i_1,\ldots,i_{\ell} \in \{1,\ldots,d\}$ and of positive integers $k_1,\ldots, k_{\ell}$ with $\ell$ an integer $\geq 2$  satisfying  
\begin{eqnarray}\label{Hanoi1}
\delta_d \in \mathcal{B}_{k_1}/z_{i_1}, \quad 1-\delta_d \in \mathcal{B}_{k_{\ell}}/z_{i_{\ell}} \quad
\mbox{and} \quad  \mathcal{B}_{k_j}^+/z_{i_j} \in \mathcal{B}_{k_{j+1}}/z_{i_{j+1}} \quad \forall j=1, \ldots, \ell-1.
\end{eqnarray}
\begin{eqnarray}\label{Hanoi1bis}
\delta_d <  \mathcal{B}_{k_{2}}^-/z_{i_{2}}, \quad  \mathcal{B}_{k_{\ell-1}}^+/z_{i_{\ell-1}} < 1-\delta_d
\end{eqnarray}
and
\begin{eqnarray}\label{Hanoi2}
 \mathcal{B}_{k_{j}}^+/z_{i_{j}} <  \mathcal{B}_{k_{j+2}}^-/z_{i_{j+2}}  \quad \forall j=1, \ldots, \ell-2,
\end{eqnarray}
provided  $\ell \geq 3$.
\end{definition}

\vspace{0.2cm}
\vspace{0.2cm}

Note that by definition a chain $(i_1\cdots i_{\ell}\vert k_1\cdots k_{\ell})$ is minimal in the sense  that $\mathcal{K}_0$ is not covered if we remove one of the $\mathcal{B}_{k_j}/z_{i_j}$, that is,
\[
\mathcal{K}_0 \subsetneq \bigcup_{j\in \{1,\ldots,\ell\} \setminus \{\bar{j}\}} \mathcal{B}_{k_j}/z_{i_j} \qquad \forall \bar{j} \in \{1,\ldots,\ell\}.
\]
The proofs of Theorems \ref{THM1} and \ref{THM2} are based on the following result. We denote by $\mathcal{S}^d$ the closed convex set of tuples $(z_1,\ldots,z_d)$ such that $1\leq z_1\leq \ldots \leq z_d$.

\begin{proposition}\label{PROP3oct1}
If the covering property $\mathcal{O} \cap \mathcal{S}^d \subset \mathcal{K}^{d} / \mathcal{K}_{0}$ does not hold for an open set $\mathcal{O} \subset (0,\infty)^d$, then there exist a $\mathcal{K}_0$-covering $(z_1, \ldots, z_d) \in \mathcal{O}$ along with a chain $(i_1\cdots i_{\ell}\vert k_1\cdots k_{\ell})$.
\end{proposition}

\begin{proof}[Proof of Proposition \ref{PROP3oct1}]
Let $\mathcal{O} \subset (0,\infty)^d$ be an open set such that 
\[
\mathcal{O} \cap \mathcal{S}^d \subsetneq \mathcal{K}^{d} / \mathcal{K}_{0}.
\]
The set $(\mathcal{O} \cap \mathcal{S}^d) \setminus (\mathcal{K}^{d} / \mathcal{K}_{0})$ is a non-empty open set of $\mathcal{S}^d$, so there is a tuple $(z_1,\ldots, z_d) \in \mathcal{O}$ such that $z_1<\ldots <z_d$ and $(z_1,\ldots,z_d) \notin \mathcal{K}^{d} / \mathcal{K}_{0}$. By Proposition \ref{PROPcaracx}, we infer that 
\[
\mathcal{K}_0 \subset \bigcup_{i=1}^d \mathcal{B}(\delta)/z_i,
\]
so that $(z_1,\ldots, z_d)$ is a $\mathcal{K}_0$-covering. Then, there is a finite set $\mathcal{F}\subset \N^* \times \{1,\ldots,d\}$ such that
\[
\mathcal{K}_0 \subset \bigcup_{(k,i)\in \mathcal{F}} \mathcal{B}(\delta)/z_i
\]
and
\[
\mathcal{K}_0 \subsetneq \bigcup_{(k,i) \in \{1,\ldots,\ell\} \setminus \{(\bar{k},\bar{i})\}} \mathcal{B}_{k}/z_{i} \qquad \forall (\bar{k},\bar{i}) \in \mathcal{F}.
\]
We get a chain by ordering the set of $\mathcal{B}_k^-/z_i$ for $(k,i)\in \mathcal{F}$.
\end{proof}

The following result provides a few obstructions on chains associated with $\mathcal{K}_0$-coverings.

\begin{proposition}\label{LEMTHMmain2}
Let $(i_1\cdots i_{\ell}\vert k_1\cdots k_{\ell})$  be a chain associated with a   $\mathcal{K}_0$-covering  $(z_1,\ldots, z_d)$, then the following properties hold:
\begin{itemize}
\item[(i)] For every $j=1, \ldots, \ell-1$, $i_j\neq i_{j+1}$.
\item[(ii)] If $d\geq 2$, then $l\geq3$ and if $i_j=i_{j+2}$ for some $j\in \{1, \ldots, \ell -2\}$ then $i_{j+1}<i_j$.
\end{itemize}
\end{proposition}

\begin{proof}[Proof of Proposition \ref{LEMTHMmain2}]
Let us consider $(i_1\cdots i_{\ell}\vert k_1\cdots k_{\ell})$ a chain associated with a   $\mathcal{K}_0$-covering  $(z_1,\ldots, z_d)$. To prove (i) suppose for contradiction that there is  some $j\in \{1,\ldots, \ell-1\}$  such that $i_j= i_{j+1}$. If $\ell=2$, then we have $\mathcal{K}_0 \subset  \mathcal{B}/z_{i_j}$ which is prohibited by Proposition \ref{PROPFeather} (ii). Otherwise, we have $\ell\geq 3$ and (\ref{Hanoi1})-(\ref{Hanoi2}). Therefore, we have  $\mathcal{B}_{k_{j+1}}^-/z_{i_{j}} = \mathcal{B}_{k_{j+1}}^-/z_{i_{j+1}} <   \mathcal{B}_{k_{j}}^+/z_{i_{j}}$ and  $\mathcal{B}_{k_{j}}^+/z_{i_{j}} <   \mathcal{B}_{k_{j+1}}^+/z_{i_{j+1}}= \mathcal{B}_{k_{j+1}}^+/z_{i_{j}}$ which imply respectively $k_{j+1}\leq k_j$ and $k_j<k_{j+1}$, a contradiction. To prove (ii), assume that $i_j=i_{j+2}$ for some $j\in \{1,\ldots, \ell-2\}$ (with $\ell\geq 3$), then by (\ref{Hanoi1})-(\ref{Hanoi2}) we have 
$$
 \mathcal{B}_{k_{j+1}}^-/z_{i_{j+1}} <   \mathcal{B}_{k_{j}}^+/z_{i_{j}} <  \mathcal{B}_{k_{j+2}}^-/z_{i_{j}} <  \mathcal{B}_{k_{j+1}}^+/z_{k_{j+1}},
$$
which means that the length of $\mathcal{B}_{k_{j+1}}/z_{i_{j+1}}$, which is equal to $\b_{{i_{j+1}}}=2\delta_d/z_{i_{j+1}}$, is strictly larger than $\k_{i_j}:=(1-2\delta_d)/z_{i_j}$. If $d\geq2$, this implies that $z_{i_{j+1}}<z_{i_j}$, that is, $i_{j+1}<i_j$ (because  $(z_1, \ldots, z_d)$ is strictly well-ordered).
\end{proof}

In the sequel, we often denote chains by $(i_1\cdots i_{\ell})$ instead of $(i_1\cdots i_{\ell}\vert k_1\cdots k_{\ell})$. We also use the word subchain, that we denote by $[i_j\cdots i_{j+r}]$ to avoid any confusion, to  speak of consecutive subsequences of chains. Moreover, we make use of the following definition.
 
\begin{definition}\label{DEF1subchain}
Given a $\mathcal{K}_0$-covering $(z_1, \ldots, z_d)$ associated with a chain $(i_1\cdots i_{\ell}\vert k_1\cdots k_{\ell})$, we call $1$-subchain of length $L\in \N^*$ any sequence of indices of the form 
\[
[i_{1}'\cdots i_{r}'\vert k_{1}'\cdots k_{r}'],
\]
with $r$ an integer $\geq 3$, satisfying the following properties:
\begin{itemize}
\item[(i)] $i_1'=i_{r}'=1$ and $k_r'-k_{1}'=L$.
\item[(ii)] $\mathcal{B}_{k_1'}^+/z_1, \mathcal{B}_{k_{r}'}^-/z_1 \in \mathcal{K}_0$.
\item[(iii)] $\mathcal{B}_{k_1'}^+/z_1\in \mathcal{B}_{k_{2}'}/z_{i_{2}'}$ and $\mathcal{B}_{k_{r}'}^-/z_1\in \mathcal{B}_{k_{r-1}'}/z_{i_{r-1}'}$.
\item[(iv)] For every $k\in \N\cap (k_1,k_{r})$, there is $s\in \{2,r-1\}$ such that $(i_{s},k_{s})=(1,k)$. 
\item[(v)] $[i_{2}'\ldots i_{r-1}'\vert k_{2}'\ldots k_{r-1}']$ is a subchain of $(i_1\cdots i_{\ell}\vert k_1\cdots k_{\ell})$.
\end{itemize}
Moreover, we say that a $\mathcal{K}_0$-covering $(z_1, \ldots, z_d)$  is $1$-nice if for every pair $\bar{k},\hat{k}\in \N^*$ such that $\hat{k}>\bar{k}$ and $\mathcal{B}_{\bar{k}}^+/z_1, \mathcal{B}_{\hat{k}}^-/z_1 \in \mathcal{K}_0$, any chain $(i_1\cdots i_{\ell}\vert k_1\cdots k_{\ell})$ associated with $(z_1, \ldots, z_d)$ admits a subchain $(i_j\cdots i_{j'}\vert k_j\cdots k_{j'})$ with $j<j'$ in $\{1,\ldots, \ell\}$ such that  
\[
[1 i_j\cdots i_{j'}1\vert \bar{k}k_{j}\cdots k_{j'}\hat{k}]
\]
is a $1$-subchain of length $\hat{k}-\bar{k}$.
\end{definition}

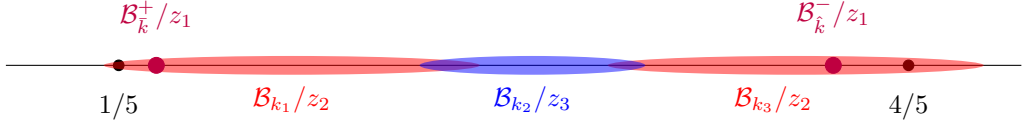
\begin{figure}[H]
\centering
\begin{tikzpicture}
\draw (0,2) -- (13.5,2);
\filldraw[black] (1.5,2)circle(2pt);
\draw (1.5,1.75) node[below]{\color{black}{$1/5$}};
\filldraw[black] (12,2)circle(2pt);
\draw (12,1.75) node[below]{\color{black}{$4/5$}};
\fill [color=red, opacity=0.5] (3.8,2) ellipse (2.5cm and 0.125 cm);
\draw (3.8,1.25) node[above]{\color{red}{$\mathcal{B}_{k_1}/z_2$}};
\fill [color=red, opacity=0.5] (10.5,2) ellipse (2.5cm and 0.125 cm);
\draw (10.2,1.25) node[above]{\color{red}{$\mathcal{B}_{k_3}/z_2$}};
\fill [color=blue, opacity=0.5] (7,2) ellipse (1.5cm and 0.125 cm);
\draw (7,1.25) node[above]{\color{blue}{$\mathcal{B}_{k_2}/z_3$}};
\filldraw[purple] (2,2)circle(3pt);
\draw (2,2.3) node[above]{\color{purple}{$\mathcal{B}_{\bar{k}}^+/z_1$}};
\filldraw[purple] (11,2)circle(3pt);
\draw (11,2.3) node[above]{\color{purple}{$\mathcal{B}_{\hat{k}}^-/z_1$}};
\end{tikzpicture}
\caption{The points $\mathcal{B}_{\bar{k}}^+/z_1$ and $\mathcal{B}_{\hat{k}}^-/z_1$ belong to $\mathcal{K}_0$ but do not appear in the chain $(k_1k_2k_3\vert232)$ associated with the $\mathcal{K}_0$-covering $(z_1,z_2,z_3)$\label{fignice}}
\end{figure}

In fact, the above definition is justified by the fact that we may have a $\mathcal{K}_0$-covering $(z_1, \ldots, z_d)$ associated with a chain $(i_1\cdots i_{\ell}\vert k_1\cdots k_{\ell})$ such that we may have a pair $\bar{k},\hat{k}\in \N^*$ satisfying $\hat{k}>\bar{k}$ and $\mathcal{B}_{\bar{k}}^+/z_1, \mathcal{B}_{\hat{k}}^-/z_1 \in \mathcal{K}_0$, but where no $j\in \{1,\ldots, \ell\}$ verifies $(i_j,k_j)=(1,\bar{k})$ or $(i_j,k_j)=(1,\hat{k})$ (see Figure \ref{fignice}).
 
\vspace{0.2cm}
\vspace{0.2cm}

\begin{remark}\label{REM2}
In other words, a $\mathcal{K}_0$-covering $(z_1, \ldots, z_d)$  is $1$-nice if when we consider $\bar{k},\hat{k}\in \N^*$ satisfying $\hat{k}>\bar{k}$ and $\mathcal{B}_{\bar{k}}^+/z_1, \mathcal{B}_{\hat{k}}^-/z_1 \in \mathcal{K}_0$, all instances of $\mathcal{B}_k/z_1$ for $k\in \N \cap (\bar{k},\hat{k})$ appears in any chain associated with $(z_1, \ldots, z_d)$. Therefore, to show that a $\mathcal{K}_0$-covering $(z_1, \ldots, z_d)$ is $1$-nice it is sufficient to prove that an interval of the form 
\begin{eqnarray}\label{SetREM2}
I = \left[ \mathcal{K}_k^-/z_1, \mathcal{K}^+_{k+1}/z_1\right] \subset \mathcal{K}_0 \quad \mbox{of length} \quad 2\k_1 + \b_1 
\end{eqnarray}
can not be covered by the union of $\mathcal{B}_i/z_i$ for $i=2,\ldots,d$. In particular, to his aim we can apply the properties given in Proposition \ref{LEMTHMmain2}.
\end{remark}

\section{Proof of Theorem \ref{THM1}}\label{SECProofTHM1}

For each case $d=1,2,3$, we suppose for contradiction that $[1,\infty)^d \subsetneq \mathcal{K}^{d} / \mathcal{K}_{0}$ and consider, thanks to Proposition \ref{PROP3oct1}, a $\mathcal{K}_0$-covering $(z_1, \ldots, z_d)$ associated with a chain $(i_1\cdots i_{\ell}\vert k_1\cdots k_{\ell})$ with $\ell\geq 2$. 

\subsection{Case $d=1$} 

We have $\delta_1=1/3$ and $\mathcal{K}_0=[1/3,2/3]$. If $\mathcal{K}_0\subset \mathcal{B}/z$ for some $z>1$, then by Proposition \ref{LEMTHMmain2} (i) there is $k\in \N$ such that $\mathcal{K}_0\subset \mathcal{B}_k/z$, which contradicts Proposition \ref{PROPFeather} (ii). 

\subsection{Case $d=2$} 

We have $\delta_2=1/4$, $\mathcal{K}_0=[1/4,3/4]$ and $\b(z)=\k(z)=1/(2z)$ for all $z\geq 1$. We know by Proposition \ref{LEMTHMmain2} (ii) that $\ell \geq 3$. If $\ell =3$,  then the chain has to be $[212\vert k_j k_{j+1} k_{j+2}]$ (by Proposition \ref{LEMTHMmain2} (ii)). But we observe that if $1/4$ belongs to $\mathcal{B}_k/z$ for some $z> 1$ and $k\in \N^*$, then we have $ \mathcal{B}_1^-/z=3/(4z)\leq  \mathcal{B}_k^-/z< 1/4$ which implies that  $z>3$ and 
\[
\mathcal{B}_{k+1}^+/z= \mathcal{B}_{k}^-/z + \b(z)+\frac{1}{z} < \frac{1}{4}+ \frac{1}{2z} + \frac{1}{z} = \frac{1}{4} + \frac{3}{2z} < \frac{1}{4} + \frac{1}{2} = \frac{3}{4}.
\]
This shows that if $1/4$ belongs to $\mathcal{B}_{k_j}/z_2$ then $\mathcal{B}^+_{k_j+1}/z_2\in \mathcal{K}_0$ and as a consequence $k_{j+2}-k_j\geq 2$. Then, we have 
\[
\frac{1}{2z_1} = \b_1 > \mathcal{B}^-_{k_j+2}/z_2 -  \mathcal{B}^+_{k_j}/z_2 \geq \frac{2}{z_2} - \b_2 = \frac{3}{2z_2},
\] 
which yields $z_2>3z_1$ and 
\[
\frac{3}{4z_1} = \mathcal{B}_{1}^-/z_1\leq \mathcal{B}_{k_2}^-/z_1 < \frac{1}{4} + \b_2 < \frac{1}{4} +\frac{1}{6z_1}.
\]
This means that $z_2>3z_1>7$ which gives $2\b_2+\b_1 <1/7+3/14 = 5/14<1/2$, a contradiction (to the fact that $[212\vert k_j k_{j+1} k_{j+2}]$ is a chain). In conclusion, the chain $(i_{j}\cdots i_{j+r}\vert k_{j}\cdots k_{j+r})$ associated with the $\mathcal{K}_0$-covering $(z_1,z_2)$ has length $\geq 4$ and so  has the form $(i_1i_2i_1i_2\cdots)$, but this contradicts Proposition  \ref{LEMTHMmain2} (ii). Therefore there is no $\mathcal{K}_0$-covering.

\subsection{Case $d=3$}

\begin{figure}[H]
\centering
\begin{tikzpicture}
\draw (0,2) -- (13.5,2);
\fill [color=black, opacity=0.5]  (2,2) ellipse  (1cm and 0.125 cm);
\draw (2,2.25) node[above]{$\mathcal{B}_{k_j}/z_1$};
\fill [color=black, opacity=0.5]  (7,2) ellipse  (1cm and 0.125 cm);
\draw (7,2.25) node[above]{$\mathcal{B}_{k_{j+3}}/z_1$};
\fill [color=black, opacity=0.5]  (12,2) ellipse  (1cm and 0.125 cm);
\draw (12,2.25) node[above]{$\mathcal{B}_{k_{j+6}}/z_1$};
\fill [color=red, opacity=0.5] (3.75,2) ellipse (0.95cm and 0.15 cm);
\draw (3.75,1.25) node[above]{\color{red}{$\mathcal{B}_{k_{j+1}}/z_2$}};
\fill [color=blue, opacity=0.5] (5.3,2) ellipse (0.9cm and 0.15 cm);
\draw (5.3,1.225) node[above]{\color{blue}{$\mathcal{B}_{k_{j+2}}/z_3$}};
\fill [color=red, opacity=0.5] (8.8,2) ellipse (0.95cm and 0.15 cm);
\draw (8.8,1.25) node[above]{\color{red}{$\mathcal{B}_{k_{j+4}}/z_2$}};
\fill [color=blue, opacity=0.5] (10.3,2) ellipse (0.9cm and 0.15 cm);
\draw (11.5,1.225) node[above]{\color{blue}{$\mathcal{B}_{k_{j+5}}/z_3$}};
\end{tikzpicture}
\caption{A $1$-subchain of the form $[1231231\vert k_j\cdots k_{j+6}]$ with $k_{j+6}=k_{j+3}+1=k_j+2$ and $k_{j+4}=k_{j+1}+1$\label{subchain1231231}}
\end{figure}
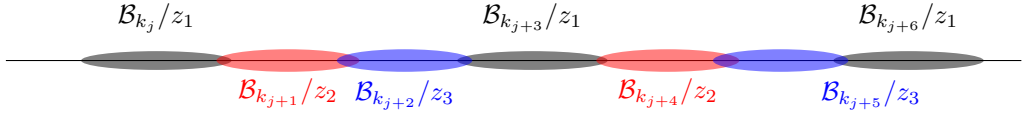

\begin{figure}[H]
\centering
\begin{tikzpicture}
\draw (0,2) -- (13.5,2);
\fill [color=black, opacity=0.5]  (2,2) ellipse  (1cm and 0.125 cm);
\draw (2,2.25) node[above]{$\mathcal{B}_{k_j}/z_1$};
\fill [color=black, opacity=0.5]  (7,2) ellipse  (1cm and 0.125 cm);
\draw (7,2.25) node[above]{$\mathcal{B}_{k_{j+3}}/z_1$};
\fill [color=black, opacity=0.5]  (12,2) ellipse  (1cm and 0.125 cm);
\draw (12,2.25) node[above]{$\mathcal{B}_{k_{j+6}}/z_1$};
\fill [color=blue, opacity=0.5] (3.75,2) ellipse (0.95cm and 0.15 cm);
\draw (3.75,1.25) node[above]{\color{blue}{$\mathcal{B}_{k_{j+1}}/z_3$}};
\fill [color=red, opacity=0.5] (5.3,2) ellipse (0.9cm and 0.15 cm);
\draw (5.3,1.225) node[above]{\color{red}{$\mathcal{B}_{k_{j+2}}/z_2$}};
\fill [color=blue, opacity=0.5] (8.8,2) ellipse (0.95cm and 0.15 cm);
\draw (8.8,1.25) node[above]{\color{blue}{$\mathcal{B}_{k_{j+4}}/z_3$}};
\fill [color=red, opacity=0.5] (10.3,2) ellipse (0.9cm and 0.15 cm);
\draw (11.5,1.225) node[above]{\color{red}{$\mathcal{B}_{k_{j+5}}/z_2$}};
\end{tikzpicture}
\caption{A $1$-subchain of the form $[1321321\vert k_j\cdots k_{j+6}]$ with $k_{j+6}=k_{j+3}+1=k_j+2$ and $k_{j+4}=k_{j+1}+1$\label{subchain1321321}}
\end{figure}
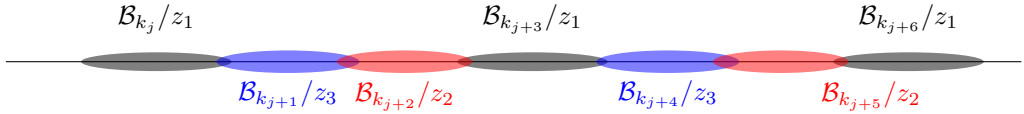

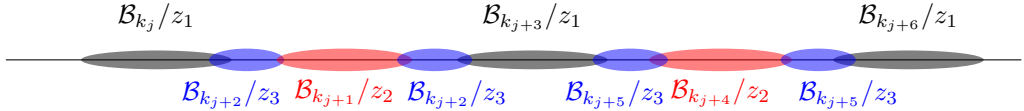
\begin{figure}[H]
\centering
\begin{tikzpicture}
\draw (0,2) -- (13.5,2);
\fill [color=black, opacity=0.5]  (2,2) ellipse  (1cm and 0.125 cm);
\draw (2,2.25) node[above]{$\mathcal{B}_{k_j}/z_1$};
\fill [color=black, opacity=0.5]  (7,2) ellipse  (1cm and 0.125 cm);
\draw (7,2.25) node[above]{$\mathcal{B}_{k_{j+3}}/z_1$};
\fill [color=black, opacity=0.5]  (12,2) ellipse  (1cm and 0.125 cm);
\draw (12,2.25) node[above]{$\mathcal{B}_{k_{j+6}}/z_1$};
\fill [color=blue, opacity=0.5] (3.2,2) ellipse (0.5cm and 0.15 cm);
\draw (3,1.225) node[above]{\color{blue}{$\mathcal{B}_{k_{j+2}}/z_3$}};
\fill [color=red, opacity=0.5] (4.5,2) ellipse (0.9cm and 0.15 cm);
\draw (4.5,1.25) node[above]{\color{red}{$\mathcal{B}_{k_{j+1}}/z_2$}};
\fill [color=blue, opacity=0.5] (5.7,2) ellipse (0.5cm and 0.15 cm);
\draw (6,1.225) node[above]{\color{blue}{$\mathcal{B}_{k_{j+2}}/z_3$}};
\fill [color=red, opacity=0.5] (9.5,2) ellipse (0.95cm and 0.15 cm);
\draw (9.5,1.25) node[above]{\color{red}{$\mathcal{B}_{k_{j+4}}/z_2$}};
\fill [color=blue, opacity=0.5] (8.3,2) ellipse (0.5cm and 0.15 cm);
\draw (8.1,1.225) node[above]{\color{blue}{$\mathcal{B}_{k_{j+5}}/z_3$}};
\fill [color=blue, opacity=0.5] (10.8,2) ellipse (0.5cm and 0.15 cm);
\draw (11.2,1.225) node[above]{\color{blue}{$\mathcal{B}_{k_{j+5}}/z_3$}};
\end{tikzpicture}
\caption{A $1$-subchain of the form $[132313231\vert k_j\cdots k_{j+8}]$ with $k_{j+8}=k_{j+4}+1=k_j+2$, $k_{j+3}=k_{j+1}+1$ and $k_{j+6}=k_{j+2}+1$\label{subchain132313231}}
\end{figure}

We have $\delta_3=1/5$, $\mathcal{K}_0=[1/5,4/5]$ and $\b(z)=2/(5z), \k(z)=3/(5z)$ for all $z\geq 1$. Recall that, arguing by contradiction, we consider a $\mathcal{K}_0$-covering $(z_1, z_2, z_3)$ associated with a chain $(i_1\cdots i_{\ell}\vert k_1\cdots k_{\ell})$ with $\ell\geq 2$. We note that, by Remark \ref{REM2} and Proposition \ref{LEMTHMmain2}, $(z_1, z_2, z_3)$ is $1$-nice because if a set $I$ of the form (\ref{SetREM2}) is covered by the union of $2$-bridges and $3$-bridges then the associated subchain has to be $[23]$, $[32]$ or $[323]$, but those subchains have lengths given respectively by $\b_2+\b_3$, $\b_2+\b_3$ and $2\b_3+\b_2<14/(15z_2)$ (because  $[323]$ yields $\b_2>\k_3 \Leftrightarrow z_3>3z_2/2$) and  all those quantities are strictly less than the length of $I$ given by $2\k_1 + \b_1 =8/(5z_1)$. As a consequence, all instances of $1$-bridges have to appear in any subchain connecting two $1$-bridges. Moreover, by Proposition \ref{LEMTHMmain2}, any $1$-subchain of length $1$ has to be in the set $\{[1231], [1321], [13231]\}$. The following result describes precisely $1$-subchains of length $2$ (see Figures \ref{subchain1231231},  \ref{subchain1321321},  \ref{subchain132313231}):

\begin{lemma}\label{LEM3Subcoverings}
Any $1$-subchain of length $2$ has one of the following forms: 
\begin{itemize}
\item[-] $[1231231\vert k_j\cdots k_{j+6}]$ with $k_{j+6}=k_{j+3}+1=k_j+2$ and $k_{j+4}=k_{j+1}+1$,
\item[-] $[1321321\vert k_j\cdots k_{j+6}]$ with $k_{j+6}=k_{j+3}+1=k_j+2$ and $k_{j+4}=k_{j+1}+1$,
\item[-] $[132313231\vert k_j\cdots k_{j+8}]$ with $k_{j+8}=k_{j+4}+1=k_j+2$, $k_{j+3}=k_{j+1}+1$ and $k_{j+6}=k_{j+2}+1$.
\end{itemize}
\end{lemma}

\vspace{0.2cm}
\vspace{0.2cm}

\vspace{0.2cm}
\vspace{0.2cm}

\vspace{0.2cm}
\vspace{0.2cm}

\begin{proof}[Proof of Lemma \ref{LEM3Subcoverings}]
Any $1$-subchain of length 2 has to be in the set
\begin{multline*}
\Bigl\{ [1231231], [1231321], [12313231], [1321231], [1321321], \\
[13213231], [13231231], [13231321], [132313231] \Bigr\}.
\end{multline*}
We first prove that the subchains $[1231321]$, $[12313231]$, $[1321231]$, $[13213231]$, $[13231231]$, $[13231321]$ cannot occur, for each  case  we assume that it is the case and infer a contradiction:\\

\noindent - $[1231321\vert k_j \cdots k_{j+6}]$: We have $k_{j+6}=k_{j+3}+1=k_j+2$. The subchain $[313\vert k_{j+2}k_{j+3}k_{j+4}]$ gives $z_3>3z_1/2$ (because $2/(5z_1)=\b_1>\k_3=3/(5z_3)$). Moreover, the inequalities 
\[
\mathcal{B}_{k_{j+1}}^-/z_2 < \mathcal{B}_{k_j}^+/z_1<\mathcal{B}_{k_{j+3}}^+/z_1< \mathcal{B}_{k_{j+5}}^-/z_2
\]
 yield  $k_{j+5} \geq k_{j+1}+2$, because otherwise we have $k_{j+5}=k_{j+1}+1$ which implies $1/z_2=\mathcal{B}_{k_{j+5}}^-/z_2- \mathcal{B}_{k_{j+1}}^-/z_2>\mathcal{B}_{k_{j+3}}^+/z_1-\mathcal{B}_{k_j}^+/z_1=1/z_1$, a contradiction. Thus, we have 
 \[
 \frac{8}{5z_2}= \frac{2}{z_2} -\b_2\leq \mathcal{B}_{k_{j+5}}^-/z_2- \mathcal{B}_{k_{j+1}}^+/z_2 < \b_1+2\b_3 < \b_1+\frac{8}{15z_1}= \frac{14}{15z_1}, 
 \]
 which gives $z_2>12z_1/7$,  then $\b_2+\b_3< 7/(30z_1)+ 4/(15z_1) = 1/(2z_1)<\k_1$, a contradiction. \\

\noindent - $[12313231\vert k_j \cdots k_{j+7}]$:  The proof is the same as for the previous case, it it left to the reader.\\

\noindent - $[1321231\vert k_j \cdots k_{j+6}]$: The subchain $[212]$ gives  $z_3>z_2>3z_1/2$, which implies $\b_2+\b_3 < 8/(15z_1)<3/(5z_1)=\k_1$, a contradiction. \\

\noindent - $[13213231\vert k_j \cdots k_{j+7}]$: The subchain $[323]$ gives $z_3>3z_2/2>3z_1/2$. Let us show that $k_{j+4}\geq k_{j+1}+2$. Otherwise,  we have $k_{j+4}= k_{j+1}+1$ and we know that $ \mathcal{B}_{k_{j+1}}^-/z_3 < \mathcal{B}_{k_{j}}^+/z_1 <  \mathcal{B}_{k_{j+3}}^+/z_1<  \mathcal{B}_{k_{j+4}}^+/z_3$. Then, we have
\[
\frac{1}{z_1}= \mathcal{B}_{k_{j+3}}^+/z_1 - \mathcal{B}_{k_{j}}^+/z_1 < \mathcal{B}_{k_{j+4}}^+/z_3 - \mathcal{B}_{k_{j+1}}^-/z_3 = \frac{1}{z_3} + \b_3 =  \frac{7}{5z_3},
\]
which contradicts $z_3>3z_1/2$. Thus  the subchain $[3213\vert k_{j+1}\cdots k_{j+4} ]$ gives
\[
\frac{4}{5z_1} =2\b_1 > \b_1+\b_2 > \mathcal{B}_{k_{j+4}}^-/z_3 - \mathcal{B}_{k_{j+1}}^+/z_3 \geq \frac{2}{z_3} - \b_3  = \frac{8}{5z_3},
\]
which implies $z_3>2z_1$ and then $\b_2+\b_3 < \b_1 + \b_1/2< \k_1$, a contradiction.\\

\noindent - $[13231231\vert k_j \cdots k_{j+7}]$. The subchain $[323]$ gives $z_3>3z_2/2>3z_1/2$. Let us show that $k_{j+6}\geq k_{j+3}+2$. Otherwise,  we have $k_{j+6}= k_{j+3}+1$ and we know that $ \mathcal{B}_{k_{j+3}}^-/z_3 < \mathcal{B}_{k_{j+4}}^-/z_1 <  \mathcal{B}_{k_{j+7}}^-/z_1<  \mathcal{B}_{k_{j+6}}^+/z_3$. Then, we have
\[
 \frac{1}{z_1} = \mathcal{B}_{k_{j+7}}^-/z_1 - \mathcal{B}_{k_{j+4}}^-/z_1 < \mathcal{B}_{k_{j+6}}^+/z_3 - \mathcal{B}_{k_{j+3}}^-/z_3 = \frac{1}{z_3} + \b_3 =  \frac{7}{5z_3},
\]
which contradicts $z_3>3z_1/2$. Thus  the subchain $[3123\vert k_{j+3}\cdots k_{j+6} ]$ gives
\[
\frac{4}{5z_1} =2\b_1 > \b_1+\b_2 > \mathcal{B}_{k_{j+6}}^-/z_3 - \mathcal{B}_{k_{j+3}}^+/z_3 \geq \frac{2}{z_3} - \b_3  = \frac{8}{5z_3},
\]
which implies $z_3>2z_1$ and then $\b_2+\b_3 < \b_1 + \b_1/2< \k_1$, a contradiction.\\

\noindent - $[13231321\vert k_j \cdots k_{j+7}]$: The subchain $[323]$ gives $z_3>3z_2/2>3z_1/2$. Moreover, the inequalities 
\[
\mathcal{B}_{k_{j+2}}^+/z_2 < \mathcal{B}_{k_{j+4}}^-/z_1<\mathcal{B}_{k_{j+7}}^-/z_1< \mathcal{B}_{k_{j+6}}^+/z_2
\]
 yield  $k_{j+6} \geq k_{j+2}+2$, because otherwise we have $k_{j+6}=k_{j+2}+1$ which implies $1/z_2=\mathcal{B}_{k_{j+6}}^+/z_2- \mathcal{B}_{k_{j+2}}^+/z_2>\mathcal{B}_{k_{j+7}}^-/z_1-\mathcal{B}_{k_{j+4}}^-/z_1=1/z_1$, a contradiction. Then, we have 
 \[
 \frac{8}{5z_2}= \frac{2}{z_2} -\b_2\leq \mathcal{B}_{k_{j+6}}^-/z_2- \mathcal{B}_{k_{j+2}}^+/z_2 < \b_1+2\b_3 < \frac{2}{5z_1}+\frac{8}{15z_2}, 
 \]
 which implies $z_2>8z_1/3$ and $\b_2+\b_3 < 3/(20z_1)+4/(15z_1)=5/(12z_1)<\k_1$, a contradiction. \\

We have shown that the only admissible $1$-subchains of length $2$ are $[1231231\vert k_j \cdots k_{j+6}]$, $[1321321\vert k_j \cdots k_{j+6}]$ and  $[132313231\vert k_j \cdots k_{j+8}]$,  it remains to check what can be $k_j \ldots, k_{j+r}$ ($r=6$ or $8$) in each case: \\

\noindent - $[1231231\vert k_j \cdots k_{j+6}]$:  If $k_{j+4} \geq k_{j+1}+2$, then $2 \b_1 > \b_3 + \b_1 > 2/z_2- \b_2$, which gives $z_2>2z_1$ and so $\b_3+ \b_2 <2\b_2 <  2/(5z_1)<\k_1$, a contradiction.\\

\noindent - $[1321321\vert k_j \cdots k_{j+6}]$: If $k_{j+4} \geq k_{j+1}+2$, then $2 \b_1 > \b_2 + \b_1 > 2/z_3- \b_3$, which gives $z_3>2z_1$ and so $\b_3+ \b_2 < \b_3+\b_1 < 3/(5z_1)<\k_1$, a contradiction.\\

\noindent -$[132313231\vert k_j \cdots k_{j+8}]$: If $k_{j+3} \geq k_{j+1}+2$, then $\b_2  > 2/z_3- \b_3$ which gives $z_3>4z_2$ and so $2\b_3+ \b_2 <3/(5z_2)=\k_1$, a contradiction. If  $k_{j+6} \geq k_{j+2}+2$, then, since the subchain $[323]$ gives $z_3>3z_2/2$, the subchain $[23132]$ gives $2/z_2-\b_2<\b_1+2\b_3<\b_1+4/(15z_2)$ which implies that $z_2>10z_1/3$ and so $ \b_2+2\b_3<3/(25z_1)+4/(25z_1)=7/(25z_1)<\k_1$, a contradiction. 
\end{proof}

Define the two integers $a, b \geq 1$ by (note that $z_1>1$)
\[
a = \left\lceil \frac{z_1-1}{5} \right\rceil \quad \mbox{and} \quad b= \left\lfloor \frac{4z_1+1}{5} \right\rfloor.
\] 
By construction, $a$ is the least integer $k$ such that $\mathcal{B}_k^+/z_1\geq \delta_3=1/5$ and $b$ is the largest integer $k$ such that $\mathcal{B}_k^-/z_1\leq 1-\delta_3=4/5$. 

\begin{lemma}\label{LEM3oct1}
If $z_1\geq 7/2$, then we have $b-a <  \lceil 4a/5 + 4/25 \rceil +1$.
\end{lemma}

\begin{proof}[Proof of Lemma \ref{LEM3oct1}]
We suppose for contradiction that $z_1\geq 7/2$ and $b-a \geq  \lceil 4a/5 + 4/25 \rceil +1$. Then we have $a\geq \lceil (7/2-1)/5 \rceil=1$,  $b\geq \lfloor (4\cdot 7/2+1)/5 \rfloor=3$, $b-a\geq  \lceil 4/5 + 4/25 \rceil +1=2$ and the set $\mathcal{K}_0$ contains the $b-a$  consecutive $1$-kwais 
\[
\mathcal{K}_{a}/z_1,  \quad \mathcal{K}_{a+1}/z_1, \quad \cdots, \quad  \mathcal{K}_{b-1}/z_1.
\]
Thanks to Lemma \ref{LEM3Subcoverings} we know that all consecutive $1$-kwais above are covered by unions of $2$-bridges and $3$-bridges in the same manner. It means that there are $\bar{i} \in \{2,3\}$ and  $j\in \{1,\ell\}$ such that 
\[
\mathcal{K}_a^+/z_1 \in \mathcal{B}_{k_j}/z_{\bar{i}}
\]
and we have, either $\bar{i}=2$ and $[i_ji_{j+1} i_{j+2}]=[231]$, or  $\bar{i}=3$ and $[i_ji_{j+1} i_{j+2}]=[321]$, or $\bar{i}=3$ and $[i_j \cdots i_{j+3}]= [3231]$. Let us treat each case separately.\\

\noindent - $[i_ji_{j+1} i_{j+2}]=[231]$: We have
$$
 \mathcal{B}_{a}^-/z_1  =      \frac{5a-1}{5z_1} < \mathcal{B}_{k_j}^-/z_2 = \frac{5k_j-1}{5z_2} < \mathcal{B}_{a}^+/z_1  =      \frac{5a+1}{5z_1} < \mathcal{B}_{k_j}^+/z_2  = \frac{5k_j+1}{5z_2},
$$
which gives, since $z_2>z_1$, $k_j \geq a+1$ and as a consequence
\[
z_2 > \left( \frac{5k_j-1}{5a+1} \right) z_1 \geq  \left( \frac{5a+4}{5a+1} \right)z_1 \geq (1+\kappa) z_1 \quad \mbox{with} \quad \kappa := \frac{3}{5a+1}.
\] 
Define the sequence $\{u_k\}_{k\in \N}$ by 
$$
u_k := \mathcal{B}_{k_j+k}^+/z_2 - \mathcal{B}_{a+k}^+/z_1 =  \mathcal{B}_{k_j}^+/z_2 - \mathcal{B}_{a}^+/z_1 + k \left(\frac{1}{z_2}-\frac{1}{z_1}\right) \qquad \forall k\in \N.
$$
By construction, we have (because $u_0<\b_2$ and  $u_0+\b_3>\k_1$)
$$
u_0 =  \mathcal{B}_{k_j}^+/z_2 - \mathcal{B}_{a}^+/z_1 \in (1/(5z_2),2/(5z_2))
$$
and
$$
u_{k+1}-u_k = \frac{1}{z_2}- \frac{1}{z_1} <0 \qquad \forall k\in \N. 
$$
Therefore, $\{u_k\}_{k\in \N}$ is decreasing and there is $k\in \N^*$ such that 
\[
u_k<\frac{1}{5z_2} \quad \mbox{and} \quad 
k\leq  \left\lceil \frac{1/(5z_2)}{1/z_1-1/z_2} \right\rceil \leq  \left\lceil \frac{1}{5\kappa} \right\rceil=  \left\lceil\frac{a}{3} + \frac{1}{15}  \right\rceil.
\]
If all the $1$-kwais from $\mathcal{K}_{a}/z_1$ to $\mathcal{K}_{a+k}/z_1$ are contained in $\mathcal{K}_0$, then by Lemma \ref{LEM3Subcoverings}, we infer that the $1$-kwai $\mathcal{K}_{a+k}/z_1$ satisfies
\[
\mathcal{K}_{a+k}/z_1 \subset \mathcal{B}_{k_{j}+k }/z_2 \cup \mathcal{B}_{k_{j+1+3k} }/z_3 \quad \mbox{with} \quad \mathcal{B}_{k_{j}+k }^+/z_2<\mathcal{K}_{a+k}^-/z_1 + \frac{1}{5z_2}
\]
which is impossible. The above argument applies provided  $\mathcal{K}_{a+k}/z_1\subset \mathcal{K}_0$ with $1\leq k\leq \lceil a/3 + 1/15\rceil$ which is satisfied if  $a+\lceil a/3 + 1/15\rceil+1\leq b$ (note that $\lceil 4a/5 + 4/25 \rceil \geq \lceil a/3 + 1/15\rceil$).\\

\noindent - $[i_ji_{j+1} i_{j+2}]=[321]$: We proof is exactly the same as before by replacing $2$ by $3$, it is left to the reader.\\

\noindent  - $[i_ji_{j+1} i_{j+2}i_{j+3}]=[3231]$:  The subchain $[323]$ yields $z_3>3z_2/2$, so we have
$$
 \mathcal{B}_{a}^+/z_1  =      \frac{5a+1}{5z_1} < \mathcal{B}_{k_{j+1}}^-/z_2 = \frac{5k_{j+1}-1}{5z_2} < \mathcal{B}_{a}^+/z_1 + \frac{2}{5z_3} <  \mathcal{B}_{a}^+/z_1 + \frac{4}{15z_2}
$$
which gives, since $z_2>z_1$, $k_{j+1} \geq a+1$ and $z_2 \geq (1+\kappa) z_1$ with $\kappa := 5/(15a+3)$. Then the sequence  $\{u_k\}_{k\in \N}$ defined by 
$$
u_k := \mathcal{B}_{k_{j+1}+k}^-/z_2 - \mathcal{B}_{a+k}^+/z_1 =  \mathcal{B}_{k_{j+1}}^-/z_2 - \mathcal{B}_{a}^+/z_1 + k \left(\frac{1}{z_2}-\frac{1}{z_1} \right) \qquad \forall k\in \N,
$$
verifies 
\[
u_0\in ]0,4/(15z_2)] \quad \mbox{and} \quad u_{k+1}-u_k = \frac{1}{z_2}-\frac{1}{z_1}  <0 \quad \forall k\in \N.
\]
Consequently, there is $k\in \N^*$ such that
\[
u_k<0 \quad \mbox{and} \quad k\leq  \left\lceil \frac{4/(15z_2)}{1/z_1-1/z_2} \right\rceil \leq  \left\lceil \frac{4}{15\kappa} \right\rceil=  \left\lceil\frac{4a}{5} + \frac{4}{25}  \right\rceil.
\]
If all the $1$-kwais from $\mathcal{K}_{a}/z_1$ to $\mathcal{K}_{a+k}/z_1$ are contained in $\mathcal{K}_0$, then by Lemma \ref{LEM3Subcoverings}, we infer that the $1$-kwai $\mathcal{K}_{a+k}/z_1$ satisfies
\[
\mathcal{K}_{a+k}/z_1 \subset \mathcal{B}_{k_{j+3k} }/z_3  \cup \mathcal{B}_{k_{j+1}+k }/z_2 \cup \mathcal{B}_{k_{j+2+3k} }/z_3 \quad \mbox{with} \quad \mathcal{B}_{k_{j+1}+k }^-/z_2>  \mathcal{K}_{a+k}^-/z_1 = \mathcal{B}_{a+k}^+/z_1
\]
which contradicts $u_k<0$.  The above argument applies provided  $\mathcal{K}_{a+k}/z_1\subset \mathcal{K}_0$ with $1\leq k\leq \lceil 4a/5 + 4/25\rceil$ which is satisfied if  $a+\lceil 4a/5 + 4/25\rceil+1\leq b$.\\

In conclusion, we have proved that whatever the starting $1$-subchain of length $1$, we always obtain a contradiction if $b-a \geq  \lceil 4a/5 + 4/25 \rceil +1$. This shows that $b-a <  \lceil 4a/5 + 4/25 \rceil +1$. 
\end{proof}

\begin{lemma}\label{LEMHanoiairport}
$z_1\geq 7/2 \Rightarrow b-a \geq  \lceil 4a/5 +4/25\rceil +1$.
\end{lemma}

\begin{proof}[Proof of Lemma \ref{LEMHanoiairport}]
Let $z_1\geq 7/2$ be fixed. If $z_1\geq 6$, then there are $q\in \N$ with $q\geq 1$ and $r\in [0,5)$ such that $z_1=1+5q+r$. Then we have 
\[
a = \left\lceil \frac{z_1-1}{5} \right\rceil =  \left\lceil q+\frac{r}{5} \right\rceil = q+  \left\lceil \frac{r}{5} \right\rceil, \quad b =  \left\lfloor \frac{4z_1+1}{5} \right\rfloor =  \left\lfloor 4q+1+ \frac{4r}{5} \right\rfloor =  4q+1 +  \left\lfloor \frac{4r}{5} \right\rfloor
\]
and 
\[
\left\lceil \frac{4a}{5} + \frac{4}{25} \right\rceil +1 = \left\lceil \frac{4q}{5} + \frac{4}{5}   \left\lceil \frac{r}{5} \right\rceil +\frac{4}{25}  \right\rceil +1 \leq  \left\lceil \frac{4q}{5} + \frac{4}{5}   +\frac{4}{25}  \right\rceil  +1 <  \frac{4q}{5} + \frac{74}{25}.
\]
Then we have $b-a\geq 3q $ which is $\geq  4q/5 +4/25  $ since $q\geq 1$. If $z_1<6$, then we have $a=1$ hence  $b-a \geq  \lceil 4a/5 +4/25\rceil +1$ is equivalent to $b\geq 3$, that is, $z_1\geq 7/2$.
\end{proof}

In conclusion, Lemmas \ref{LEM3oct1} and \ref{LEMHanoiairport} show that the $\mathcal{K}_0$-covering $(z_1,z_2,z_3)$ satisfies $z_1<7/2$. Let us now prove that no $\mathcal{K}_0$-covering satisfies this property. If $z_1<7/2$, then we have
\[
\mathcal{B}_1^+/z_1> \mathcal{B}_1^-/z_1=\frac{4}{5z_1} > \frac{8}{35} > \frac{1}{5} \quad \mbox{and} \quad \mathcal{B}_3^-/z_1=\frac{14}{5z_1}> \frac{4}{5},
\]
so that  we may have $\mathcal{K}_1/z_1 \subset \mathcal{K}_0$ or $\mathcal{K}_1/z_1 \subsetneq \mathcal{K}_0$. We treat the two cases separately:  \\

\noindent Case 1: $\mathcal{K}_1/z_1 \subset \mathcal{K}_0$.\\
The $1$-kwai $\mathcal{K}_1/z_1$ can be covered by the following subchains: $[23]$, $[32]$, or $[323]$, we treat each case separately. \\

\noindent Case 1.1: $\mathcal{K}_1/z_1$ is covered by $[23]$.\\
Then there are $k_j,k_{j+1}\in \N^*$ such that
\begin{multline*}
\mathcal{B}_{k_j}^-/z_2 = \frac{5k_j-1}{5z_2} < \mathcal{B}_{1}^+/z_1  =  \frac{6}{5z_1} < \mathcal{B}_{k_{j+1}}^-/z_3 = \frac{5k_{j+1}-1}{5z_3}\\
 < \mathcal{B}_{k_j}^+/z_2 = \frac{5k_j+1}{5z_2} <  \mathcal{B}_{2}^-/z_1  = \frac{9}{5z_1} < \mathcal{B}_{k_{j+1}}^+/z_3  = \frac{5k_{j+1}+1}{5z_3}.
\end{multline*}
Then we have $(5k_j+1)>6z_2/z_1$ which gives $k_j\geq 2$ and as a consequence $z_2/z_1>(5k_j-1)/6\geq 3/2$. But we have 
\[
\frac{4}{5z_2} > \b_2+\b_3=\frac{2}{5z_2} + \frac{2}{5z_3} > \k_1 = \frac{3}{5z_1},
\]
so that $z_2/z_1<4/3$. We obtain a contradiction because $3/2>4/3$. \\

\noindent Case 1.2: $\mathcal{K}_1/z_1$ is covered by $[32]$.\\
Then there are $k_j,k_{j+1}\in \N^*$ such that
\begin{multline*}
\mathcal{B}_{k_j}^-/z_3 = \frac{5k_j-1}{5z_3} < \mathcal{B}_{1}^+/z_1  =  \frac{6}{5z_1} < \mathcal{B}_{k_{j+1}}^-/z_2 = \frac{5k_{j+1}-1}{5z_2}\\
 < \mathcal{B}_{k_j}^+/z_3 = \frac{5k_j+1}{5z_3} <  \mathcal{B}_{2}^-/z_1  =      \frac{9}{5z_1} < \mathcal{B}_{k_{j+1}}^+/z_2  = \frac{5k_{j+1}+1}{5z_2}.
\end{multline*}
Then, $5k_{j}+1>6z_3/z_1>6$ gives $k_j\geq 2$,  $5k_{j+1}+1>9z_2/z_1>9$ gives $k_{j+1}\geq 2$ and $5k_j+1>(5k_{j+1}-1)z_3/z_2>5k_{j+1}-1$ gives $k_{j+1}\leq k_j$. But we have 
\[
\frac{2}{5z_1} + \frac{2}{5z_3} > \frac{2}{5z_2} + \frac{2}{5z_3} = \b_2+\b_3> \k_1=\frac{3}{5z_1},
\]
which implies $z_3/z_1<2$. If $k_j\geq 3$, we have $z_3/z_1>(5k_j-1)/6\geq 7/3>2$, a contradiction. Thus, we have $k_j=k_{j+1}=2$ and
\[
\frac{9}{z_3} <   \frac{6}{z_1} <  \frac{9}{z_2} < \frac{11}{z_3} <   \frac{9}{z_1} <  \frac{11}{z_2}.
\]
So we have on the one hand $z_3/z_1>9/6=3/2$ and on the other hand $z_3/z_1 = (z_3/z_2)\cdot(z_2/z_1)< 11^2/9^2<3/2$, a contradiction. \\

\noindent Case 1.3: $\mathcal{K}_1/z_1$ is covered by $[323]$.\\
Then there are $k_j,k_{j+1},k_{j+2} \in \N^*$ such that
\begin{multline*}
\mathcal{B}_{k_j}^-/z_3 = \frac{5k_j-1}{5z_3} < \mathcal{B}_{1}^+/z_1  =  \frac{6}{5z_1} < \mathcal{B}_{k_{j+1}}^-/z_2 = \frac{5k_{j+1}-1}{5z_2}\\
 < \mathcal{B}_{k_j}^+/z_3 = \frac{5k_j+1}{5z_3}  < \mathcal{B}_{k_{j+2}}^-/z_3 = \frac{5k_{j+2}-1}{5z_3} \\
 <  \mathcal{B}_{k_{j+1}}^+/z_2  = \frac{5k_{j+1}+1}{5z_2}< \mathcal{B}_{2}^-/z_1  =      \frac{9}{5z_1} < \mathcal{B}_{k_{j+2}}^+/z_3  = \frac{5k_{j+2}+1}{5z_3}.
\end{multline*}
Then, $5k_j+1>(5k_{j+1}-1)z_3/z_1>5k_{j+1}-1$ gives $k_{j}\geq k_{j+1}$ and $5k_{j+1}-1>6z_2/z_1>6$ gives $k_{j+1}\geq 2$. If $k_j\geq 5$, then $z_3/z_1>(5k_j-1)/6\geq 4$, a contradiction (because $\b_2+2\b_3< 2/(5z_1)+1/(5z_3)<\k_1$). So we have $k_j\in \{2,3,4\}$. But we have
\[
\frac{5k_{j+2}-1}{5k_{j+1}+1} < \frac{z_3}{z_2} < \frac{5k_j+1}{5k_{j+1}-1}.
\]
So the admissible pairs for $(k_j,k_{j+1},k_{j+2})$ are $(3,2,4)$ and $(4,2,5)$. In the first case,  $\mathcal{B}_{k_j}^-/z_3 < \mathcal{B}_{1}^+/z_1$ and $\mathcal{B}_{2}^-/z_1  < \mathcal{B}_{k_{j+2}}^+/z_3$ give respectively $z_3/z_1>7/3$ and  $z_3/z_1<7/3$, a contradiction. In the second case, the same inequalities imply $19/6<7/3<26/9$, a contradiction again.\\

\noindent Case 2: $\mathcal{K}_1/z_1 \subsetneq \mathcal{K}_0$.\\
Since $\mathcal{B}_1^-/z_1 > \frac{1}{5}$, the assumption means that $\mathcal{B}_{2}^-/z_1  =  9/(5z_1)>4/5$, which gives $z_1<9/4$. By Proposition \ref{LEMTHMmain2}, the chain associated to the $\mathcal{K}_0$-covering $(z_1,z_2,z_3)$ has to be in the following list: $(21)$, $(212)$, $(2123)$, $(2132)$, $(21323)$, $(31)$, $(312)$, $(3123)$, $(3132)$, $(31323)$, $(231)$, $(2312)$, $(23123)$, $(23132)$, $(231323)$, $(321)$, $(3212)$, $(32123)$, $(32132)$, $(321323)$, $(3231)$, $(32312)$, $(323123)$, $(323132)$, $(3231323)$. We need to distinguish several cases. \\

\noindent Case 2.1: The chain has the form $(2\cdots)$.\\
Then $1/5 \in \mathcal{B}/z_2$ and there is $k\in \N^*$ such that $1/5\in [\mathcal{B}_k^-/z_2,\mathcal{B}_k^+/z_2 ]=[(5k-1)/(5z_2),(5k+1)/(5z_2)]$, which shows that $z_3>z_2>4$ and gives $\b_3< \b_2 <1/10$. The interval $[1/5,\mathcal{B}_1^-/z_1]=[1/5,4/(5z_1)]$ can be covered either by a $z_2$-bridge or by the union of a $z_2$-bridge and a $z_3$-bridge, so we have $4/(5z_1)-1/5 < \b_2 + \b_3 < 1/5$, which gives $z_1>2$ and as a consequence $4/5-\mathcal{B}_1^+/z_1>1/5$. Thus the interval $[\mathcal{B}_1^+/z_1,4/5]$  has to be covered by a subchain $[323]$. So we have $z_3>3z_2/2>6$ and $\b_3 < 1/15$ which implies that $4/(5z_1)=\mathcal{B}_1^-/z_1<1/5+\b_2+\b_3<11/30$, that is, $z_1>24/11$. Then the interval $[\mathcal{B}_1^+/z_1,4/5]$ has length $4/5-\mathcal{B}_1^+/z_1>1/4$ and cannot be covered by a subchain $[323]$ (because $\b_2+2\b_3<7/30<1/4$), a contradiction.\\

\noindent Case 2.2: The chain has the form $(3\cdots)$.\\
Then $1/5 \in \mathcal{B}/z_3$ and as in the previous case, we infer that $z_3>4$ and $\b_3<1/10$. Since $4/(5z_1)-1/5>7/45>\b_3$, the interval $[1/5,\mathcal{B}_1^-/z_1]$ can be covered either by a subchain $[32]$ or by a subchain $[323]$ and so we can leave out the cases  $(31)$, $(312)$, $(3123)$, $(3132)$, $(31323)$. Let us now distinguish two cases.\\

\noindent Case 2.2.1:  The chain has the form $(321\cdots)$.\\
The inequality $(5k_2-1)/(5z_2)=\mathcal{B}_{k_2}^-/z_2<1/5+\b_3<3/10$ gives $z_2>2(5k_2-1)/3$. If $k_2\geq 2$, then we have $z_3>z_2> 6$ which implies $\b_3+\b_2<2/15<7/45$, a contradiction (remember that $4/(5z_1)-1/5>7/45$). Hence we have $k_2=1$ and $z_2>8/3$. We also have $6/(5z_2)=\mathcal{B}_{k_2}^+/z_2>\mathcal{B}_1^-/z_1=4/(5z_1)$ which gives $z_2/z_1<3/2$ and $z_1>2z_2/3>16/9$. So there is no subchain of the form $[212]$ (which implies $z_2/z_1>3/2$) and our chain has to be  $(3213\vert k_1k_2k_3 k_4)$, $(32132\vert k_1k_2\cdots k_5)$ or $(321323\vert k_1k_2\cdots k_6)$. If it is $(3213\vert k_1k_2k_3 k_4)$, then we have $4/5-6/(5z_1)=4/5-\mathcal{B}_1^+/z_1<\b_3<1/10$, which gives $z_1<12/7<16/9$, a contradiction.  If the chain is $(32132\vert k_1k_2\cdots k_5)$ or $(321323\vert k_1k_2\cdots k_6)$, then we have
\[
\frac{5k_1-1}{5z_3} < \frac{1}{5} < \frac{4}{5z_2} < \frac{5k_1+1}{5z_3}  < \frac{4}{5z_1}  < \frac{6}{5z_2} <  \frac{5k_4-1}{5z_3}<    \frac{6}{5z_1} < \frac{5k_5-1}{5z_2} <  \frac{5k_4+1}{5z_3}  < \frac{4}{5},
\]
which gives $z_2>(5k_5-1)/4$ with $k_5\geq k_2+1\geq 2$. Note that we have necessarily $k_5=2$ because otherwise we would have $z_2>7/2>(3/2)\cdot (9/4)>3z_1/2$, a contradiction. If $k_1\geq 2$, then $1/5\in \mathcal{B}_{k_1}/z_3$ implies $z_3>9$ and $\b_3<2/45$, then $4/(5z_2)=\mathcal{B}_{k_2}^-/z_2<1/5+\b_3$ gives $z_2>36/11$ and $\b_2<11/90$, so that $4/(5z_1)-1/5<\b_3+\b_2$ yields $z_1>24/11$ and then the interval $[\mathcal{B}_1^+/z_1,4/5]$ has length $4/5-6/(5z_1)>1/4$ and cannot be filled with the union of a $z_2$-bridge and two $z_3$-bridges ($\b_2+2\b_3<19/90$). In conclusion, we have $k_1=1$. We now claim that $k_4=3$. Since $z_2<3z_1/2<27/8$, we have
\[
\frac{32}{27} < \frac{z_3}{z_2} < \frac{5k_4+1}{9}
\]
which implies $k_4\geq 3$. If $k_4=4$, then we have $z_3/z_2>19/4$ which implies $z_3>19z_2/4>171/16>9$ and as before we obtain a contradiction. In conclusion, our chain has to be $(32132\vert 11132)$ or $(321323\vert 11132k_6)$, so that we have
\[
\frac{4}{z_3} < 1 < \frac{4}{z_2} < \frac{6}{z_3}  < \frac{4}{z_1}  < \frac{6}{z_2} <  \frac{14}{z_3}<  \frac{6}{z_1} < \frac{9}{z_2} <  \frac{16}{z_3}  < 4.
\]
Therefore, we have $z_3/z_2<3/2$,   $z_2/z_1>2/3$, so that $z_3/z_1<9/4$. But we also have $z_3/z_1>7/3$, this is a contradiction because $7/3>9/4$.\\

\noindent Case 2.2.2:  The chain has the form $(3231\cdots)$.\\
The subchain $[323]$ gives $z_3>3z_2/2$. The inequality $(5k_2-1)/(5z_2)=\mathcal{B}_{k_2}^-/z_2<1/5+\b_3<3/10$ gives $z_2>2(5k_2-1)/3$. If $k_2\geq 2$, then we have $z_2> 6$ and $z_3>9$ which implies $2\b_3+\b_2<2/15<7/45$, a contradiction (remember that $4/(5z_1)-1/5>7/45$). Hence we have $k_2=1$ and $z_2>8/3$. We also have $6/(5z_2)=\mathcal{B}_{k_2}^+/z_2>\mathcal{B}_1^-/z_1=4/(5z_1)$ which gives $z_2/z_1<3/2$ and $z_1>2z_2/3>16/9$. So there is no subchain of the form $[212]$ (which implies $z_2/z_1>3/2$) and our chain has to be  $(32313\vert k_1\cdots k_5)$, $(323132\vert k_1\cdots k_6)$ or $(3231323\vert k_1\cdots k_7)$. If it is $(32313\vert k_1\cdots k_()$, then we have $4/5-6/(5z_1)=4/5-\mathcal{B}_1^+/z_1<\b_3<1/10$, which gives $z_1<12/7<16/9$, a contradiction.  If the chain is $(323132\vert k_1\cdots k_6)$ or $(3231323\vert k_1\cdots k_7)$,  then we have
\begin{multline*}
\frac{5k_1-1}{5z_3} < \frac{1}{5} < \frac{4}{5z_2} < \frac{5k_1+1}{5z_3} < \frac{5k_3-1}{5z_3} < \frac{6}{5z_2} < \frac{4}{5z_1} \\
 < \frac{5k_3+1}{5z_3} <  \frac{5k_5-1}{5z_3}<    \frac{6}{5z_1} < \frac{5k_6-1}{5z_2} <  \frac{5k_5+1}{5z_3}  < \frac{4}{5},
\end{multline*}
which gives $z_2>(5k_6-1)/4$ with $k_6\geq k_2+1\geq 2$. Note that we have necessarily $k_6=2$ because otherwise we would have $z_2>7/2>(3/2)\cdot (9/4)>3z_1/2$, a contradiction. If $k_3\geq 3$, then $z_3>7z_2/3$ which implies $4/(5z_2)=\mathcal{B}_{k_2}^-/z_2<1/5+\b_3<1/5+6/(35z_2)$ which gives $z_2>34/7>27/8$, a contradiction.  So we have $k_1=1$ and $k_3=2$.  Therefore, we have
\[
\frac{4}{z_3} < 1 < \frac{4}{z_2} < \frac{6}{z_3}  <  \frac{9}{z_3} < \frac{6}{z_2} < \frac{4}{z_1}   <  \frac{11}{z_3},
\]
which gives both $z_3/z_2>3/2$ and $z_3/z_2<3/2$, a contradiction. 

\section{Proof of Theorem \ref{THM2}}\label{SECProofTHM2}

We have $\delta_4=1/6$, $\mathcal{K}_0=[1/6,5/6]$, $\mathcal{K}_1=[7/6,11/6]$, $\mathcal{K}_{\llbracket 0,1\rrbracket}=\mathcal{K}_0 \cup \mathcal{K}_1$, $\b(z)=2/(6z)=1/(3z)$, $\k(z)=4/(6z)=2/(3z)$ for all $z\geq 1$ and we need to show that 
\[
 [1,\infty)^4 \subset \mathcal{K}^{4}\left(\delta_4\right) / \mathcal{K}_{\llbracket 0,1\rrbracket}  \quad \mbox{and} \quad  [1,\infty)^4 \subsetneq \mathcal{K}^{4}\left(\delta_4\right) / \mathcal{K}_{0}. 
 \]
 We note that by symmetry, the first inclusion is equivalent to 
 \[
  \mathcal{S}^4 \subset \mathcal{K}^{4}\left(\delta_4\right) / \mathcal{K}_{\llbracket 0,1\rrbracket},
 \]
 where we recall that $\mathcal{S}^4$ stands for the closed convex set of tuples $(z_1,z_2,z_3,z_4)$ such that $1\leq z_1\leq z_2\leq z_3 \leq z_4$. Our proof is divided in three parts as follows:\\

\noindent - In Section \ref{SEC4Dinfinity}, we show that there is a compact set $P\subset \mathcal{S}^4$ such that
\[
\mathcal{S}^4 \setminus P \subset \mathcal{K}^{4}\left(\delta_4\right) / \mathcal{K}_{0}.
\]
As in the proof of Theorem \ref{THM1} for $d=3$, the result follows from particular dynamical properties satisfied by $1$-subchains of length $1$. However, unlike the $d=3$ case, the dynamics of $1$-subchains of length $1$ is not trivial, so we perform its study will help of a computer. \\

\noindent - In Section \ref{SEC4Dcompact}, we check that the required covering property is satisfied over $P$, that is,
\[
 P \subset \mathcal{K}^{4}\left(\delta_4\right) /  \mathcal{K}_{\llbracket 0,1\rrbracket}.
\]
In fact, we explain how to check covering properties over compact sets with computer-assisted proofs and apply the method to our set $P$. \\

\noindent - In Section \ref{SEC4Dnot}, we verify that 
\[
[1,\infty)^4 \subsetneq \mathcal{K}^{4}\left(\delta_4\right) / \mathcal{K}_{0}. 
\]
To do this, we simply check by hand that some points of $[1,\infty)^4$ given by a formula are not covered by $\mathcal{K}^{4}\left(\delta_4\right) / \mathcal{K}_{0}$. 

\subsection{Proof of the covering property at infinity}\label{SEC4Dinfinity}
\paragraph{Introduction.}  The aim of Section \ref{SEC4Dinfinity} is to prove the following result:

\begin{proposition}\label{PROP4Dinfinity}
We have 
\[
 \mathcal{S}^4 \setminus P \subset \mathcal{K}^{4}\left(\delta_4\right) /  \mathcal{K}_{0},
\]
where $P\subset \mathcal{S}^4$ is the compact set consisting of the points $(z_1,z_2,z_3,z_4)\in \mathcal{S}^4$ satisfying 
\[
5z_1 \leq 47, \quad 2z_2\leq 5z_1, \quad z_2z_3+z_1z_3 \leq 8z_1z_2, \quad z_3z_4+z_2z_4\leq 10z_2z_3. 
\]
\end{proposition}

Our strategy to prove Proposition \ref{PROP4Dinfinity} is to show by Proposition \ref{PROP3oct1} that any $\mathcal{K}_0$-covering has to be in $P$. So, let us consider a $\mathcal{K}_0$-covering $(z_1, z_2, z_3,z_4)$ associated with a chain $(i_1\cdots i_{\ell}\vert k_1\cdots k_{\ell})$ with $\ell\geq 2$. First of all, we note that the property $N_3=1$ given by Theorem \ref{THM1} allows to prove the following result:

\begin{lemma}\label{LEM12oct}
We have 
\begin{eqnarray}\label{12oct1}
\frac{10}{z_4}>\frac{1}{z_2} + \frac{1}{z_3}   \quad \mbox{and} \quad \frac{8}{z_3} >\frac{1}{z_1} + \frac{1}{z_2}.
\end{eqnarray}
\end{lemma}

\begin{proof}[Proof of Lemma \ref{LEM12oct}]
Since $N_3=1$, we know that  $(z_1,z_2,z_3)$ is not a covering of $\mathcal{K}_0(\delta_3)=\mathcal{K}_0(1/5)=[1/5,4/5]$. So, by Proposition \ref{PROPcaracx}, the set 
\[
S := \mathcal{K}_0(1/5)\cap \mathcal{K}(1/5)/z_1 \cap \mathcal{K}(1/5)/z_2  \cap \mathcal{K}(1/5)/z_3
\]
is not empty. Let $\lambda\in S$ be fixed. Since 
\[
S\subset  \mathcal{K}_0(1/6)\cap \mathcal{K}(1/6)/z_1 \cap \mathcal{K}(1/6)/z_2  \cap \mathcal{K}(1/6)/z_3
\]
and $(z_1,z_2,z_3,z_4)$ is assumed to be a $\mathcal{K}_0(1/6)$-covering, $\lambda$ belongs to some $\mathcal{B}_{k_j}(1/6)/z_4$. We distinguish several cases:\\

\noindent Case 1: $j=1$.\\
Then we have 
\[
\mathcal{B}_{k_1}^-(1/6)/z_4 <\frac{1}{6} < \frac{1}{5} \leq \lambda < \mathcal{B}_{k_1}^+(1/6)/z_4
\]
and, since $\lambda \in \mathcal{K}_0(1/5) \setminus \mathcal{B}_{k_2}(1/5)$,
\[
\mathcal{B}_{k_1}^-(1/6)/z_4 <\lambda < \mathcal{B}_{k_2}^-(1/5)/z_{i_2}< \mathcal{B}_{k_2}^-(1/6)/z_{i_2} < \mathcal{B}_{k_1}^+(1/6)/z_4 \quad \mbox{with} \quad i_2 \in \{1,2,3\}.
\]
We infer that
\begin{eqnarray*}
\frac{1}{3z_4} = \b_4 & = & \mathcal{B}_{k_1}^+(1/6)/z_4 - \mathcal{B}_{k_1}^-(1/6)/z_4\\
& = & \left( \mathcal{B}_{k_1}^+(1/6)/z_4 -\lambda \right) + \left( \lambda - \mathcal{B}_{k_1}^-(1/6)/z_4\right)\\
& > &  \left( \frac{1}{5} -\frac{1}{6} \right) \frac{1}{z_{i_2}} + \left( \frac{1}{5} -\frac{1}{6} \right)\geq  \frac{1}{30} \left( 1+\frac{1}{z_3}\right).
\end{eqnarray*}

\noindent Case 2: $j=\ell$.\\
The proof follows the same line as Case 1, it is left to the reader.\\

\noindent Case 3: $j\notin \{1,\ell\}$.\\
Then we have 
\[
\mathcal{B}_{k_{j-1}}^-(1/6)/z_{i_{j-1}} < \mathcal{B}_{k_j}^-(1/6)/z_4 < \mathcal{B}_{k_{j-1}}^+(1/6)/z_{i_{j-1}} < \mathcal{B}_{k_{j-1}}^+(1/5)/z_{i_{j-1}} \leq \lambda
\]
and
\[
\lambda \leq   \mathcal{B}_{k_{j+1}}^-(1/5)/z_{i_{j+1}}<  \mathcal{B}_{k_{j+1}}^-(1/6)/z_{i_{j+1}} < \mathcal{B}_{k_j}^+(1/6)/z_4,
\]
with $i_{j-1},i_{j+1} \in \{1,2,3\}$ and $i_{j-1}\neq i_{j+1}$ (see Proposition \ref{LEMTHMmain2}). We infer that
\begin{eqnarray*}
\frac{1}{3z_4} = \b_4 & = & \mathcal{B}_{k_1}^+(1/6)/z_4 - \mathcal{B}_{k_1}^-(1/6)/z_4\\
& = & \left( \mathcal{B}_{k_1}^+(1/6)/z_4 -\lambda \right) + \left( \lambda - \mathcal{B}_{k_1}^-(1/6)/z_4\right)\\
& > &  \left( \frac{1}{5} -\frac{1}{6} \right) \frac{1}{z_{i_{j+1}}} + \left( \frac{1}{5} -\frac{1}{6} \right) \frac{1}{z_{i_{j-1}}}\geq  \frac{1}{30} \left( \frac{1}{z_2}+\frac{1}{z_3}\right),
\end{eqnarray*}
because $i_{j-1},i_{j+1} \in \{1,2,3\}$ and $i_{j-1}\neq i_{j+1}$.\\

To prove the second inequality of (\ref{12oct1}), we note that since $(z_1,z_2)$ is not a covering of $\mathcal{K}_0(\delta_2)=\mathcal{K}_0(1/4)=[1/4,3/4]$, the set 
\[
\mathcal{K}_0(1/6) \setminus \left(  \mathcal{K}_0(1/4)\cap \mathcal{K}(1/4)/z_1 \cap \mathcal{K}(1/4)/z_2  \right)
\]
contains a closed interval of length at least $(1/4-1/6)(1/z_1+1/z_2)=1/12(1/z_1+1/z_2)$. Such an interval cannot be covered by unions of $3$-bridges and $4$-bridges if $2\b_3\leq 1/12(1/z_1+1/z_2)$, that is, if $8/z_3\leq 1/z_1+1/z_2$.
\end{proof}

Now we observe that the $\mathcal{K}_0$-covering $(z_1, z_2, z_3,z_4)$ is $1$-nice. As a matter of fact, a set $I$ of the form (\ref{SetREM2}) necessarily contains a $2$-kwai which has to be covered by the union of $3$-bridges and $4$-bridges. But, by Proposition \ref{LEMTHMmain2}, the subchain corresponding to the covering of the $2$-kwai has to be $[2342]$, $[2342]$ or $[24342]$. In the first two cases we have $\b_3 +\b_4=1/(3z_3)+1/(3z_4)<2/(3z_2)=\k_2$ which makes the covering property impossible and in the third case the subchain $[434]$ yields $\b_3>\k_4 \Leftrightarrow z_4>2z_3$ which implies $\b_3+2\b_4<\k_2$ and makes again the covering property impossible. In conclusion, all instances of $1$-bridges have to appear in any subchain connecting two $1$-bridges. Moreover, by Proposition \ref{LEMTHMmain2}, any $1$-subchain of length $1$ has to be in the set  
\begin{multline*}
\mathcal{L} = \left\{ [12341], [12431], [124341], [13241], [132341], [132431], [1324341], [14231], [142341],\right. \\
[142431], [1424341], [13421], [134231], [1342341], [134241], [1342431], [13424341], \\
[14321], [143231],  [1432341], [143241], [1432431], [14324341], [143421], [1434231],\\
\left. [14342341], [1434241], [14342431], [143424341]\right\}.
\end{multline*}

We need now to study the dynamics of $1$-subchains of length $1$, that is, to understand how they can succeed to each other. A way to do this is to relax the notion of chain by allowing translations to the bridges.

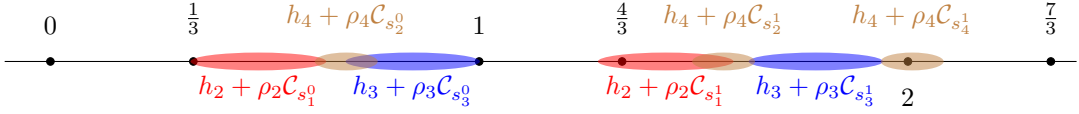
\begin{figure}[H]
\centering
\begin{tikzpicture}[scale=0.95]
\draw (-1,2) -- (14,2);
\filldraw[black] (-0.36,2)circle(1.5pt);
\draw (-0.36,2.25) node[above]{\color{black}{$0$}};
\filldraw[black] (1.64,2)circle(1.5pt);
\draw (1.64,2.25) node[above]{\color{black}{$\frac{1}{3}$}};
\filldraw[black] (5.64,2)circle(1.5pt);
\draw (5.64,2.25) node[above]{\color{black}{$1$}};
\filldraw[black] (7.64,2)circle(1.5pt);
\draw (7.64,2.25) node[above]{\color{black}{$\frac{4}{3}$}};
\filldraw[black] (11.64,2)circle(1.5pt);
\draw (11.64,1.75) node[below]{\color{black}{$2$}};
\filldraw[black] (13.64,2)circle(1.5pt);
\draw (13.64,2.25) node[above]{\color{black}{$\frac{7}{3}$}};
\fill [color=red, opacity=0.5] (2.55,2) ellipse (0.95cm and 0.125 cm);
\draw (2.55,1.25) node[above]{\color{red}{$h_{2} + \rho_{2} \mathcal{C}_{s_{1}^0}$}};
\fill [color=red, opacity=0.5] (8.25,2) ellipse (0.95cm and 0.125 cm);
\draw (8.25,1.25) node[above]{\color{red}{$h_{2} + \rho_{2} \mathcal{C}_{s_{1}^1}$}};
\fill [color=blue, opacity=0.5] (4.71,2) ellipse (0.94cm and 0.125 cm);
\draw (4.71,1.25) node[above]{\color{blue}{$h_{3} + \rho_{3} \mathcal{C}_{s_{3}^0}$}};
\fill [color=blue, opacity=0.5] (10.35,2) ellipse (0.94cm and 0.125 cm);
\draw (10.35,1.25) node[above]{\color{blue}{$h_{3} + \rho_{3} \mathcal{C}_{s_{3}^1}$}};
\fill [color=brown, opacity=0.5] (3.78,2) ellipse (0.44cm and 0.125 cm);
\draw (3.78,2.25) node[above]{\color{brown}{$h_{4} + \rho_{4} \mathcal{C}_{s_{2}^0}$}};
\fill [color=brown, opacity=0.5] (9.06,2) ellipse (0.44cm and 0.125 cm);
\draw (9.06,2.25) node[above]{\color{brown}{$h_{4} + \rho_{4} \mathcal{C}_{s_{2}^1}$}};
\fill [color=brown, opacity=0.5] (11.7,2) ellipse (0.44cm and 0.125 cm);
\draw (11.7,2.25) node[above]{\color{brown}{$h_{4} + \rho_{4} \mathcal{C}_{s_{4}^1}$}};
\end{tikzpicture}
\caption{A weak $2$-chain $\langle243\m 2434\vert s_1^0s_2^0s_3^0 \m s_1^1s_2^1s_3^1s_4^1\rangle $ associated with tuples $(\rho_2, \rho_3,\rho_4)$ and $(h_2,h_3,h_4)$ \label{FIGweakchain}}
\end{figure}

\paragraph{Admissible weak chains.} Given a tuple $(\rho_2,\rho_3,\rho_4)\in \R^3$ such that
\begin{eqnarray}\label{tuplerho}
0 < \rho_4 < \rho_3 < \rho_2 < 1
\end{eqnarray}
and a positive integer $L$, we call weak $L$-chain any family denoted as
\begin{eqnarray}\label{weakchain}
\langle i_1^0\cdots i_{\ell_0}^0 \m \cdots \m  i_1^{L-1}\cdots i_{\ell_{L-1}}^{L-1} \vert s_1^0\cdots s_{\ell_0}^0\m \cdots \m s_1^{L-1}\cdots s_{\ell_{L-1}}^{L-1} \rangle 
\end{eqnarray}
which satisfies the following properties (see Figure \ref{FIGweakchain}):
\begin{itemize}
\item[(P1)] For every $r=0, \ldots, L-1$, $\ell_r \in \N^*$ and $\ell_r\geq 3$.
\item[(P2)] For every $r=0, \ldots, L-1$ and every $j=1,\ldots, \ell_r$, $i_j^r\in \{2,3,4\}$ and $s_j^r\in \N$.
\item[(P3)] For every $j\in \{1,\cdots, \ell_0\}$ such that $i_k\neq i_j$ for all $k\in \llbracket 1,j-1\rrbracket$, we have $s_j=0$. 
\item[(P4)] There is $(h_2, h_3, h_4) \in \R^3$ such that we have for all $r=0, \ldots, L-1$, 
\[
r +2\delta_4 \in h_{i_1^r} + \rho_{i_1^r} \mathcal{C}_{s_{1}^r}, \quad (r+1) \in  h_{i_{\ell_r}^r } + \rho_{i_{\ell_r}^r} \mathcal{C}_{s_{\ell_r}^r},
\]
\[ 
h_{i_j^r} + \rho_{i_j^r} \mathcal{C}^+_{s_j^r} \in h_{i_{j+1}^r} + \rho_{i_{j+1}^r} \mathcal{C}_{s_{j+1}^r} \qquad \forall j=1, \ldots, \ell_r -1,
 \]
 \[
 r +2\delta_4 < h_{i_{2}^r} + \rho_{i_{2}^r} \mathcal{C}_{s_{2}^r}^-, \quad  h_{i_{\ell_r -1}^r} + \rho_{i_{\ell_r -1}^r}  \mathcal{C}_{s_{\ell_r -1}^r}^+ < r+1
 \]
 and
\[
h_{i_j^r} + \rho_{i_j^r}  \mathcal{C}_{s_{j}^r}^+ < h_{i_{j+2}^r} + \rho_{i_{j+2}^r} \mathcal{C}_{s_{j+2}^r}^-  \qquad \forall j=1, \ldots, \ell_r-2,
\]
where for every $s\in \N$, $\mathcal{C}_{s}, \mathcal{C}_{s}^-, \mathcal{C}_{s}^+$ are defined by
\[
\mathcal{C}_s  = \left( \mathcal{C}_{s}^-, \mathcal{C}_{s}^+\right) =  \left(s,s+2\delta_4\right).
\]
\end{itemize}

\vspace{0.2cm}
\vspace{0.2cm}

The sets $\rho_i\mathcal{C}_s(h_i)$ play exactly the same role as the $i$-bridges, their length is given by $2\delta_4\rho_i=\rho_i/3$ and two consecutive such sets are separated by an interval of length $\rho_i(1-2\delta_4)=2\rho_i/3$ (which plays the role of a $i$-kwai). 
Then, given a tuple a tuple $(\rho_2,\rho_3,\rho_4)\in \R^4$ satisfying (\ref{tuplerho}) and a weak chain as in (\ref{weakchain}), we set for every $r=0,\ldots, L-1$ and every $i\in \{2,3,4\}$,
 \begin{multline*}
 m_r(i) := \min \Bigl\{ s_j^r \, \vert \, j \in \{1,\ldots, \ell_r\}  \mbox{ s.t. }  i_j=i\Bigr\} \\
  \mbox{and} \quad M_r(i) := \max \Bigl\{ s_j^r \, \vert \, j \in \{1,\ldots, \ell_r\}  \mbox{ s.t. }  i_j=i\Bigr\}.
 \end{multline*}
We note that since the covering property
\[
\left[r+2\delta_4, r+1\right] \subset \bigcup_{j=1}^{\ell_r} \left( h_{i_j^r} + \rho_{i_j}^r \mathcal{C}_{s_{j}^r} \right)
\]
is satisfied for every $r=1,\ldots, L-1$, the same arguments as in the proof of Proposition \ref{LEMTHMmain2} show that for each $r=1,\ldots, L-1$, the set of $i_j^r$ with $j\in \{1,\ldots, \ell_r\}$ is equal to $\{2,3,4\}$, so that the min and maximum above are taken 
over non-empty sets. Moreover, if $L\geq 2$, then for every $r=0,\ldots, L-2$ and every $i\in \{2,3,4\}$, we call $i$-jump from $r$ to $r+1$ the quantity
\[
J_r(i) := m_{r+1}(i) - M_r(i).
\]
The following result explains how $1$-subchains are connected with weak chains. 

 \begin{lemma}\label{LEM20oct0}
 If there are $a,b\in \N$ such that $L:=b-a \geq 2$ and 
 \[
\mathcal{K}_{a}/z_1,  \quad \mathcal{K}_{a+1}/z_1, \quad \cdots, \quad  \mathcal{K}_{b-1}/z_1 \subset \mathcal{K}_0,
\]
then the tuple $(\rho_2,\rho_3,\rho_4)\in \R^3$  defined by 
\[
\rho_i := \frac{z_1}{z_i} \qquad \forall i=2,3,4
\]
admits a weak $L$-chain of the form (\ref{weakchain}) such that 
\begin{eqnarray}\label{11nov1}
\left[1 i_1^0\cdots i_{\ell_0}^0 1 \cdots 1 i_1^{L-1}\cdots i_{\ell_{L-1}}^{L-1} 1 \Big\vert s^0 \tilde{s}_1^0\cdots \tilde{s}_{\ell_0}^0s^1 \cdots \tilde{s}^{L-1} \tilde{s}_1^{L-1}\cdots \tilde{s}_{\ell_{L-1}}^{L-1} s^L\right],
\end{eqnarray}
where $s^0,\ldots, s^L$ and $\tilde{s}_1^0,\ldots, \tilde{s}_{\ell_0}^0, \ldots, \tilde{s}_1^{L-1}, \ldots, \tilde{s}_{\ell_{L-1}}^{L-1}$ are defined by
\[
s^r := a+r \quad \mbox{and} \quad \tilde{s}_j^r := k_{\bar{j}(i_j)} + s_j^r \qquad \forall r=0,\ldots, L-1, \, \forall j \in \{1,\ldots,\ell_r\},
\]
with 
\[
 \bar{j}(i) :=  \min \Bigl\{ j \in \{1,\ldots, \ell\} \, \vert \, i_j=i \, \mbox{ and } \, \mathcal{K}_a^-/z_1 < \mathcal{B}_{k_{j}}^+/z_i \Bigr\} \qquad \forall i\in \{2,3,4\},
\]
is a $1$-subchain associated with $(z_1,z_2,z_3,z_4)$ of length $L$.  Moreover, we have
  \begin{eqnarray}\label{10nov1}
5\rho_2 >2, \quad  10\rho_4>\rho_2 + \rho_3, \quad 8\rho_3 >1 + \rho_2,
\end{eqnarray}
\begin{eqnarray}\label{10nov2}
m_r(2)=M_r(2), \quad M_r(3)-m_r(3) \leq 4,  \quad M_r(4)-m_r(4) \leq 11 \qquad \forall r =0,\ldots, L-1 
\end{eqnarray}
and
\begin{eqnarray}\label{10nov3}
J_r(2) \leq 4, \quad J_r(3) \leq 9, \quad J_r(4) \leq 29  \qquad \forall r =0,\ldots, L-2.
\end{eqnarray} 
 \end{lemma}

\begin{proof}[Proof of Lemma \ref{LEM20oct0}]
Since $(z_1,z_2,z_3,z_4)$ is $1$-nice (see Definition \ref{DEF1subchain} and Remark \ref{REM2}) and 
\[
\mathcal{K}_{a}^-/z_1=\mathcal{B}^+_a/z_1, \mathcal{K}_{b-1}^+/z_1=\mathcal{B}^-_b/z_1 \in \mathcal{K}_0,
\] 
it admits a $1$-subchain of length $L=b-a$ of the form
\[
\left[1 i_j\cdots i_{j'}1\vert ak_{j}\cdots k_{j'}b\right],
\]
with $j<j'$ in $\{1,\ldots, \ell\}$. Then, we define the tuple $(h_1,h_2,h_3)\in \R^3$ by
\[
h_i :=  \frac{z_1}{z_i} \, \mathcal{B}^-_{k_{\bar{j}(i)}} - a+\delta_4.
\]
Properties (i)-(v) along with the fact that $1$-kwais can only be covered by subchains of length at least $3$ (we mean here with at least three indices $i_s$) imply that the chain defined by (\ref{11nov1}) is a weak $L$-chain.

To prove that $5\rho_2>2$ we note that the interval $I=[2\delta_4=1/3,1]$  is covered by an union of sets of the form $h_i + \rho_i \mathcal{C}_s$ with $i\in \{2,3,4\}$ and $s\in \N$ that we call now weak $i$-bridges and whose lengths are given by $\c_i:=2\rho_i \delta_4= \rho_i/3$. As we said above, indices of weak chains satisfy the same properties as indices of classical chains (see Proposition \ref{LEMTHMmain2}). As a consequence, $I$ can be covered by weak subchains of the same form as the ones given in $\mathcal{L}$, in other words it can be covered by one of the following weak suchains: $\langle 234 \rangle$, $\langle 243\rangle$, $\langle 2434\rangle$, $\langle 324\rangle$, $\langle 3234\rangle$, $\langle 3243\rangle$, $\langle 32434\rangle$, $\langle 423\rangle$, $\langle 4234\rangle$, $\langle 4243\rangle$, $\langle 42434\rangle$, $\langle 342\rangle$, $\langle 3423\rangle$, $\langle 34234\rangle$, $\langle 3424\rangle$, $\langle 34243\rangle$, $\langle 342434\rangle$, $\langle 432\rangle$, $\langle 4323\rangle$,  $\langle 43234\rangle$, $\langle 4324\rangle$, $\langle 43243\rangle$, $\langle 432434\rangle$, $\langle 4342\rangle$, $\langle 43423\rangle$, $\langle 434234\rangle$, $\langle 43424\rangle$, $\langle 434243\rangle$, $\langle 4342434\rangle$. In any case, it can be covered from left to right by an interval $I_{3,4}$ given by an union of $3$ and $4$-weak bridges, a weak $2$-bridge $I_2$ and another interval $I_{3,4}'$ given by an union of $3$ and $4$-weak bridges. Then we have 
\[
\left| I \right| = \frac{2}{3} < \left| I_{3,4} \right| +   \left| I_{2} \right| +  \left| I_{3,4}' \right|
\]
with 
\[
 \left| I_{2} \right| < \frac{\rho_2}{3} \quad \mbox{and} \quad \left| I_{3,4} \right|, \left| I_{3,4}' \right| <  \frac{2\rho_3}{3} <  \frac{2\rho_2}{3}.
 \]
 We infer that $5\rho_2>2$. The two other inequalities of (\ref{10nov1}) follow from (\ref{12oct1}) and the formulas of the $\rho_i$'s. The equality $m_r(2)=M_r(2)$ for all $r=0,\ldots, L-1$ is a consequence of the fact that $2$ appears only once in each interval of the form $[r+2\delta_4,r+1]$. To prove the other inequalities, we note that (\ref{10nov1}) implies $\rho_3>7/40$ and $\rho_4>23/400$. Then, fix $r\in \{0,\ldots, L-1\}$ and set $m_i:=m_r(i), M_i:=M_i(r)$ for $i=3,4$. We have for every $i=3,4$,
 \[
 r+2\delta_4= r+\frac{1}{3} < h_i + \mathcal{C}^+_{m_i} < h_i + \mathcal{C}^-_{M_i} < r+1 
 \]
which gives
\[
\rho_i \left( M_i - m_i -\frac{1}{3} \right)=\mathcal{C}^-_{M_i} - \mathcal{C}^+_{m_i} < \frac{2}{3},
\]
so that
\[
M_i - m_i  \leq \frac{2}{3\rho_i} + \frac{1}{3}.
\] 
We infer the two remaining inequalities of (\ref{10nov2}) by using that $\rho_3>7/40$ and $\rho_4>23/400$. To prove (\ref{10nov3}), we fix $r$ in $\{0,\ldots, N-2\}$, set $M_i:=M_r(i)$ and $m_i':=m_{r+1}(i)$ for $i=2,3,4$ and note that we have 
\[
 r+2\delta_4= r+\frac{1}{3} < h_i + \mathcal{C}^+_{M_i} < h_i + \mathcal{C}^-_{m_i'} < r+2. 
 \]
 We infer that
  \[
\rho_i \left( m_i' - M_i -\frac{1}{3} \right)=\mathcal{C}^-_{m_i'} - \mathcal{C}^+_{M_i} < \frac{5}{3},
\]
which gives 
\[
m_i' - M_i  \leq \frac{5}{3\rho_i} + \frac{1}{3} \qquad \forall i=2,3,4.
\]
We conclude by the inequalities $\rho_2>2/5$, $\rho_3>7/40$ and $\rho_4>23/400$.
\end{proof}

 From now on, we say that a weak $L$-chain is admissible if it satisfies the properties (\ref{10nov1})-(\ref{10nov3}). The following lemma provides the exhaustive list of all weak $1$-chains. As explained below, its proof consists only in checking if some convex polytopes in $\R^6$ have non-empty interior. 
 
\begin{lemma}\label{LEMweak1chains}
The list of admissible weak $1$-chains is given by (if $a$ is an integer, then $\bar{a}$ stands for $a+1$): 
\begin{longtable}{|p{0.8cm}|p{5.15cm}|p{7.15cm}|}
\hline
\hspace{-0.15cm} 1 & $\langle 234\vert 000\rangle$ &  \\ 
\hline 
\hspace{-0.15cm} 2 & $\langle 243\vert 000\rangle$ & \\
\hline
\hspace{-0.15cm}  3 & $\langle324\vert 000\rangle$ &  \\
\hline
\hspace{-0.15cm}  4 & $\langle 342\vert 000\rangle$ &  \\
 \hline 
\hspace{-0.15cm}  5 & $\langle423\vert 000\rangle$ &  \\ 
\hline 
\hspace{-0.15cm}  6 & $\langle 432\vert 000\rangle$ & \\
\hline
\hspace{-0.15cm}  7 & $\langle 2434\vert 0001\rangle$ &  \\
\hline
\hspace{-0.15cm}  8 & $\langle 3234 \vert 0010\rangle $ &  \\
 \hline 
\hspace{-0.15cm}  9 & $\langle 3243 \vert 0001\rangle$ &  \\ 
\hline 
\hspace{-0.15cm}  10-12 & $\langle 4234\vert 000a\rangle$  & $a\in \{1,2,3\}$\\
\hline
\hspace{-0.15cm}  13-14 & $\langle 4243\vert 00a0\rangle$ & $a\in \{1,2\}$ \\
\hline
\hspace{-0.15cm}  15 & $\langle 4323\vert 0001\rangle$ &  \\
 \hline 
\hspace{-0.15cm}  16 & $\langle 3423\vert 0001\rangle $  & \\
\hline
\hspace{-0.15cm}  17-19 & $\langle4324\vert 000a\rangle$  & $a\in \{1,2,3\}$ \\
\hline
\hspace{-0.15cm}  20-21 & $\langle 3424\vert 000a\rangle$ & $a\in \{1,2\}$ \\
 \hline 
\hspace{-0.15cm}  22 & $\langle 4342\vert 0010\rangle$ &  \\
 \hline 
\hspace{-0.15cm}  23 & $\langle 32434\vert 00011\rangle$  & \\
\hline
\hspace{-0.15cm}  24 & $\langle 42434\vert 00102\rangle$  & \\
\hline
\hspace{-0.15cm}  25-27 & $\langle 34234\vert 0001a\rangle$ & $a\in \{1,2,3\}$ \\
 \hline 
\hspace{-0.15cm}  28-29 & $\langle 34243\vert 000a1\rangle$  & $a\in \{1,2\}$ \\
\hline
\hspace{-0.15cm}  30-32 & $\langle 43234\vert 0001a\rangle$   & $a\in \{2,3,4\}$\\
\hline
\hspace{-0.15cm}  33-35 & $\langle 43243\vert 000a1\rangle$ & $a\in \{1,2,3\}$ \\
 \hline 
\hspace{-0.15cm}  36 & $\langle 43423\vert 00101\rangle$  &  \\
\hline
\hspace{-0.15cm}  37 &  $\langle 43424\vert 00102\rangle$   & \\
\hline
\hspace{-0.15cm}  38 & $\langle 342434\vert 000213\rangle$ &  \\
 \hline 
\hspace{-0.15cm}  39-40 & $\langle 432434\vert 000a1\bar{a}\rangle$  & $a\in \{2,3\}$ \\
\hline
\hspace{-0.15cm}  41-42 &  $\langle 434234\vert 00101a\rangle$    & $a\in \{3,4\}$ \\
\hline
\hspace{-0.15cm}  43 & $\langle 434243\vert 001031\rangle$ &  \\
 \hline 
 \hspace{-0.15cm}  44 & $\langle 4342434\vert 0010314\rangle$ &  \\
 \hline 
\end{longtable}
\end{lemma}

\begin{proof}[Proof of Lemma \ref{LEMweak1chains}]
We verify the result with a Sage program \cite{SageMath}. Let us explain how to do it for example with the $1$-chain $\langle 2434\vert 000a\rangle$ where $a\in \N^*$ is an unknown parameter. We need to check what are the admissible chains of the form  $\langle 2434\vert 000a\rangle$, where we know thanks to (\ref {10nov2}), that the parameter $a$ can be taken in the range $\llbracket 1,3\rrbracket$. The corresponding set of inequalities is given by 
\[
 h_2 <  \frac{1}{3} < h_4 <  h_2+\frac{\rho_2}{3} < h_3  <  h_4+\frac{\rho_4}{3}  < h_4+a\rho_4<h_3+\frac{\rho_3}{3}<1 < h_4+a\rho_4+\frac{\rho_4}{3},
 \]
or equivalently by the system of inequalities with integer coefficients 
\[
\left\{ 
\begin{array}{rcl}
0 & < & 1-3 h_2 \\
0 & < & -1 + 3h_4 \\
0  & < & \rho_2 +3h_2 - 3h_4  \\
0  & < & -\rho_2 -3h_2 + 3h_3  \\
0  & < & \rho_4 -3h_3 + 3h_4  \\
0 & < & \rho_3 -3a\rho_4 +3h_3-3h_4\\
0 & < & 3 - \rho_3 -3 h_3\\
0  & < & -3 + (3a+1)\rho_4 +3h_4.   \\
\end{array}
\right.
\]
Thus, the admissible chains are those chains for which the interior of the corresponding convex polytope is empty. This verification can be made easily with Sage.
\end{proof}

Then, the following lemma provides the exhaustive list of all weak $2$-chains. Its proof consists again in checking the non-emptyness of convex (open) polytopes with Sage \cite{SageMath} (note that (\ref{10nov3}) is useful to restrict the set of parameters to consider for jumps).

\begin{lemma}\label{LEMweak2chains}
The list of admissible weak $2$-chains is given by (if $a$ is an integer, then $\bar{a}$ stands for $a+1$): 
\begin{longtable}{|p{1.2cm}|p{4.85cm}|p{7.05cm}|}
\hline
\hspace{-0.15cm} 1-5 & $\langle234 \hspace{-0.075cm}\shortmid \hspace{-0.075cm}234\vert 000 \m 11a\rangle$ & $a\in \llbracket 1,5\rrbracket$ \\ 
\hline \hline
\hspace{-0.15cm}  6-10 & $\langle243\m 234\vert 000\m 11a\rangle$ & $a\in \llbracket2,6\rrbracket$ \\
\hline
\hspace{-0.15cm}  11-15 & $\langle243\m 243\vert 000\m 1a1\rangle$ & $a\in  \llbracket1,5 \rrbracket$ \\
\hline
\hspace{-0.15cm}  16-18 & $\langle243\m 2434\vert 000\m 1a1\bar{a}\rangle$ & $a\in  \{2,3,4\}$ \\
 \hline \hline
\hspace{-0.15cm}  19 & $\langle324\m 243\vert 000\m 112\rangle$ &  \\
\hline
\hspace{-0.15cm}  20-24 & $\langle324\m 324\vert 000\m 11a\rangle$ & $a\in  \llbracket1,5 \rrbracket$ \\
\hline
\hspace{-0.15cm}  25 & $\langle324\m 423\vert 000\m 112\rangle$ & \\
\hline
\hspace{-0.15cm}  26 & $\langle324\m 3424 \vert 000\m 1213\rangle$ &  \\
\hline
\hspace{-0.15cm}  27 & $\langle324\m 43243 \vert 000\m 12123\rangle$ & \\
 \hline \hline
\hspace{-0.15cm}  28-32 & $\langle342\m 342 \vert 000\m 1a1\rangle$ & $a\in  \llbracket1,5 \rrbracket$ \\
\hline
\hspace{-0.15cm}  33 & $\langle342\m 423 \vert 000\m 112\rangle$ &  \\
\hline
\hspace{-0.15cm}  34 & $\langle342\m 4234 \vert 000\m 1122\rangle$ &  \\
\hline
\hspace{-0.15cm}  35 & $\langle342\m 3423 \vert 000\m 1112\rangle$ &  \\
\hline
\hspace{-0.15cm}  36-39 & $\langle342\m 3424 \vert 000\m 1a1b\rangle$ & $(a,b)\in \{(2,3), (3,4), (4,5), (5,7)\} $ \\
\hline
\hspace{-0.15cm}  40 & $\langle342\m 34243 \vert 000\m 12132\rangle$ &  \\
 \hline \hline
\hspace{-0.15cm}  41-47 & $\langle423\m 234 \vert 000\m 11a\rangle$ & $a\in  \llbracket2,8 \rrbracket$ \\
\hline
\hspace{-0.15cm}  48-52 & $\langle423\m 423 \vert 000\m a11\rangle$ & $a\in  \llbracket1,5 \rrbracket$ \\
\hline
\hspace{-0.15cm}  53-57 & $\langle423\m 4234 \vert 000\m a11b\rangle$ & $(a,b) \in \{(1,2), (2,3), (3,5), (4,7), (5,8)\}$ \\
\hline
\hspace{-0.15cm}  58 & $\langle423\m 4243 \vert 000\m 2132\rangle$ &  \\
\hline
\hspace{-0.15cm}  59 & $\langle423\m 43243 \vert 000\m 21132\rangle$ &  \\
 \hline
\hspace{-0.15cm}  60-66 & $\langle432\m 324 \vert 000\m 11a\rangle$ & $a\in  \llbracket2,8 \rrbracket$  \\
\hline
\hspace{-0.15cm}  67-71 & $\langle432\m 342 \vert 000\m 1a1\rangle$ & $a\in  \llbracket2,6 \rrbracket$  \\
\hline
\hspace{-0.15cm}  72-76 &$\langle432\m 432 \vert 000\m a11\rangle$ & $a\in  \llbracket1,5 \rrbracket$  \\
\hline
\hspace{-0.15cm}  77 & $\langle432\m 3243 \vert 000\m 1122\rangle$ &   \\
\hline
\hspace{-0.15cm}  78 & $\langle432\m 3423 \vert 000\m 1212\rangle$ &   \\
\hline
\hspace{-0.15cm}  79-83 & $\langle432\m 4324 \vert 000\m a11b\rangle$ & $(a,b) \in \{(1,2), (2,3), (3,5), (4,7), (5,8)\}$  \\
\hline
\hspace{-0.15cm}  84-86 & $\langle432\m 3424 \vert 000\m 1a1\bar{a}\rangle$ & $a \in  \{3,4,5\}$  \\
\hline
\hspace{-0.15cm}  87-89 & $\langle432\m 4342 \vert 000\m a1\bar{a}1\rangle$ & $a\in  \{2,3,4\}$  \\
\hline
\hspace{-0.15cm}  90 & $\langle432\m 43424 \vert 000\m 31415\rangle$ &   \\
 \hline \hline
 \hspace{-0.15cm}  91-93 & $\langle2434\m 234 \vert 0001\m 11a\rangle$ & $a \in  \{3,4,5\}$  \\
\hline
\hspace{-0.15cm}  94 & $\langle2434\m 2434 \vert 0001\m 11314\rangle$ &   \\
 \hline \hline
\hspace{-0.15cm}  95-97 & $\langle3234\m 34234 \vert 0010\m 2a13b\rangle$ & $(a,b) \in \{(1,2), (2,4), (3,6)\}$  \\
 \hline \hline
\hspace{-0.15cm}  98 & $\langle3243\m 234 \vert 0001\m 122\rangle$ &   \\
\hline
\hspace{-0.15cm}  99 & $\langle3243\m 243 \vert 0001\m112\rangle$ &   \\
\hline
\hspace{-0.15cm}  100 & $\langle3243\m 4234 \vert 0001\m 1122\rangle$ &   \\
\hline
\hspace{-0.15cm}  101 & $\langle3243\m 4324 \vert 0001\m 1212\rangle$ &   \\
\hline
\hspace{-0.15cm}  102-104 & $\langle3243\m 43243 \vert 0001\m a21b3\rangle$ &  $(a,b) \in \{(1,2), (2,4), (3,6)\}$ \\
 \hline \hline
\hspace{-0.15cm}  105-109 & $\langle4234\m 234 \vert 000a\m 11b\rangle$ & $(a,b) \in \{(1,2),(1,3),(2,5),(3,7),(3,8) \}$  \\
\hline
\hspace{-0.15cm}  110 & $\langle4234\m 3423 \vert 0001\m1212\rangle$ &   \\
\hline
\hspace{-0.15cm}  111 & $\langle4234\m 34234 \vert 0001\m 12123\rangle$ &   \\
\hline
\hspace{-0.15cm}  112 & $\langle4234\m 4234 \vert 0002\m 3115\rangle$ &   \\
\hline \hline
\hspace{-0.15cm}  113-115 & $\langle4243\m 234 \vert 0010\m 11a\rangle$ & $a\in  \{4,5,6\}$  \\
\hline
\hspace{-0.15cm}  116-119 & $\langle4243\m 243 \vert 00a0\m 1b1\rangle$ & $(a,b) \in \{(1,3),(1,4),(1,5),(2,7) \}$  \\
\hline
\hspace{-0.15cm}  120 & $\langle4243\m 423 \vert 0010\m 311\rangle$ &   \\
\hline
\hspace{-0.15cm}  121-123 & $\langle4243\m 2434 \vert 0010\m 1a1\bar{a}\rangle$ & $a\in  \{3,4,5 \}$  \\
\hline
\hspace{-0.15cm}  124 & $\langle4243\m 3234 \vert 0010\m 1124\rangle$ &   \\
\hline
\hspace{-0.15cm}  125 & $\langle4243\m 4234 \vert 0010\m 3115\rangle$ &   \\
\hline
\hspace{-0.15cm}  126 & $\langle4243\m 4243 \vert 0010\m 3141\rangle$ &   \\
\hline
\hspace{-0.15cm}  127 & $\langle4243\m 42434 \vert 0010\m 31415\rangle$ &   \\
\hline
\hspace{-0.15cm}  128 & $\langle4243\m 43243 \vert 0010\m 21132\rangle$ &   \\
 \hline \hline
\hspace{-0.15cm}  129 & $\langle4323\m 3424 \vert 0001\m 2314\rangle$ &   \\
\hline
\hspace{-0.15cm}  130-135 & $\langle4323\m 34234 \vert 0001\m 2a13b\rangle$ & \hspace{-0.25cm} $(a,b) \in \{(3,4), (4,6), (5,7), (6,9), (7,10), (8,11)\}$ \\
\hline
\hspace{-0.15cm}  136-138 & $\langle4323\m 43423 \vert 0001\m a2\bar{a}13\rangle$ &  $a \in  \{4,5,6 \}$ \\
\hline
\hspace{-0.15cm}  139-140 & $\langle4323\m 434234 \vert 0001\m a2\bar{a}13b\rangle$ &  $(a,b) \in \{(4,7), (6,10) \}$ \\
 \hline \hline
\hspace{-0.15cm}  141 & $\langle3423\m 234 \vert 0001\m 122\rangle$ &   \\
\hline
\hspace{-0.15cm}  142 & $\langle3423\m 4243 \vert 0001\m 2133\rangle$ &   \\
\hline
\hspace{-0.15cm}  143-146 & $\langle3423\m 43243 \vert 0001\m a21b3\rangle$ & $(a,b) \in \{(2,3), (3,5), (4,6), (5,8) \}$  \\
\hline
\hspace{-0.15cm}  147-148 & $\langle3423\m 43423 \vert 0001\m a2\bar{a}13\rangle$ & $a \in \{4,5\}$  \\
\hline
\hspace{-0.15cm}  149 & $\langle3423\m 434243 \vert 0001\m 526183\rangle$ &   \\
 \hline \hline
\hspace{-0.15cm}  150 & $\langle4324\m 243 \vert 0001\m 122\rangle$ &   \\
\hline
\hspace{-0.15cm}  151-155 & $\langle4324\m 324 \vert 000a\m 11b\rangle$ & $(a,b)  \in \{(1,2),(1,3), (2,5), (3,7),(3,8)\}$  \\
\hline
\hspace{-0.15cm}  156 & $\langle4324\m 3423 \vert 0001\m 1212\rangle$ &   \\
\hline
\hspace{-0.15cm}  157 & $\langle4324\m 34234 \vert 0001\m 12123\rangle$ &   \\
\hline
\hspace{-0.15cm}  158 & $\langle4324\m 4324 \vert 0002\m 3115\rangle$ &   \\
\hline
\hspace{-0.15cm}  159 & $\langle4324\m 3424 \vert 0002\m 1415\rangle$ &   \\
\hline
\hspace{-0.15cm}  160 & $\langle4324\m 43424 \vert 0002\m 31415\rangle$ &   \\
 \hline \hline
\hspace{-0.15cm}  161 & $\langle3424\m 324 \vert 0001\m 213\rangle$ &   \\
\hline
\hspace{-0.15cm}  162 & $\langle3424\m 3243 \vert 0001\m 2133\rangle$ &   \\
\hline
\hspace{-0.15cm}  163 & $\langle3424\m 3424 \vert 0001\m 1314\rangle$ &   \\
\hline
\hspace{-0.15cm}  164 & $\langle3424\m 43243 \vert 0001\m 22133\rangle$ &   \\
 \hline \hline
\hspace{-0.15cm}  165-167 & $\langle4342\m 342 \vert 0010\m 1a1\rangle$ & $a \in  \{3,4,5 \}$  \\
\hline
\hspace{-0.15cm}  168-170 & $\langle4342\m 3424 \vert 0010\m 1a1\bar{a}\rangle$ & $a \in  \{3,4,5 \}$  \\
\hline
\hspace{-0.15cm}  171 & $\langle4342\m 4342 \vert 0010\m 3141\rangle$ &   \\
\hline
\hspace{-0.15cm}  172 & $\langle4342\m 43424 \vert 0010\m 31415\rangle$ &   \\
 \hline \hline
\hspace{-0.15cm}  173-175 & $\langle32434\m 3234 \vert 00011\m 213a\rangle$ & $a \in  \{5,6,7 \}$  \\
\hline
\hspace{-0.15cm}  176-177 & $\langle32434\m 3243 \vert 00011\m 21a3\rangle$ & $a \in  \{5,6 \}$  \\
\hline
\hspace{-0.15cm}  178-179 & $\langle32434\m 32434 \vert 00011\m 21a3\bar{a}\rangle$ & $a \in \{5,6\}$  \\
\hline
\hspace{-0.15cm}  180-181 & $\langle32434\m 34234 \vert 00011\m 2a13b\rangle$ & $(a,b) \in \{(3,5), (4,7)\}$  \\
\hline
\hspace{-0.15cm}  182 & $\langle32434\m 43243 \vert 00011\m 32163\rangle$ &   \\
 \hline \hline
\hspace{-0.15cm}  183 & $\langle42434\m 234 \vert 00102\m 115\rangle$ &   \\
\hline
\hspace{-0.15cm}  184 & $\langle42434\m 2434 \vert 00102\m 1415\rangle$ &   \\
\hline
\hspace{-0.15cm}  185 & $\langle42434\m 4234 \vert 00102\m 3115\rangle$ &   \\
\hline \hline
\hspace{-0.15cm}  186 & $\langle34234\m 324 \vert 00011\m 213\rangle$ &   \\
\hline
\hspace{-0.15cm}  187 & $\langle34234\m 423 \vert 00011\m 213\rangle$ &   \\
\hline
\hspace{-0.15cm}  188-191 & $\langle34234\m 3243 \vert 0001a\m 21b3\rangle$ &  $(a,b) \in \{(1,3),(2,5),(2,6),(3,8)\}$  \\
\hline
\hspace{-0.15cm}  192 & $\langle34234\m 4243 \vert 00011\m 2133\rangle$ &   \\
\hline
\hspace{-0.15cm}  193-195 & $\langle34234\m 3423 \vert 0001a\m 2b13\rangle$ &   $(a,b) \in \{(1,2),(2,4),(3,6)\}$ \\
\hline
\hspace{-0.15cm}  196 & $\langle34234\m 3424 \vert 00011\m 2213\rangle$ &   \\
\hline
\hspace{-0.15cm}  197-199 & $\langle34234\m 34234 \vert 0001a\m 2b13c\rangle$ & $(a,b,c) \in \{(1,2,3),(2,4,6),(3,6,9)\}$  \\
\hline
\hspace{-0.15cm}  200-202 & $\langle34234\m 34243 \vert 0001a\m 2b1c3\rangle$ & $(a,b,c) \in \{(1,2,3),(2,4,5),(3,6,8)\}$  \\
\hline
\hspace{-0.15cm}  203-205 & $\langle34234\m 43243 \vert 0001a\m b21c3\rangle$ & $(a,b,c) \in \{(1,2,3),(2,3,5),(3,5,8)\}$   \\
\hline
\hspace{-0.15cm}  206 & $\langle34234\m 32434 \vert 00013\m 21839\rangle$ &   \\
\hline
\hspace{-0.15cm}  207 & $\langle34234\m 43423 \vert 00013\m 52613\rangle$ &   \\
\hline
\hspace{-0.15cm}  208 & $\langle34234\m 342434 \vert 00013\m 261839\rangle$ &   \\
\hline
\hspace{-0.15cm}  209 & $\langle34234\m 434243 \vert 00013\m 526183\rangle$ &   \\
 \hline \hline
\hspace{-0.15cm}  210 & $\langle34243\m 243 \vert 00011\m 132\rangle$ &   \\
\hline
\hspace{-0.15cm}  211-213 & $\langle34243\m 43243 \vert 000a1\m b21c3\rangle$ & $(a,b,c) \in \{(1,2,3),(1,3,5),(2,5,8)\}$   \\
 \hline \hline
\hspace{-0.15cm}  214-215 & $\langle43234\m 34234 \vert 0001a\m 2b13c\rangle$ & $(a,b,c) \in \{(3,5,7),(4,7,10)\}$  \\
 \hline \hline
\hspace{-0.15cm}  216-221 & $\langle43243\m 3234 \vert 000a1\m 213b\rangle$ &  \hspace{-0.25cm} $(a,b) \in \{(1,4),(2,6),(2,7),(3,9),(3,10),(3,11)\}$  \\
\hline
\hspace{-0.15cm}  222 & $\langle43243\m 4234 \vert 00011\m 2123\rangle$ &   \\
\hline
\hspace{-0.15cm}  223-225 & $\langle43243\m 4323 \vert 000a1\m b213\rangle$ & $(a,b) \in \{ (1,2), (2,4),(3,6) \}$  \\
\hline
\hspace{-0.15cm}  226 & $\langle43243\m 4324 \vert 00011\m 2213\rangle$ &   \\
\hline
\hspace{-0.15cm}  227-230 & $\langle43243\m 34234 \vert 000a1\m 2b13c\rangle$ & $(a,b,c) \in \{(1,3,4),(2,4,6),(2,5,7),(3,7,10)\}$  \\
\hline
\hspace{-0.15cm}  231-233 & $\langle43243\m 43243 \vert 000a1\m b21c3\rangle$ & $(a,b,c) \in \{(1,2,3),(2,4,6),(3,6,9)\}$  \\
\hline
\hspace{-0.15cm}  234-235 & $\langle43243\m 32434 \vert 000a1\m 21b3c\rangle$ & $(a,b,c) \in \{(2,6,7),(3,9,10)\}$  \\
\hline
\hspace{-0.15cm}  236-237 & $\langle43243\m 43234 \vert 000a1\m b213c\rangle$ & $(a,b,c) \in \{(2,4,7),(3,6,10)\}$  \\
\hline
\hspace{-0.15cm}  238-239 & $\langle43243\m 43423 \vert 000a1\m b2c13\rangle$ & $(a,b,c) \in \{(2,4,5),((3,6,7)\}$  \\
\hline
\hspace{-0.15cm}  240-241 & $\langle43243\m 432434 \vert 000a1\m b21c3d\rangle$ & $(a,b,c,d) \in \{(2,4,6,7),(3,6,9,10)\}$  \\
\hline
\hspace{-0.15cm}  242-243 & $\langle43243\m 434234 \vert 000a1\m b2c13d\rangle$ &  $(a,b,c,d) \in \{(2,4,5,7),(3,6,7,10)\}$ \\
\hline \hline
\hspace{-0.15cm}  244-247 & $\langle43423\m 32434 \vert 00101\m 21a3\bar{a}\rangle$ & $a\in  \llbracket7,10 \rrbracket$   \\
\hline
\hspace{-0.15cm}  248-249 & $\langle43423\m 34234 \vert 00101\m 2a13b\rangle$ & $(a,b) \in \{(5,7),(7,10)\}$  \\
\hline
\hspace{-0.15cm}  250 & $\langle43423\m 43243 \vert 00101\m 62193\rangle$ &    \\
\hline
\hspace{-0.15cm}  251-252 & $\langle43423\m 43423 \vert 00101\m a2\bar{a}13\rangle$ & $a\in \{5,6\}$  \\
\hline
\hspace{-0.15cm}  253 & $\langle43423\m 342434 \vert 00101\m 27193a\rangle$ & $a=10$  \\
\hline
\hspace{-0.15cm}  254 & $\langle43423\m 432434 \vert 00101\m 62193a\rangle$ & $a=10$  \\
\hline
\hspace{-0.15cm}  255 & $\langle43423\m 434234 \vert 00101\m 62713a\rangle$ & $a=10$  \\
\hline
\hspace{-0.15cm}  256 & $\langle43423\m 434243 \vert 00101\m 627193\rangle$ &   \\
\hline
\hspace{-0.15cm}  257 & $\langle43423\m 4342434 \vert 00101\m 627193a\rangle$ &  $a=10$ \\
\hline \hline
\hspace{-0.15cm}  258 & $\langle43424\m 3424 \vert 00102\m 1415\rangle$ &   \\
\hline \hline
\hspace{-0.15cm}  259 & $\langle342434\m 3243 \vert 000213\m 2183\rangle$ &   \\
\hline
\hspace{-0.15cm}  260 & $\langle342434\m 32434 \vert 000213\m 21839\rangle$ &   \\
\hline
\hspace{-0.15cm}  261 & $\langle342434\m 43243 \vert 000213\m 52183\rangle$ &   \\
\hline \hline
\hspace{-0.15cm}  262-263 & $\langle432434\m 3234 \vert 000a1b\m 213c\rangle$ & $(a,b,c) \in \{(2,3,7), (3,4,10) \}$  \\
\hline
\hspace{-0.15cm}  264-265 & $\langle432434\m 34234 \vert 000a1b\m 2c13d\rangle$ &  $(a,b,c,d) \in \{(2,3,5,7), (3,4,7,10) \}$  \\
\hline
\hspace{-0.15cm}  266 & $\langle432434\m 32434 \vert 000314\m 2193a\rangle$ &   $a=10$  \\
\hline \hline
\hspace{-0.15cm}  267-268 & $\langle434234\m 34234 \vert 00101a\m 2b13c\rangle$ &   $(a,b,c) \in \{(3,5,7), (4,7,10) \}$   \\
\hline
\hspace{-0.15cm}  269 & $\langle434234\m 32434 \vert 001014\m 2193a\rangle$ & $a=10$   \\
\hline
\hspace{-0.15cm}  270 & $\langle434234\m 342434 \vert 001014\m 27193a\rangle$ & $a=10$  \\
\hline \hline
\hspace{-0.15cm}  271 & $\langle434243\m 32434 \vert 001031\m 2193a\rangle$ & $a=10$  \\
\hline
\hspace{-0.15cm}  272 & $\langle434243\m 43243 \vert 001031\m 62193\rangle$ &   \\
\hline
\hspace{-0.15cm}  273 & $\langle434243\m 432434 \vert 001031\m 62193a\rangle$ & $a=10$  \\
\hline \hline
\hspace{-0.15cm}  274 & $\langle4342434\m 32434 \vert 0010314\m 2193a\rangle$ & $a=10$  \\
\hline
\end{longtable}
\end{lemma}

We observe that the index $2$ appears always twice in the admissible weak $2$-chains and that the second occurence of $2$ is always associated with a jump equal to $1$. This property allows us to restrict now our attention to weak $2$-chains associated with tuples $(h_1,h_2,h_3)\in \R^3$ satisfying the additional constraints
\begin{eqnarray}\label{EQ11nov1}
h_2 > \frac{1}{3} + \frac{\rho_2}{3}.
\end{eqnarray}
The following result can be checked easily with Sage \cite{SageMath}.

\begin{lemma}\label{LEMweak2chainsb2}
The list of admissible weak $2$-chains which can be extended into admissible weak $6$-chains and which satisfy the additional constraint (\ref{EQ11nov1}) is given  by (if $a$ is an integer, then $\bar{a}$ stands for $a+1$): 
\begin{longtable}{|p{1.2cm}|p{4.85cm}|p{7.05cm}|}
\hline
\hspace{-0.15cm} 1-5 & $\langle342 \hspace{-0.075cm}\shortmid \hspace{-0.075cm}342\vert 000 \m 1a1\rangle$ & $a\in \llbracket 1,5\rrbracket$ \\ 
\hline 
\hspace{-0.15cm}  6 & $\langle 342\m 3424\vert 000\m 1314\rangle$ &  \\
\hline \hline
\hspace{-0.15cm}  7-11 & $\langle 432\m 324\vert 000\m 11a\rangle$ & $a\in \{2,3,5,7,8\}$ \\
\hline
\hspace{-0.15cm}  12-14 & $\langle 432\m 342\vert 000\m 1a1\rangle$ & $a\in \llbracket 3,5\rrbracket$ \\
\hline
\hspace{-0.15cm}  15-19 & $\langle 432\m 432\vert 000\m a11\rangle$ & $a\in \llbracket 1,5\rrbracket$ \\
\hline
\hspace{-0.15cm}  20-24 & $\langle 432\m 4324\vert 000\m a11b\rangle$ & $(a,b) \in  \{(1,2),(2,3),(3,5),(4,7),(5,8)\}$ \\
\hline
\hspace{-0.15cm}  25 & $\langle 432\m 3424\vert 000\m 1415\rangle$ &  \\
\hline
\hspace{-0.15cm}  26-28 & $\langle 432\m 4342\vert 000\m a1\bar{a}1\rangle$ & $a\in \llbracket 2,4\rrbracket$ \\
\hline
\hspace{-0.15cm}  29 & $\langle 432\m 43424\vert 000\m 31415\rangle$ &  \\
\hline \hline
\hspace{-0.15cm}  30-34 & $\langle 4324 \m 324\vert 000a\m 11b\rangle$ & $(a,b) \in  \{(1,2),(1,3),(2,5),(3,7),(3,8)\}$  \\
\hline
\hspace{-0.15cm}  35 & $\langle 4324 \m 4324\vert 0002\m 3115\rangle$ &   \\
\hline
\hspace{-0.15cm}  36 & $\langle 4324 \m 3424\vert 0002\m 1415\rangle$ &   \\
\hline
\hspace{-0.15cm}  37 & $\langle 4324 \m 43424\vert 0002\m 31415\rangle$ &   \\
\hline \hline
\hspace{-0.15cm}  38 & $\langle 3424 \m 3424\vert 0001\m 1314\rangle$ &   \\
\hline \hline
\hspace{-0.15cm}  39-41 & $\langle 4342\m 342\vert 0010\m 1a1\rangle$ & $a\in \llbracket 3,5\rrbracket$ \\
\hline 
\hspace{-0.15cm}  42 & $\langle 4342\m 3424\vert 0010\m 1415\rangle$ &  \\
\hline 
\hspace{-0.15cm}  43 & $\langle 4342\m 4342\vert 0010\m 3141 \rangle$ &  \\
\hline 
\hspace{-0.15cm}  44 & $\langle 4342\m 43424\vert 0010\m 31415 \rangle$ &  \\
\hline \hline
\hspace{-0.15cm}  45 & $\langle 43424\m 3424\vert 00102\m 1415 \rangle$ &  \\
\hline \hline
\hspace{-0.15cm}  46 & $\langle 34234\m 34234\vert 00011\m 22133 \rangle$ &  \\
\hline \hline
\hspace{-0.15cm}  47 & $\langle 43243\m 43243\vert 00011\m 22133 \rangle$ &  \\
\hline 
\end{longtable}
\end{lemma}

We found 438 admissible weak $3$-chains of length 3, but we do not need all of them. We only need those starting from the weak  $2$-chains listed in Lemma \ref{LEMweak2chainsb2}. Given two weak $2$-chains  
\begin{multline*}
C_1= \langle i_1^0\cdots i_{\ell_0}^0 \m   i_1^{1}\cdots i_{\ell_{1}}^1 \vert s_1^0\cdots s_{\ell_0}^0 \m s_1^{1}\cdots s_{\ell_{1}}^{1} \rangle \\
 \mbox{and} \quad C_2= \langle \tilde{i}_1^0\cdots \tilde{i}_{\tilde{\ell}_0}^0 \m   \tilde{i}_1^{1}\cdots \tilde{i}_{\tilde{\ell}_{1}}^1 \vert \tilde{s}_1^0\cdots \tilde{s}_{\tilde{\ell}_0}^0 \m \tilde{s}_1^{1}\cdots \tilde{s}_{\tilde{\ell}_{1}}^{1} \rangle,
\end{multline*}
we say that $C_1$ is transfered to $C_2$ and we write 
\[
 \langle i_1^0\cdots i_{\ell_0}^0 \m   i_1^{1}\cdots i_{\ell_{1}}^1 \vert s_1^0\cdots s_{\ell_0}^0 \m s_1^{1}\cdots s_{\ell_{1}}^{1} \rangle \quad \longrightarrow \quad  \langle \tilde{i}_1^0\cdots \tilde{i}_{\tilde{\ell}_0}^0 \m   \tilde{i}_1^{1}\cdots \tilde{i}_{\tilde{\ell}_{1}}^1 \vert \tilde{s}_1^0\cdots \tilde{s}_{\tilde{\ell}_0}^0 \m \tilde{s}_1^{1}\cdots \tilde{s}_{\tilde{\ell}_{1}}^{1} \rangle
\]
if we have (the function $m_1$ is associated with the weak $2$-chain $C_1$) 
\[
\ell_1= \tilde{\ell}_0, \quad i_j^1 = \tilde{i}^0_j, \quad \tilde{s}_j^0 = s^1_j - m_1(i_j) \qquad \forall j=1, \ldots, \ell_1
\]
and if the weak $3$-chain given by
\[
 \langle i_1^0\cdots i_{\ell_0}^0 \m   i_1^{1}\cdots i_{\ell_{1}}^1  \m   \tilde{i}_1^{1}\cdots \tilde{i}_{\tilde{\ell}_{1}}^1 \vert s_1^0\cdots s_{\ell_0}^0 \m s_1^{1}\cdots s_{\ell_{1}}^{1} \m s_1^{2}\cdots s_{\ell_{2}}^{2} \rangle,
 \]
 where $\ell_2, s_1^{2}, \ldots, s_{\ell_{2}}^{2} \in \N$ are defined by
 \[
\ell_2 := \tilde{\ell}_1 \quad \mbox{and} \quad s_j^2 = m_1(i_j) + \tilde{s}_j^1 \qquad \forall j=1, \ldots, \ell_2,
 \]
 is admissible. Moreover, we say that a weak $2$-chain $C'$ is reachable from another  weak $2$-chain $C$ if there is a finite sequence of weak $2$-chains $C_1, \ldots, C_n$ such that 
 \[
 C \longrightarrow C_1\longrightarrow  \cdots \longrightarrow C_n \longrightarrow C'.
 \] 
 The following lemma can be checked with Sage \cite{SageMath} by verifying the admissibility of all weak $3$-chains reachable from the weak $2$-chains listed in Lemma \ref{LEMweak2chainsb2}.

\begin{lemma}\label{LEMweak2chainstransfers}
The list of all transfers from weak $2$-chains reachable from the weak $2$-chains from 1 to 45 in Lemma \ref{LEMweak2chainsb2} is given in the order in which they appear by (if $a$ is an integer, then $\bar{a}$ stands for $a+1$): 
\begin{longtable}{|p{1.2cm}|p{3.8cm}|p{0.25cm}|p{3.85cm}|p{3.2cm}|}
\hline
\hspace{-0.15cm} 1-5 & $\langle 342 \m 342 \vert 000 \m 1a1 \rangle$ &  $ \rightarrow$ & $\langle 342 \m 342 \vert 000 \m 1a1 \rangle$ & $a=1,2,3,4,5$\\ 
\hline
\hspace{-0.15cm} 6 & $\langle 342 \m 342 \vert 000 \m 111\rangle$ & $ \rightarrow$ & $\langle 342 \m 423 \vert 000 \m 112\rangle$ & \\ 
\hline
\hspace{-0.15cm} 7 & $\langle 342 \m 342 \vert 000 \m 111 \rangle$ & $ \rightarrow$ & $\langle 342 \m 3423 \vert 000 \m 1112\rangle$ & \\ 
\hline
\hspace{-0.15cm} 8-10 & $\langle 342 \m 342 \vert 000 \m 1a1\rangle$ & $ \rightarrow$ & $\langle 342 \m 3424 \vert 000 \m 1a1\bar{a}\rangle$ & $a=2,3,4$\\ 
\hline
\hspace{-0.15cm} 11 & $\langle 342 \m 342 \vert 000 \m 151\rangle$ & $ \rightarrow$ & $\langle 342 \m 3424 \vert 000 \m 1517\rangle$ & \\ 
\hline \hline
\hspace{-0.15cm} 12 & $\langle 342 \m 3424 \vert 000 \m 1314 \rangle$ & $ \rightarrow$ & $\langle 3424 \m 3424 \vert 0001 \m 1314 \rangle$ & \\ 
\hline \hline
\hspace{-0.15cm} 13 & $\langle 432 \m 324 \vert 000 \m 112\rangle$ & $ \rightarrow$ & $\langle 324 \m 243 \vert 000 \m 112\rangle$ & \\ 
\hline
\hspace{-0.15cm} 14-15 & $\langle 432 \m 324 \vert 000 \m 11a \rangle$ & $ \rightarrow$ & $\langle 324 \m 324 \vert 000 \m 11b\rangle$ & \hspace{-0.2cm} $(a,b)\in \{(2,1),(3,2)\}$ \\ 
\hline
\hspace{-0.15cm} 16 & $\langle 432 \m 324 \vert 000 \m 11a \rangle$ & $ \rightarrow$ & $\langle 324 \m 324 \vert 000 \m 11b\rangle$ & \hspace{-0.2cm} $(a,b)\in \{(5,3)\}$ \\ 
\hline
\hspace{-0.15cm} 17-18 & $\langle 432 \m 324 \vert 000 \m 11a \rangle$ & $ \rightarrow$ & $\langle 324 \m 324 \vert 000 \m 11b\rangle$ & \hspace{-0.2cm} $(a,b)\in \{(7,4),(8,5)\}$ \\ 
\hline
\hspace{-0.15cm} 19 & $\langle 432 \m 324 \vert 000 \m 115 \rangle$ & $ \rightarrow$ & $\langle 324 \m 3424 \vert 000 \m 1213\rangle$ & \\ 
\hline \hline
\hspace{-0.15cm} 20-21 & $\langle 432 \m 342 \vert 000 \m 1a1\rangle$ & $ \rightarrow$ & $\langle 342 \m 342 \vert 000 \m 1b1 \rangle$ &  \hspace{-0.2cm} $(a,b) \in \{(3,2),(4,3)\}$\\ 
\hline
\hspace{-0.15cm} 22 & $\langle 432 \m 342 \vert 000 \m 151\rangle$ & $ \rightarrow$ & $\langle 342 \m 342 \vert 000 \m 141 \rangle$ & \\ 
\hline
\hspace{-0.15cm} 23-24 & $\langle 432 \m 342 \vert 000 \m 1a1 \rangle$ & $ \rightarrow$ & $\langle 342 \m 3424 \vert 000 \m 1b1\bar{b}\rangle$ & \hspace{-0.2cm}$(a,b) \in \{(3,2),(4,3)\}$\\ 
\hline
\hspace{-0.15cm} 25 & $\langle 432 \m 342 \vert 000 \m 151 \rangle$ & $ \rightarrow$ & $\langle 342 \m 3424 \vert 000 \m 1415\rangle$ & \\ 
\hline \hline
\hspace{-0.15cm} 26 & $\langle 432 \m 432 \vert 000 \m 111 \rangle$ & $ \rightarrow$ & $\langle 432 \m 324 \vert 000 \m 112 \rangle$ & \\ 
\hline
\hspace{-0.15cm} 27-28 & $\langle 432 \m 432 \vert 000 \m 211 \rangle$ & $ \rightarrow$ & $\langle 432 \m 324 \vert 000 \m 11a \rangle$ & $a=3,4$\\ 
\hline
\hspace{-0.15cm} 29 & $\langle 432 \m 432 \vert 000 \m 311 \rangle$ & $ \rightarrow$ & $\langle 432 \m 324 \vert 000 \m 115 \rangle$ & \\ 
\hline
\hspace{-0.15cm} 30-31 & $\langle 432 \m 432 \vert 000 \m 411 \rangle$ & $ \rightarrow$ & $\langle 432 \m 324 \vert 000 \m 11a \rangle$ & $a=6,7$\\ 
\hline
\hspace{-0.15cm} 32 & $\langle 432 \m 432 \vert 000 \m 511 \rangle$ & $ \rightarrow$ & $\langle 432 \m 324 \vert 000 \m 118 \rangle$ & \\ 
\hline
\hspace{-0.15cm} 33 & $\langle 432 \m 432 \vert 000 \m 111\rangle$ & $ \rightarrow$ & $\langle 432 \m 342 \vert 000 \m 121\rangle$ & \\ 
\hline
\hspace{-0.15cm} 34-35 & $\langle 432 \m 432 \vert 000 \m 211\rangle$ & $ \rightarrow$ & $\langle 432 \m 342 \vert 000 \m 1a1\rangle$ & $a=2,3$\\ 
\hline
\hspace{-0.15cm} 36 & $\langle 432 \m 432 \vert 000 \m 311\rangle$ & $ \rightarrow$ & $\langle 432 \m 342 \vert 000 \m 141\rangle$ & \\ 
\hline
\hspace{-0.15cm} 37-38 & $\langle 432 \m 432 \vert 000 \m 411\rangle$ & $ \rightarrow$ & $\langle 432 \m 342 \vert 000 \m 1a1\rangle$ & $a=5,6$\\ 
\hline
\hspace{-0.15cm} 39 & $\langle 432 \m 432 \vert 000 \m 511\rangle$ & $ \rightarrow$ & $\langle 432 \m 342 \vert 000 \m 161\rangle$ & \\ 
\hline
\hspace{-0.15cm} 40-44 & $\langle 432 \m 432 \vert 000 \m a11\rangle$ & $ \rightarrow$ & $\langle 432 \m 432 \vert 000 \m a11\rangle$ & $a=1, 2,3,4,5$\\ 
\hline
\hspace{-0.15cm} 45-46 & $\langle 432 \m 432 \vert 000 \m a11\rangle$ & $ \rightarrow$ & $\langle 432 \m 4324 \vert 000 \m a11b\rangle$ &\hspace{-0.2cm} $(a,b)\in\{(1,2),(2,3)\}$\\ 
\hline
\hspace{-0.15cm} 47 & $\langle 432 \m 432 \vert 000 \m 311\rangle$ & $ \rightarrow$ & $\langle 432 \m 4324 \vert 000 \m 3115\rangle$ & \\ 
\hline
\hspace{-0.15cm} 48-49 & $\langle 432 \m 432 \vert 000 \m a11\rangle$ & $ \rightarrow$ & $\langle 432 \m 4324 \vert 000 \m a11b\rangle$ & \hspace{-0.2cm} $(a,b)\in \{(4,7),(5,8)\}$\\ 
\hline
\hspace{-0.15cm} 50-51 & $\langle 432 \m 432 \vert 000 \m a11 \rangle$ & $ \rightarrow$ & $\langle 432 \m 3424 \vert 000 \m 1b1\bar{b}\rangle$ &\hspace{-0.2cm} $(a,b)\in \{(2,3),(3,4)\}$\\ 
\hline
\hspace{-0.15cm} 52 & $\langle 432 \m 432 \vert 000 \m 411 \rangle$ & $ \rightarrow$ & $\langle 432 \m 3424 \vert 000 \m 1516\rangle$ & \\ 
\hline
\hspace{-0.15cm} 53-54 & $\langle 432 \m 432 \vert 000 \m a11\rangle$ & $ \rightarrow$ & $\langle 432 \m 4342 \vert 000 \m a1b1\rangle$ & \hspace{-0.2cm} $(a,b)\in \{(2,3),(3,4)\}$\\ 
\hline
\hspace{-0.15cm} 55 & $\langle 432 \m 432 \vert 000 \m 411\rangle$ & $ \rightarrow$ & $\langle 432 \m 4342 \vert 000 \m 4151\rangle$ & \\ 
\hline
\hspace{-0.15cm} 56 & $\langle 432 \m 432 \vert 000 \m 311\rangle$ & $ \rightarrow$ & $\langle 432 \m 43424 \vert 000 \m  31415\rangle$ & \\ 
\hline \hline
\hspace{-0.15cm} 57 & $\langle 432 \m 4324 \vert 000 \m 1112\rangle$ & $ \rightarrow$ & $\langle 4324 \m 324 \vert 0001 \m 112 \rangle$ & \\ 
\hline
\hspace{-0.15cm} 58 & $\langle 432 \m 4324 \vert 000 \m 2113\rangle$ & $ \rightarrow$ & $\langle 4324 \m 324 \vert 0001 \m 113 \rangle$ & \\ 
\hline
\hspace{-0.15cm} 59 & $\langle 432 \m 4324 \vert 000 \m 3115\rangle$ & $ \rightarrow$ & $\langle 4324 \m 324 \vert 0002 \m 115  \rangle$ &  \\ 
\hline
\hspace{-0.15cm} 60 & $\langle 432 \m 4324 \vert 000 \m 4117\rangle$ & $ \rightarrow$ & $\langle 4324 \m 324 \vert 0003 \m 117  \rangle$ &  \\ 
\hline
\hspace{-0.15cm} 61 & $\langle 432 \m 4324 \vert 000 \m 5118\rangle$ & $ \rightarrow$ & $\langle 4324 \m 324 \vert 0003 \m 118  \rangle$ &  \\ 
\hline
\hspace{-0.15cm} 62 & $\langle 432 \m 4324 \vert 000 \m 3115\rangle$ & $ \rightarrow$ & $\langle 4324 \m 4324 \vert 0002 \m 3115 \rangle$ & \\ 
\hline
\hspace{-0.15cm} 63 & $\langle 432 \m 4324 \vert 000 \m 3115 \rangle$ & $ \rightarrow$ & $\langle 4324 \m 3424 \vert 0002 \m 1415 \rangle$ & \\ 
\hline
\hspace{-0.15cm} 64 & $\langle 432 \m 4324 \vert 000 \m 3115 \rangle$ & $ \rightarrow$ & $\langle 4324 \m 43424 \vert 0002 \m 31415 \rangle$ & \\ 
\hline \hline
\hspace{-0.15cm} 65 & $\langle 432 \m 3424 \vert 000 \m 1415\rangle$ & $ \rightarrow$ & $\langle 3424 \m 3424 \vert 0001 \m   1314\rangle$ & \\ 
\hline \hline
\hspace{-0.15cm} 66-67 & $\langle 432 \m 4342 \vert 000 \m a1\bar{a}1\rangle$ & $ \rightarrow$ & $\langle 4342 \m 342 \vert 0010\m 1b1 \rangle$  & \hspace{-0.2cm} $(a,b) \in \{(2,3),(3,4)\}$\\ 
\hline
\hspace{-0.15cm} 68 & $\langle 432 \m 4342 \vert 000 \m 4151\rangle$ & $ \rightarrow$ & $\langle 4342 \m 342 \vert 0010\m 151 \rangle$  & \\ 
\hline
\hspace{-0.15cm} 69-70 & $\langle 432 \m 4342 \vert 000 \m a1\bar{a}1\rangle$ & $ \rightarrow$ & $\langle 4342 \m 3424 \vert 0010 \m 1b1\bar{b} \rangle$ & \hspace{-0.2cm} $(a,b) \in \{(2,3),(3,4)\}$\\ 
\hline
\hspace{-0.15cm} 71 & $\langle 432 \m 4342 \vert 000 \m 4151\rangle$ & $ \rightarrow$ & $\langle 4342 \m 3424 \vert 0010 \m 1516 \rangle$ & \\ 
\hline
\hspace{-0.15cm} 72 & $\langle 432 \m 4342 \vert 000 \m 3141 \rangle$ & $ \rightarrow$ & $\langle 4342 \m 4342 \vert 0010 \m 3141 \rangle$ & \\ 
\hline
\hspace{-0.15cm} 73 & $\langle 432 \m 4342 \vert 000 \m 3141\rangle$ & $ \rightarrow$ & $\langle 4342 \m 43424 \vert 0010 \m 31415 \rangle$ & \\ 
\hline \hline
\hspace{-0.15cm} 74 & $\langle 432 \m 43424 \vert 000 \m 31415 \rangle$ & $ \rightarrow$ & $\langle 43424 \m 3424 \vert 00102 \m 1415 \rangle$ & \\ 
\hline \hline
\hspace{-0.15cm} 75-76 & $\langle 4324 \m 324 \vert 0001 \m 11a \rangle$ & $ \rightarrow$ & $\langle 324 \m 324 \vert 000 \m 11 b\rangle$ & \hspace{-0.2cm} $(a,b)\in \{(2,1),(3,2)\}$\\ 
\hline
\hspace{-0.15cm} 77 & $\langle 4324 \m 324 \vert 0002 \m 115\rangle$ & $ \rightarrow$ & $\langle 324 \m 324 \vert 000 \m 113\rangle$ & \\ 
\hline
\hspace{-0.15cm} 78-79 & $\langle 4324 \m 324 \vert 0003\m 11a \rangle$ & $ \rightarrow$ & $\langle 324 \m 324 \vert 000 \m 11b\rangle$ & \hspace{-0.2cm} $(a,b)\in \{(7,4),(8,5)\}$\\ 
\hline
\hspace{-0.15cm} 80 & $\langle 4324 \m 324 \vert 0002 \m 115 \rangle$ & $ \rightarrow$ & $\langle 324 \m 3424 \vert 000 \m 1213 \rangle$ & \\ 
\hline \hline
\hspace{-0.15cm} 81 & $\langle 4324 \m 4324 \vert 0002 \m 3115 \rangle$ & $ \rightarrow$ & $\langle 4324 \m 324 \vert 0002 \m 115 \rangle$ & \\ 
\hline
\hspace{-0.15cm} 82 & $\langle 4324 \m 4324 \vert 0002 \m 3115 \rangle$ & $ \rightarrow$ & $\langle 4324 \m 4324 \vert 0002 \m 3115 \rangle$ & \\ 
\hline
\hspace{-0.15cm} 83 & $\langle 4324 \m 4324 \vert 0002 \m 3115 \rangle$ & $ \rightarrow$ & $\langle 4324 \m 3424 \vert 0002 \m 1415  \rangle$ & \\ 
\hline
\hspace{-0.15cm} 84 & $\langle 4324 \m 4324 \vert 0002 \m 3115 \rangle$ & $ \rightarrow$ & $\langle 4324 \m 43424 \vert 0002 \m 31415 \rangle$ & \\  
\hline \hline
\hspace{-0.15cm} 85 & $\langle 4324 \m 3424 \vert 0002 \m 1415 \rangle$ & $ \rightarrow$ & $\langle 3424 \m 3424 \vert 0001 \m 1314 \rangle$ & \\ 
\hline \hline
\hspace{-0.15cm} 86 & $\langle 4324 \m 43424 \vert 0002 \m 31415 \rangle$ & $ \rightarrow$ & $\langle 43424 \m 3424 \vert 00102 \m 1415 \rangle$ &  \\ 
\hline \hline
\hspace{-0.15cm} 87 & $\langle 3424 \m 3424 \vert 0001 \m 1314 \rangle$ & $ \rightarrow$ & $\langle 3424 \m 3424 \vert 0001\m 1314 \rangle$ & \\
\hline \hline 
\hspace{-0.15cm} 88-89 & $\langle 4342 \m 342 \vert 0010\m 1a1 \rangle$ & $ \rightarrow$ & $\langle 342 \m 342 \vert 000 \m 1b1 \rangle$ & \hspace{-0.2cm} $(a,b) \in \{(3,2),(4,3)\}$ \\ 
\hline
\hspace{-0.15cm} 90 & $\langle 4342 \m 342 \vert 0010\m 151 \rangle$ & $ \rightarrow$ & $\langle 342 \m 342 \vert 000 \m 141 \rangle$ &  \\ 
\hline
\hspace{-0.15cm} 91-92 & $\langle 4342 \m 342 \vert 0010 \m 1a1 \rangle$ & $ \rightarrow$ & $\langle 342 \m 3424 \vert 000 \m 1b1\bar{b} \rangle$ & \hspace{-0.2cm} $(a,b) \in \{(3,2),(4,3)\}$\\ 
\hline
\hspace{-0.15cm} 93 & $\langle 4342 \m 342 \vert 0010 \m 151 \rangle$ & $ \rightarrow$ & $\langle 342 \m 3424 \vert 000 \m 1415 \rangle$ & \\ 
\hline \hline
\hspace{-0.15cm} 94 & $\langle 4342 \m 3424 \vert 0010 \m 1415 \rangle$ & $ \rightarrow$ & $\langle 3424 \m 3424 \vert 0001 \m 1314 \rangle$ & \\ 
\hline \hline
\hspace{-0.15cm} 95-96 & $\langle 4342 \m 3424 \vert 0010 \m 1a1\bar{a} \rangle$ & $ \rightarrow$ & no chain & $a=3,5$\\ 
\hline \hline
\hspace{-0.15cm} 97 & $\langle 4342 \m 4342 \vert 0010 \m 3141 \rangle$ & $ \rightarrow$ & $\langle 4342 \m 342 \vert 0010 \m 141 \rangle$ & \\ 
\hline
\hspace{-0.15cm} 98 & $\langle 4342 \m 4342 \vert 0010 \m 3141 \rangle$ & $ \rightarrow$ & $\langle 4342 \m 3424 \vert 0010 \m 1415 \rangle$ & \\ 
\hline
\hspace{-0.15cm} 99 & $\langle 4342 \m 4342 \vert 0010 \m 3141 \rangle$ & $ \rightarrow$ & $\langle 4342 \m 4342 \vert 0010 \m 3141 \rangle$ & \\ 
\hline
\hspace{-0.15cm} 100 & $\langle 4342 \m 4342 \vert 0010 \m 3141 \rangle$ & $ \rightarrow$ & $\langle 4342 \m 43424 \vert 0010 \m 31415 \rangle$ & \\ 
\hline \hline
\hspace{-0.15cm} 101 & $\langle 4342 \m 43424 \vert 0010 \m 31415 \rangle$ & $ \rightarrow$ & $\langle 43424 \m 3424 \m 00102 \m 1415 \rangle$ & \\ 
\hline \hline
\hspace{-0.15cm} 102 & $\langle 43424 \m 3424 \vert 00102 \m 1415 \rangle$ & $ \rightarrow$ & $\langle 3424 \m 3424 \vert 0001 \m 1314 \rangle$ & \\ 
\hline \hline
\hspace{-0.15cm} 103 & $\langle 342 \m 423 \vert 000 \m 112 \rangle$ & $ \rightarrow$ & $\langle 423 \m 234 \vert 000 \m 112\rangle$ & \\ 
\hline \hline
\hspace{-0.15cm} 104 & $\langle 342 \m 3423  \vert 000 \m 1112\rangle$ & $ \rightarrow$ & no chain & \\ 
\hline \hline
\hspace{-0.15cm} 105-106 & $\langle 342 \m 3424 \vert 000 \m 1a1\bar{a}\rangle$ & $ \rightarrow$ & no chain & $a=2,4$\\ 
\hline
\hspace{-0.15cm} 107 & $\langle 342 \m 3424 \vert 000 \m 1517\rangle$ & $ \rightarrow$ & no chain & \\ 
\hline
\hline \hline
\hspace{-0.15cm} 108 & $\langle 324 \m 243 \vert 000 \m 112 \rangle$ & $ \rightarrow$ & $\langle 243 \m 243 \vert 000 \m 111\rangle$ &  \\ 
\hline \hline
\hspace{-0.15cm} 109-113 & $\langle 324 \m 324 \vert 000 \m 11a\rangle$ & $ \rightarrow$ & $\langle 324 \m 324 \vert 000 \m 11a\rangle$ &  \hspace{-0.2cm}  $a=1,2,3,4,5$ \\ 
\hline 
\hspace{-0.15cm} 114 & $\langle 324 \m 324 \vert 000 \m 113 \rangle$ & $ \rightarrow$ & $\langle 324 \m 3424 \vert 000 \m 1213 \rangle$ & \\ 
\hline \hline
\hspace{-0.15cm} 115 & $\langle 324 \m 3424 \vert 000 \m 1213\rangle$ & $ \rightarrow$ & $\langle 3424 \m 3424 \vert 0001 \m 1314 \rangle$ & \\ 
\hline \hline
\hspace{-0.15cm} 116-117 & $\langle 432 \m 324 \vert 000 \m 11a \rangle$ & $ \rightarrow$ & $\langle 324 \m 324 \vert 000 \m 11b \rangle$ & $(a,b)=(4,2), (6,4)$\\ 
\hline \hline
\hspace{-0.15cm} 118-119 & $\langle 432 \m 342 \vert 000 \m 121 \rangle$ & $ \rightarrow$ & $\langle 342 \m 342 \vert 000 \m 1b1 \rangle$ & $b=1,2$\\ 
\hline
\hspace{-0.15cm} 120-121 & $\langle 432 \m 342 \vert 000 \m 161 \rangle$ & $ \rightarrow$ & $\langle 342 \m 342 \vert 000 \m 1b1 \rangle$ & $b=4,5$\\ 
\hline \hline
\hspace{-0.15cm} 122-124 & $\langle 432 \m 3424 \vert 000 \m 1a1\bar{a} \rangle$ & $ \rightarrow$ & no chain & $a=2,3,5$\\ 
\hline \hline
\hspace{-0.15cm} 125 & $\langle 423 \m 234 \vert 000 \m 112 \rangle$ & $ \rightarrow$ & $\langle 234 \m 234 \vert 000 \m 111\rangle$ &   \\ 
\hline \hline
\hspace{-0.15cm} 126 & $\langle 243 \m 243  \vert 000 \m 111 \rangle$ & $ \rightarrow$ & $\langle 243 \m 234 \vert 000 \m 112\rangle$ & \\ 
\hline 
\hspace{-0.15cm} 127 & $\langle 243 \m 243 \vert 000 \m 111\rangle$ & $ \rightarrow$ & $\langle 243 \m 243 \vert 000 \m 111\rangle$ &    \\ 
\hline \hline
\hspace{-0.15cm} 128 & $\langle 234 \m 234 \vert 000 \m 111 \rangle$ & $ \rightarrow$ & $\langle 234 \m 234 \vert 000 \m 111 \rangle$ &  \hspace{-0.2cm}  \\ 
\hline \hline
\hspace{-0.15cm} 129-130 & $\langle 243 \m 234 \vert 000 \m 112\rangle$ & $ \rightarrow$ & $\langle 234 \m 234 \vert 000 \m 11a \rangle$ & $a\in \{1,2\}$ \\
\hline \hline 
\hspace{-0.15cm} 131-132 & $\langle 234 \m 234 \vert 000 \m 112\rangle$ & $ \rightarrow$ & $\langle 234 \m 234 \vert 000 \m 11a \rangle$ & $a\in \{1,2\}$ \\  
\hline 
\end{longtable}
\end{lemma}

We notice that almost all weak $2$-chains appearing in the table of Lemma \ref{LEMweak2chainstransfers} have a $3$-jump equal to $1$. As a matter of fact, in the above table, only the chains in the set 
\[
\mathcal{A} := \Bigl\{ \langle 342 \m 423 \vert 000 \m 112 \rangle \langle 342 \m 3423  \vert 000 \m 1112\rangle, \langle 324 \m 243 \vert 000 \m 112 \rangle\Bigr\},
\]
have a a $3$-jump strictly larger than $1$. The following lemma will allow us to rule out the chains of $\mathcal{A}$ in our final reasoning, it can be checked with Sage \cite{SageMath}.

\begin{lemma}\label{LEMweak2chainsrulesout}
There are no admissible weak $6$-chains which contain one of the weak $2$-chains of $\mathcal{A}$ as a subchain. 
\end{lemma}

Finally the following lemma easily checkable with Sage will be useful to treat the case of chains starting with subchains corresponding to the weak $2$-chains 46 and 47 of Lemma  \ref{LEMweak2chainsb2}.

\begin{lemma}\label{LEMweak2chainsrulesoutotherchains}
If there is an admissible weak $6$-chain starting with the weak $2$-chains 46 or 47 of Lemma  \ref{LEMweak2chainsb2}, then we have $\rho_3\leq 3\rho_2/5$. 
\end{lemma}

We are now ready to finalize the proof of the covering property at infinity.

\paragraph{End of the proof.} Define the two integers $a, b \geq 1$ by (note that $z_1>1$)
\[
a = \left\lceil \frac{z_1-1}{6} \right\rceil \quad \mbox{and} \quad b= \left\lfloor \frac{5z_1+1}{6} \right\rfloor.
\] 
By construction, $a$ is the least integer $k$ such that $\mathcal{K}_k^-/z_1=\mathcal{B}_k^+/z_1\geq \delta_4=1/6$ and $b$ is the largest integer $k$ such that $\mathcal{K}_{k-1}^+/z_1=\mathcal{B}_k^-/z_1\leq 1-\delta_4=5/6$. We aim to show the following result:

\begin{lemma}\label{LEM10nov}
If $z_1\geq 47/5$, then we have $b-a < 2a+2$. 
\end{lemma}

\begin{proof}[Proof of Lemma \ref{LEM10nov}]
We suppose for contradiction that $z_1\geq 47/5$ and $L:=b-a \geq 2a+2$. Then we have $a\geq \lceil 21/15 \rceil=2$,  $b\geq 8$, $b-a\geq 2a+2\geq 6$ and the set $\mathcal{K}_0$ contains the $b-a$ consecutive $1$-kwais 
\[
\mathcal{K}_{a}/z_1,  \quad \mathcal{K}_{a+1}/z_1, \quad \cdots, \quad  \mathcal{K}_{b-1}/z_1.
\]
By Lemma \ref{LEM20oct0}, the tuple $(\rho_2,\rho_3,\rho_4)\in \R^3$  defined by $\rho_i := z_1/z_i$ for $i=2,3,4$ 
admits a weak $L$-chain of the form (\ref{weakchain}) such that 
\begin{eqnarray}\label{11nov1}
\left[1 i_1^0\cdots i_{\ell_0}^0 1 \cdots 1 i_1^{L-1}\cdots i_{\ell_{L-1}}^{L-1} 1 \Big\vert s^0 \tilde{s}_1^0\cdots \tilde{s}_{\ell_0}^0s^1 \cdots \tilde{s}^{L-1} \tilde{s}_1^{L-1}\cdots \tilde{s}_{\ell_{L-1}}^{L-1} s^L\right]
\end{eqnarray}
is a $1$-subchain associated with $(z_1,z_2,z_3,z_4)$ of length $L$ (we use the notations of Lemma \ref{LEM20oct0}). Recall that $\bar{j}(2)$ is given by 
\[
\bar{j}(2)  =  \min \Bigl\{ j \in \{1,\ldots, \ell\} \, \vert \, i_j=i \, \mbox{ and } \, \mathcal{K}_a^-/z_1 < \mathcal{B}_{k_{j}}^+/z_i \Bigr\} 
\]
and set
\[
\bar{B}_2 := \mathcal{B}_{k_{\bar{j}(2)}}^-/z_2.
\]
We distinguish two cases:\\

\noindent Case 1: $\bar{B}_2 \leq  \mathcal{B}_{a}^+/z_1  +\b_2$.\\
We have
 \[
 \mathcal{B}_{a}^-/z_1  =   \frac{6a-1}{6z_1} < \bar{B}_2= \mathcal{B}_{k_{\bar{j}(2)}}^-/z_2 = \frac{6k_{\bar{j}(2)}-1}{6z_2}  <  \mathcal{B}_{a}^+/z_1  +\b_2 =  \frac{6a+1}{6z_1} +\frac{1}{3z_2},
\]
which gives, since $z_2>z_1$, $k_{\bar{j}(2)} \geq a+1$ and as a consequence
\[
z_2 > \left(\frac{6k_{\bar{j}(2)}-3}{6a+1}\right) z_1 \geq \left( \frac{6a+3}{6a+1}\right) z_1 = (1+\kappa) z_1 \quad \mbox{with} \quad \kappa := \frac{2}{6a+1}.
\]
Define the sequence $\{u_k\}_{k\in \N}$ by 
$$
u_k := \mathcal{B}_{k_{\bar{j}(2)}+k}^+/z_2 - \mathcal{B}_{a+k}^+/z_1 =  \mathcal{B}_{k_{\bar{j}(2)}}^+/z_2 - \mathcal{B}_{a}^+/z_1 + k \left(\frac{1}{z_2}-\frac{1}{z_1}\right) \qquad \forall k\in \N.
$$
By construction, we have (because $\bar{B}_2 <  \mathcal{B}_{a}^+/z_1  +\b_2$ with $\b_2=1/(3z_2)$)
$$
u_0 =  \mathcal{B}_{k_{\bar{j}(2)}}^+/z_2 - \mathcal{B}_{a}^+/z_1 = \bar{B}_2 + \b_2 - \mathcal{B}_{a}^+/z_1  \in (0,2/(3z_2))
$$
and
$$
u_{k+1}-u_k = \frac{1}{z_2}- \frac{1}{z_1} <0 \qquad \forall k\in \N. 
$$
Therefore, $\{u_k\}_{k\in \N}$ is decreasing and there is $k\in \N^*$ such that 
\[
u_k<0 \quad \mbox{and} \quad 
k\leq  \left\lceil \frac{2/(3z_2)}{1/z_1-1/z_2} \right\rceil \leq  \left\lceil \frac{2}{3\kappa} \right\rceil=  \left\lceil 2a + \frac{1}{3}  \right\rceil = 2a+1.
\]
By Lemma \ref{LEMweak2chains}, all $2$-jumps of weak $2$-chains are equal to $1$, so our Lemma \ref{LEM20oct0} guarantees that the same result holds for the $1$-subchain (\ref{11nov1}). As a consequence, if all the $1$-kwais from $\mathcal{K}_{a}/z_1$ to $\mathcal{K}_{a+k}/z_1$ are contained in $\mathcal{K}_0$, then we infer that the $1$-kwai $\mathcal{K}_{a+k}/z_1$ satisfies
\[
\mathcal{K}_{a+k}/z_1 \subset \mathcal{B}_{k_{\bar{j}(2)+k} }/z_2 \cup \mathcal{B}/z_3 \cup \mathcal{B}/z_4 \quad \mbox{with} \quad \mathcal{B}_{k_{\bar{j}(2)+k} }^+/z_2<\mathcal{B}_{a+k}^+/z_1 = \mathcal{K}_{a+k}^-/z_1,
\]
which is impossible. The above argument applies provided  $\mathcal{K}_{a+k}/z_1\subset \mathcal{K}_0$ with $1\leq k\leq 2a+1$ which is satisfied if  $b-a \geq 2a+2$. So we obtain a contradiction.\\

\noindent Case 2: $\bar{B}_2 >  \mathcal{B}_{a}^+/z_1  +\b_2$.\\
This assumption implies that the weak $L$-chain associated with (\ref{11nov1}) by Lemma \ref{LEM20oct0} satisfies the constraints (\ref{EQ11nov1}). Thus, the weak $2$-chain associated with the subchain $[1 i_1^0\cdots i_{\ell_0}^0  1 i_1^{1}\cdots i_{\ell_{1}}^{1} 1]$  is one of the weak $2$-chains from $1$ to $47$ listed in Lemma \ref{LEMweak2chainsb2}. Let us now distinguish two subcases:\\

\noindent Subcase 2.1: The weak $2$-chain associated with the subchain $[1 i_1^0\cdots i_{\ell_0}^0  1 i_1^{1}\cdots i_{\ell_{1}}^{1} 1]$ is not one of the weak $2$-chains 46 or 47 of Lemma \ref{LEMweak2chainsb2}.\\
Let us label all the weak $2$-chains reachable from one of the weak $2$-chains $1$ to $45$ Lemma \ref{LEMweak2chainsb2} of as follows:

\begin{longtable}{|p{0.3cm}|p{4cm}|p{6cm}|}
\hline 
$A$ & $\langle 342 \m 342 \vert 000 \m 1a1 \rangle$ &  $a=1,2,3,4,5$\\
\hline
$B$ & $\langle 342 \m 3424 \vert 000 \m 1314 \rangle$ & \\
\hline
$C$ & $\langle 432 \m 324 \vert 000 \m 11a \rangle$ & $a=2,3,5,7,8$\\
\hline
$D$ & $\langle 432 \m 342 \vert 000 \m 1a1 \rangle$ & $a=3,4,5$\\
\hline
$E$ & $\langle 432 \m 432 \vert 000 \m a11\rangle$ & $a=1,2,3,4,5$\\
\hline
$F$ & $\langle 432 \m 4324 \vert 000 \m a11b\rangle$ & $(a,b)=(1,2), (2,3), (3,5),(4,7), (5,8)$\\
\hline
$G$ & $\langle 432 \m 3424 \vert 000 \m 1415\rangle$ & \\
\hline
$H$ & $\langle 432 \m 4342 \vert 000 \m a1\bar{a}1\rangle$ & $a=2,3,4$\\
\hline
$I$ & $\langle 432 \m 43424 \vert 000 \m 31415 \rangle$ & \\
\hline
$J$ & $\langle 4324 \m 324 \vert 000a \m 11b \rangle$ & $(a,b)=(1,2), (1,3), (2,5),(3,7), (3,8)$\\
\hline
$K$ & $\langle 4324 \m 4324 \vert 0002 \m 3115 \rangle$ & \\
\hline
$L$ & $\langle 4324 \m 3424 \vert 0002 \m 1415 \rangle$ & \\
\hline
$M$ & $\langle 4324 \m 43424 \vert 0002 \m 31415 \rangle$ & \\
\hline
$N$ & $\langle 3424 \m 3424 \vert 0001 \m 1314 \rangle$ & \\
\hline
$O$ & $\langle 4342 \m 342 \vert 0010\m 1a1 \rangle$ & $a=3,4,5$\\
\hline
$P$ & $\langle 4342 \m 3424 \vert 0010 \m 1415 \rangle$ & \\
\hline
$Q$ & $\langle 4342 \m 3424 \vert 0010 \m 1a1\bar{a} \rangle$ &$a=3,5$\\
\hline
$R$ & $\langle 4342 \m 4342 \vert 0010 \m 3141 \rangle$ &\\
\hline
$S$ & $\langle 4342 \m 43424 \vert 0010 \m 31415 \rangle$ &\\
\hline
$T$ & $\langle 43424 \m 3424 \vert 00102 \m 1415 \rangle$ &\\
\hline
$U$ & $\langle 342 \m 423 \vert 000 \m 112 \rangle$ &\\
\hline
$V$ & $\langle 342 \m 3423 \vert 000 \m 1112 \rangle$ &\\
\hline
$W$ & $\langle 342 \m 3424 \vert 0001 \m 1a1b \rangle$ & $(a,b)=(2,3),(4,5),(5,7)$\\
\hline
$X$ & $\langle 324 \m 243 \vert 000 \m 112 \rangle$ & \\
\hline
$Y$ & $\langle 324 \m 324 \vert 000 \m 11a\rangle$ & $a=1,2,3,4,5$\\
\hline
$Z$ & $\langle 324 \m 3424 \vert 000 \m 1213\rangle$ & \\
\hline
$\Gamma$ & $\langle 432 \m 324 \vert 000 \m 11a \rangle$ & $a=4,6$\\
\hline
$\Delta$ & $\langle 432 \m 342 \vert 000 \m 1a1 \rangle$ & $a=2,6$\\
\hline
$\Theta$ & $\langle 432 \m 3424 \vert 000 \m 1a1\bar{a}\rangle$ & $a=2,3,5$\\
\hline
$\Lambda$ & $\langle 423 \m 234 \vert 000 \m 112\rangle$ & \\
\hline
$\Xi$ & $\langle 243 \m 243 \vert 000 \m 111\rangle$ &\\
\hline
$\Pi$ & $\langle 234 \m 234 \vert 000 \m 11a \rangle$ & $a=1,2$\\
\hline
$\Sigma$ & $\langle 243 \m 234 \vert 000 \m 112 \rangle$ & \\
\hline
\end{longtable}

Lemma \ref{LEMweak2chainstransfers} shows that the graph of all transfers between the above weak $2$-chains, where we draw an arrow from a weak $2$-chain $C$ to a weak $2$-chain $C'$ if $C$ can be transfered to $C'$, is invariant (see Figure \ref{figgraph}). 

\vspace{0.2cm}
\begin{figure}[H]
\centering
\begin{tikzpicture}[scale=0.75]
\tikzstyle{sommetini}=[circle,draw,thick,fill=yellow]
\tikzstyle{sommet}=[circle,draw,thick]
\tikzstyle{sommetpb}=[circle,draw,thick,red,fill=yellow]
\tikzstyle{sommetpb2}=[circle,draw,thick,red]
\node[sommetini] (A) at (0,0) {A};
\node[sommetini] (B) at (2,0) {B};
\node[sommetini] (C) at (4,0) {C};
\node[sommetini] (D) at (6,0) {D};
\node[sommetini] (E) at (8,0) {E};
\node[sommetini] (F) at (10,0) {F};
\node[sommetini] (G) at (12,0) {G};
\node[sommetini] (H) at (1,-2) {H};
\node[sommetini] (I) at (3,-2) {I};
\node[sommetini] (J) at (5,-2) {J};
\node[sommetini] (K) at (7,-2) {K};
\node[sommetini] (L) at (9,-2) {L};
\node[sommetini] (M) at (11,-2) {M};
\node[sommetini] (N) at (1,-4) {N};
\node[sommetini] (O) at (3,-4) {O};
\node[sommetini] (P) at (5,-4) {P};
\node[sommetini] (Q) at (7,-4) {Q};
\node[sommetini] (R) at (9,-4) {R};
\node[sommetini] (S) at (11,-4) {S};
\node[sommet] (T) at (0,-6) {T};
\node[sommetpb2] (U) at (2,-6) {U};
\node[sommetpb2] (V) at (4,-6) {V};
\node[sommet] (W) at (6,-6) {W};
\node[sommetpb2] (X) at (8,-6) {X};
\node[sommet] (Y) at (10,-6) {Y};
\node[sommet] (Z) at (12,-6) {Z};
\node[sommet] (Gamma) at (0,-8) {$\Gamma$};
\node[sommet] (Delta) at (2,-8) {$\Delta$};
\node[sommet] (Theta) at (4,-8) {$\Theta$};
\node[sommet] (Lambda) at (6,-8) {$\Lambda$};
\node[sommet] (Xi) at (8,-8) {$\Xi$};
\node[sommet] (Pi) at (10,-8) {$\Pi$};
\node[sommet] (Sigma) at (12,-8) {$\Sigma$};
\draw (A) edge[loop above] (A);
\draw[->] (A) to (B);
\draw[->] (A) to[bend right] (U);
\draw[->] (A) to[bend left=22] (V);
\draw[->] (A) to[bend left=38] (W);
\draw[->] (B) to[bend right=50] (N);
\draw[->] (C) to[bend right=25] (X);
\draw[->] (C) to[bend left=35] (Y);
\draw[->] (C) to[bend right=80] (Z);
\draw[->] (D) to[bend right] (A);
\draw[->] (D) to[bend right] (B);
\draw[->] (D) to (W);
\draw[->] (E) to[bend right] (C);
\draw[->] (E) to[bend right] (D);
\draw (E) edge[loop above] (E);
\draw[->] (E) to (F);
\draw[->] (E) to [bend left](G);
\draw[->] (E) to (H);
\draw[->] (E) to (I);
\draw[->] (E) to[bend left=23] (Delta);
\draw[->] (E) to[bend left=45] (Theta);
\draw[->] (F) to (J);
\draw[->] (F) to (K);
\draw[->] (F) to (L);
\draw[->] (F) to (M);
\draw[->] (G) to [bend right=43] (N);
\draw[->] (H) to (D);
\draw (H) edge[loop above] (H);
\draw[->] (H) to (P);
\draw[->] (H) to (Q);
\draw[->] (H) to (S);
\draw[->] (I) to[bend right=40] (T);
\draw[->] (J) to[bend left=10] (Y);
\draw[->] (J) to[bend left=18] (Z);
\draw[->] (K) to (J);
\draw (K) edge[loop above] (K);
\draw[->] (K) to (L);
\draw[->] (K) to [bend left](M);
\draw[->] (L) to (N);
\draw[->] (M) to[bend right=9] (T);
\draw (N) edge[loop above] (N);
\draw[->] (O) to[bend left=40] (A);
\draw[->] (O) to (B);
\draw[->] (O) to (W);
\draw[->] (P) to[bend left] (N);
\draw[->] (R) to[bend left] (O);
\draw[->] (R) to[bend left] (P);
\draw (R) edge[loop above] (R);
\draw[->] (R) to (S);
\draw[->] (S) to (T);
\draw[->] (T) to (N);
\draw[->] (U) to (Lambda);
\draw[->] (X) to (Xi);
\draw (Y) edge[loop below] (Y);
\draw[->] (Y) to (Z);
\draw[->] (Z) to[bend left=42] (N);
\draw[->] (Gamma) to (Y);
\draw[->] (Delta) to[bend left=15] (A);
\draw[->] (Lambda) to[bend right] (Pi);
\draw (Xi) edge[loop below] (Xi);
\draw[->] (Xi) to[bend right] (Sigma);
\draw (Pi) edge[loop above] (Pi);
\draw[->] (Sigma) to (Pi);
\end{tikzpicture}
\caption{The graph of transfers between all weak $2$-chains listed in the above table with in yellow the starting weak $2$-chains from $1$ to $45$ of Lemma \ref{LEMweak2chainsb2} and in red the weak $2$-chains whose $3$-jump is equal to $2$\label{figgraph}}
\end{figure}
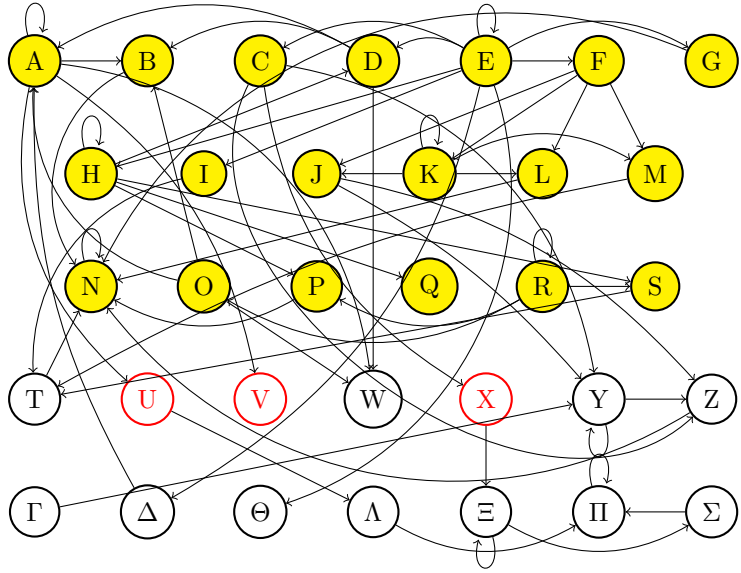
\vspace{0.2cm}
 
Moreover, Lemma \ref{LEMweak2chainsrulesout} shows that the weak $2$-chains $U$, $V$, $X$ do not belong to a weak chain of length at least $6$. Since $L=b-a\geq 6$, we infer that all $3$-jumps along the subchain (\ref{11nov1}) are equal to $1$. Consequently, since all weak $2$-chains of Lemma \ref{LEMweak2chainsb2} have the index $3$ in first and second position, we have 
\[
\bar{B}_3 := \mathcal{B}_{k_{\bar{j}(3)}}^-/z_3 < \mathcal{B}_a^+/z_1 +\b_3.
\]
Therefore, we can repeat the proof of Case 2 by replacing $z_2$ by $z_3$ to get a contradiction.\\

\noindent Subcase 2.2: The weak $2$-chain associated with the subchain $[1 i_1^0\cdots i_{\ell_0}^0  1 i_1^{1}\cdots i_{\ell_{1}}^{1} 1]$ is one of the weak $2$-chains 46 or 47 of Lemma \ref{LEMweak2chainsb2}.\\
Thanks to Lemma \ref{LEMweak2chainsrulesoutotherchains}, we know that $\b_3 \leq 3\b_2/5$. Moreover, since the corresponding weak $2$-chains start with $34234$ or $43243$, we have 
 \[
 \mathcal{B}_{a}^+/z_1  =   \frac{6a+1}{6z_1} <  \mathcal{B}_{k_{\bar{j}(2)}}^-/z_2 = \frac{6k_{\bar{j}(2)}-1}{6z_2} <  \mathcal{B}_{k_{\bar{j}(2)}}^+/z_2 = \frac{6k_{\bar{j}(2)}+1}{6z_2}  <  \mathcal{B}_{a+1}^-/z_1  =  \frac{6a+5}{6z_1}.
\]
Thus, we have $k_{\bar{j}(2)} \geq a+1$ and as a consequence
\[
z_2 > \left(\frac{6k_{\bar{j}(2)}+1}{6a+5}\right) z_1 \geq \left( \frac{6a+7}{6a+5}\right) z_1 = (1+\kappa) z_1 \quad \mbox{with} \quad \kappa := \frac{2}{6a+5}.
\]
Define the sequence $\{u_k\}_{k\in \N}$ by 
$$
u_k :=  \mathcal{B}_{a+k+1}^-/z_1 - \mathcal{B}_{k_{\bar{j}(2)}+k}^+/z_2  = \mathcal{B}_{a+1}^-/z_1- \mathcal{B}_{k_{\bar{j}(2)}}^+/z_2  + k \left(\frac{1}{z_1}-\frac{1}{z_2}\right) \qquad \forall k\in \N.
$$
By construction, we have $u_0>0$ and $u_{k+1}-u_k = 1/z_1- 1/z_2 >0$ for all $k\in \N$. Therefore, $\{u_k\}_{k\in \N}$ is increasing and there is $k\in \N^*$ such that 
\[
u_k\geq \frac{2}{5z_2} \quad \mbox{and} \quad 
k\leq  \left\lceil \frac{2/(5z_2)}{1/z_1-1/z_2} \right\rceil \leq  \left\lceil \frac{2}{5\kappa} \right\rceil=  \left\lceil \frac{6a}{5} + 1  \right\rceil = \left\lceil \frac{6a}{5}  \right\rceil+1.
\]
If the kwai $\mathcal{K}_{a+k}/z_1$ is contained in $\mathcal{K}_0$, that is if $a+k\leq b-1$, then the inequality $u_k\geq 2/(5z_2)$ implies, because $\b_3 \leq 3\b_2/5$, that the interval from  $\mathcal{B}_{k_{\bar{j}(2)}+k}^+/z_2$ to $\mathcal{B}_{a+k+1}^-/z_1$ has length at least $2/(5z_2)\geq 2\b_3$. But such an interval cannot be covered by an union of $z_3$-bridges and $z_4$-bridges, so we get a contradiction to the fact that $(z_1, z_2, z_3,z_4)$ is a $\mathcal{K}_0$-covering. This argument works provided $b-a\geq \left\lceil \frac{6a}{5}  \right\rceil+2$ which is satisfied because we assumed that $b-a\geq 2a+2$.
 \end{proof}

\begin{lemma}\label{LEMHanoiairport2}
$z_1\geq 47/5 \Rightarrow b-a\geq  2a+2$.
\end{lemma}

\begin{proof}[Proof of Lemma \ref{LEMHanoiairport2}]
Let $z_1\geq 13$ then there are $q\in \N$ with $q\geq 2$ and $r\in [0,6)$ such that $z_1=1+6q+r$. Then we have 
\[
a = \left\lceil \frac{z_1-1}{6} \right\rceil =  \left\lceil q+\frac{r}{6} \right\rceil = q+  \left\lceil \frac{r}{6} \right\rceil, \quad b =  \left\lfloor \frac{5z_1+1}{6} \right\rfloor =  \left\lfloor 5q+1+ \frac{5r}{6} \right\rfloor =  5q+1 +  \left\lfloor \frac{5r}{6} \right\rfloor. 
\]
and 
\[
b-a \geq 4q \geq 2q+4 \geq 2 q + 2 \left\lceil \frac{r}{6} \right\rceil +2 = 2a+2,
\]
because $q\geq 2$. If $z_1\in (7,13)$, then we have $a=1+\lceil r/6\rceil=2$, $b=6+\lfloor 5r/6\rfloor$, so that $b-a = 4+  \lfloor 5r/6\rfloor \geq 6=2a+2$ if and only if $ \lfloor 5r/6\rfloor\geq 2$, which is equivalent to $r\geq 12/5$ and implies that $z_1\geq 7+12/5=47/5$ verifies $b-a\geq 2a+2$.  \end{proof}

 \subsection{Checking the covering property over compact sets}\label{SEC4Dcompact} 

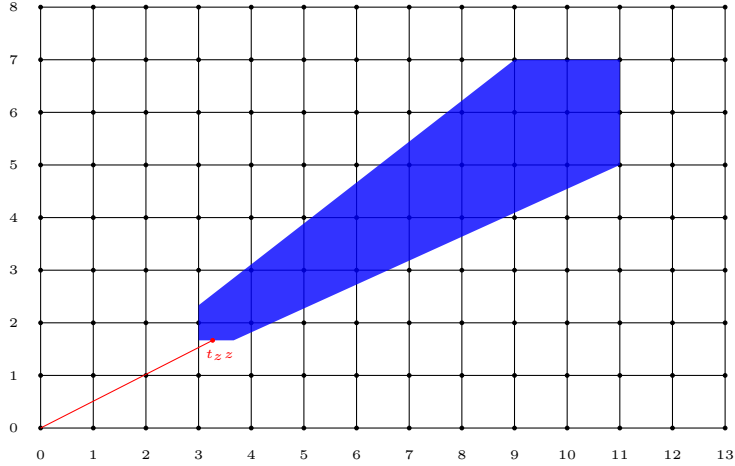
\begin{figure}[H]
\centering
\begin{tikzpicture}[scale=0.7]
\draw (0,0) grid(13,8);
\foreach \x in {0,...,13}
{
	\draw (\x,-0.25) node[below]{\tiny \color{black}{$\x$}};
}
\foreach \y in {0,...,8}
{
	\draw (-0.25,\y) node[left]{\tiny \color{black}{$\y$}};
}
\foreach \x in {0,...,13}
{
	\foreach \y in {0,...,8}
	{
		\filldraw[black] (\x,\y)circle(1pt);
	}
}
\foreach \k in {0,...,0}
{
	\pgfmathsetmacro\a{\k + 0.75};
	\pgfmathsetmacro\b{\k + 0.25};
	\fill[color=blue, opacity=.8] ( 2.25 / \a,1.25/ \a ) -- ( 2.75 / \a,1.25 / \a) -- ( 2.75 / \b,1.25 / \b ) -- ( 2.75 / \b,1.75 / \b ) -- ( 2.25 / \b,1.75 / \b) -- ( 2.25 / \a,1.75 / \a ) -- cycle; 
}
\draw [color=red] (0,0) -- (2.45/0.75,1.25/0.75);
\filldraw[red] (2.45/0.75,1.25/0.75)circle(1pt);
\draw (2.55/0.75,1.25/0.75) node[below]{\tiny \color{red}{$t_z z$}};
\end{tikzpicture}
\caption{If $[1,\infty)^d  \subsetneq \mathcal{K}^{d}(\delta_d) / \mathcal{K}_{\llbracket 0,N-1 \rrbracket}(\delta_d)$ does not hold, then there is a point of the form $t_zz$ satisfying (\ref{EQ15nov})  such that $tz\notin \mathcal{K}^{d}(\delta_d) / \mathcal{K}_{\llbracket 0,N-1 \rrbracket}(\delta_d)$ for all $t\in [1,t_z)$ \label{FigLowerface}}
\end{figure}

The aim of this section is to provide a method for checking a covering property of the form 
\begin{eqnarray}\label{16oct0}
[1,\infty)^d \subset \mathcal{K}^{d}\left( \delta_d \right) / \mathcal{K}_{\llbracket 0,N-1 \rrbracket} \left( \delta_d \right)
\end{eqnarray}
with a computer program if we already know that the property is satisfied outside of a compact set and to apply it to conclude the proof of Theorem \ref{THM2}. So, before checking our covering property in dimension $4$, we fix $d\in \N^*$ and $N\in \N^*$ and prove the validity of the method in dimension $d$. Our strategy consists in showing that the covering property can be verified by checking it on finitely many points. Let $P\subset \R^d$ be a compact set such that 
\begin{eqnarray}\label{16oct1}
[1,\infty)^d \setminus P \subset \mathcal{K}^{d}(\delta_d) / \mathcal{K}_{\llbracket 0,N-1 \rrbracket}(\delta_d),
\end{eqnarray}
if the covering property does not hold globally, that is, if 
\[
[1,\infty)^d  \subsetneq \mathcal{K}^{d}(\delta_d) / \mathcal{K}_{\llbracket 0,N-1 \rrbracket}(\delta_d), 
\]
then the open set (note that the family of feathers $\{F_{k,l}^d(\delta_d) =  \mathcal{K}_k^d(\delta_d)/\mathcal{K}_l(\delta_d)\}_{k\in \N^d, l\in \llbracket 0,N-1 \rrbracket}$ is locally finite),
\[
\mathcal{O}:= (1,\infty)^d\setminus \mathcal{K}^{d}(\delta_d) / \mathcal{K}_{\llbracket 0,N-1 \rrbracket}(\delta_d)
\]
contains a point $z=(z_1,\ldots,z_d)$. Then, since by (\ref{16oct1}) we know that $tz\notin \mathcal{O}$ for $t>0$ large enough, the supremum $t_z$ of $t>0$ such that $tz\in \mathcal{O}$ is well-defined and we have (see Figure \ref{FigLowerface})
\begin{eqnarray}\label{EQ15nov}
t_z \, z \in \bar{O} \, \cap  \, \left( \left( \bigcup_{k\in \N^d, l\in \llbracket 0,N-1 \rrbracket} \partial^- F_{k,l}^d(\delta_d) \right)\setminus \left( \bigcup_{k\in \N^d, l\in \llbracket 0,N-1 \rrbracket} \left[ F_{k,l}^{d}(\delta_d)\right]^{>} \right) \right),
\end{eqnarray}
where for each $k\in \N^d$ and $l\in \llbracket 0,N-1 \rrbracket$, $ \partial^- F_{k,l}^d(\delta_d)$ denotes the lower face of $F_{k,l}^d(\delta_d)$ given by (see Proposition \ref{PROPFeather})
\begin{multline*}
 \partial^- F_{k,l}^d(\delta_d) =  \Bigl\{ z \in F_{k,l}^d(\delta_d) \, \vert \, tz\notin F_{k,l}^d(\delta_d) \, \forall t \in [0,1)\Bigr\} \\
  = F_{k,l}^d(\delta_d) \, \cap \, \bigcup_{i=1}^d \left\{ z \in \R^d \, \vert \, z_i = \frac{k_i+\delta_d}{l+1-\delta_d}\right\} 
\end{multline*}
and $[F_{k,l}^{d}(\delta_d)]^{>}$ denotes its complement in $F_{k,l}^d(\delta_d)$, that is, 
\[
\left[ F_{k,l}^{d}(\delta_d)\right]^{>} := F_{k,l}^d(\delta_d) \setminus \partial^- F_{k,l}^d(\delta_d).
\]

As a consequence, if we show that for every $k\in \N^d$, every $l\in \llbracket 0,N-1 \rrbracket$ and every $i=1,\ldots,d$, the hyperplan
\[
H_i^{k,l} :=\left\{ z \in \R^d \, \vert \, z_i = \frac{k_i+\delta_d}{l+1-\delta_d}\right\} 
\]
satisfies
\[
H_i^{k,l} \cap \bar{\mathcal{O}} \subset \bigcup_{k'\in \N^d, l'\in \llbracket 0,N-1 \rrbracket} \left( H_i^{k,l} \cap \left[ F_{k',l'}^{d}(\delta_d)\right]^{>}\right),
\]
then we can infer that (\ref{16oct0}) is satisfied. For every $k\in \N^d$, $l\in \llbracket 0,N-1 \rrbracket$, $i=1,\ldots,d$, and every $k'\in \N^d$, $l'\in \llbracket 0,N-1 \rrbracket$, the trace of $F_{k',l'}^{d}(\delta_d)$ on $H_i^{k,l}$, given by  $H_i^{k,l} \cap F_{k',l'}^{d}(\delta_d)$, is a feather of dimension $d-1$ in $H_i^{k,l}$. So we can repeat this argument and show that it is sufficient to verify the covering property over a finite set in order to check that it holds over a given compact set. Since it is more convenient, we express below this idea in term of kwais. \\

\begin{figure}[H]
\centering
\begin{tikzpicture}
\tikzstyle{root}=[circle,draw,thick,fill=yellow]
\tikzstyle{childone}=[circle,draw,thick, fill=blue]
\tikzstyle{childtwo}=[circle,draw,thick,fill=red]
\node[childtwo] (CC1) at (-1,-4) {};
\node[childtwo] (CC2) at (0,-4) {};
\node[childtwo] (CC3) at (2,-4) {};
\node[childtwo] (CC4) at (3,-4) {};
\node[childtwo] (CC5) at (4,-4) {};
\node[childtwo] (CC6) at (5.5,-4) {};
\node[childtwo] (CC7) at (6.5,-4) {};
\node[childtwo] (CC8) at (8.5,-4) {};
\node[childtwo] (CC9) at (9.5,-4) {};
\node[childtwo] (CC10) at (12,-4) {};
\node[childone] (C1) at (0,-6) {};
\node[childone] (C2) at (3,-6) {};
\node[childone] (C3) at (6,-6) {};
\node[childone] (C4) at (9,-6) {};
\node[childone] (C5) at (12,-6) {};
\node[root] (R) at (6,-8) {$\mathfrak{T}_0$};
\draw[-] (R) to (C1);
\draw[-] (R) to (C2);
\draw[-] (R) to (C3);
\draw[-] (R) to (C4);
\draw[-] (R) to (C5);
\draw[-] (CC1) to (C1);
\draw[-] (CC2) to (C1);
\draw[-] (CC3) to (C2);
\draw[-] (CC4) to (C2);
\draw[-] (CC5) to (C2);
\draw[-] (CC6) to (C3);
\draw[-] (CC7) to (C3);
\draw[-] (CC8) to (C4);
\draw[-] (CC9) to (C4);
\draw[-] (CC10) to (C5);
\end{tikzpicture}
\caption{A rooted tree with the root or $0$-child in yellow, the $1$-children in blue and $2$-children in red\label{rootedtree}}
\end{figure}
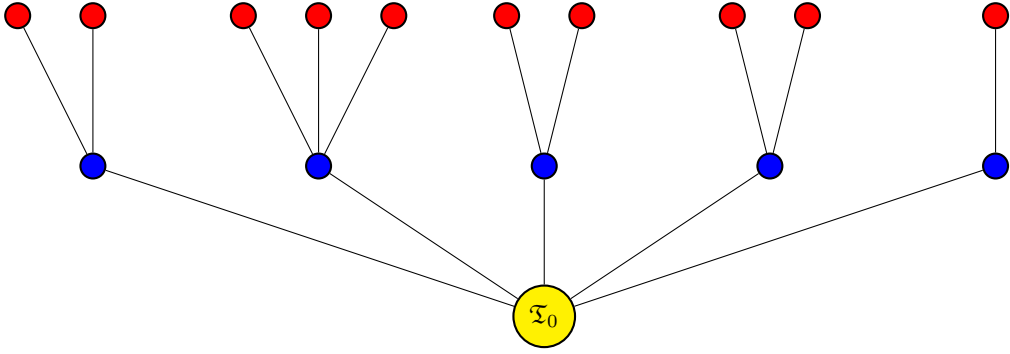

We need to introduce some notations. First we recall that $\mathcal{S}^d$ stands for the closed convex set of tuples $(z_1,\ldots,z_d) \in \R^d$ satisfying 
\[
1\leq z_1\leq \cdots \leq z_d.
\] 
Then, we denote by  $\mathfrak{S}_d$ the set of permutations of the set $\{1,\ldots,d\}$ and  given a compact set $P\subset \mathcal{S}^d$, we define the set $\mathfrak{S}_d(P)\subset [1,\infty)^d$ by 
\[
\mathfrak{S}_d(P):= \Bigl\{ \left(z_{\sigma(1)},\ldots, z_{\sigma(d)}\right) \, \vert \, (z_1,\ldots,z_d)\in P, \, \sigma \in \mathfrak{S}_d \Bigr\},
\]
and we say that $P$ is symmetric if $P=\mathfrak{S}_d(P)$. Moreover, we denote by $\I$ the set of disjoint unions of finitely many closed intervals in $\R$ and if $A \in \I$ is the disjoint union of intervals $[a_l,b_l] \subset \R$  with $l=1, \ldots, L$, then we define the sets $A^-, A^+$ and $A^>$ by  
\[
A^-:= \Bigl\{a_l\,\vert \, l=1,\ldots, L\Bigr\}, \quad
A^+:= \Bigl\{b_l\,\vert \, l=1,\ldots, L\Bigr\}
\]
and
\[
A^> := \bigcup_{l=1,\ldots,L} (a_l,b_l],
\]
with the convention that $(a_l,b_l]=\emptyset$ if $a_l=b_l$. Then, we denote by $\mbox{Clos}(A)$ the closure of a subset $A$ of $\R$. We recall the a rooted tree is a graph which is a tree in which a special node is singled out, we call this node root or $0$-child. Then, if $\mathfrak{T}$ is a rooted tree with root $\mathfrak{T}_0$, we call $l$-child with $l\in \N$ any node which is at distance $l$ from $\mathfrak{T}_0$ (see Figure \ref{rootedtree}). 

\vspace{0.2cm}
\vspace{0.2cm}

In the following result, we deal with a rooted tree $\mathfrak{T}$ with root $\mathfrak{T}_0$ whose nodes belong to $\R\times \I$. Given a $d$-child $(z,\mathfrak{K})$, we define its $\R$-ancestry as 
\[
\mbox{Anc}_{\R}(z,\mathfrak{K}) := (z_1,\cdots, z_{d-1},z),
\]
where $z_1, \ldots, z_{d-1}\in \R^d$ are associated with $\mathfrak{K}_1,\ldots, \mathfrak{K}_{d-1}\in \I$ in such a way that each $(z_i,\mathfrak{K}_i)$ is a $i$-child (with $i=1,\ldots, d-1$) and the edges from $\mathfrak{T}_0$ to $(z,\mathfrak{K})$
passing through the nodes $(z_1,\mathfrak{K}_1), \ldots, (z_{d-1},\mathfrak{K}_{d-1})$ connect $\mathfrak{T}_0$ to $(z,\mathfrak{K})$.

\begin{proposition}\label{3August1}
Let $d, N \in \N^*$ and $C>1$ be fixed and let $P\subset \mathcal{S}^d$ be a compact set such that  
\[
P\subset \mathcal{C}:= [1,C]^d \quad \mbox{and} \quad \mathcal{S}^d  \setminus P \subset \mathcal{K}^{d}(\delta_d) / \mathcal{K}_{\llbracket 0,N-1 \rrbracket}(\delta_d). 
\]
Let $\mathfrak{T}$ be the rooted tree whose root in $\R\times \I$ is given by
\[
\mathfrak{T}_0 := \left( 1, \mathfrak{K}_0\right) \quad \mbox{with} \quad \mathfrak{K}_0:= \mathcal{K}_{\llbracket 0,N-1 \rrbracket}(\delta_d)
\]
and such that the set of children of a $i$-child $(z,\mathfrak{K})$ with $i\in \{0,\ldots,d-1\}$ is given by the set of all $(z',\mathfrak{K}')\in \R\times \I$ such that
\[
z' \in \left( \mathcal{K}^-(\delta_d) / \mathfrak{K}^+\right) \cap (1,C] \quad \mbox{and} \quad \mathfrak{K}' := \mbox{Clos} \left( \mathfrak{K} \cap \mathcal{K}^>(\delta_d)/z'\right).
\]
Then, the following properties are equivalent:
\begin{itemize}
\item[(i)] $[1,\infty)^d \subset \mathcal{K}^{d}(\delta_d) / \mathcal{K}_{\llbracket 0,N-1 \rrbracket}(\delta_d)$.
\item[(ii)] For every $(z_1,\ldots, z_d) \in [1,\infty)^d$, $\mathcal{K}_{\llbracket 0,N-1\rrbracket}(\delta_d) \cap \left( \cap_{i=1}^d \mathcal{K}(\delta_d)/z_i \right)\neq \emptyset$.
\item[(iii)]  For every $d$-child $(z,\mathfrak{K})$ of $\mathfrak{T}$ such that $\mbox{Anc}_{\R}(z,\mathfrak{K})\in \mathfrak{S}_d(P)$, we have $\mathfrak{K}\neq \emptyset$.
\end{itemize}
\end{proposition}

\begin{proof}[Proof of Proposition \ref{3August1}]
We start with the following lemma. 

\begin{lemma}\label{LEM24oct}
Let $c$ be an integer $\geq 1$, $C>1$, $\mathfrak{K}\in \I$ and $Q\in [1,\infty)^c$ be a symmetric compact set  such that  
\begin{eqnarray}\label{EQ13nov1}
Q \subset [1,C]^c \quad \mbox{and} \quad \mathfrak{K} \cap \left( \cap_{i=1}^{c} \mathcal{K}(\delta_d)/z_i \right) \neq \emptyset \qquad \forall (z_1,\ldots, z_{c}) \in [1,\infty)^{c} \setminus Q.
\end{eqnarray}
Then the following properties are equivalent:
\begin{itemize}
\item[(i)] For every $(z_1,\ldots,z_{c})\in [1,\infty)^{c}$, $\mathfrak{K} \cap \left( \cap_{i=1}^{c} \mathcal{K}(\delta_d)/z_i \right) \neq \emptyset$.
\item[(ii)] If $c\geq 2$, then for every $(z_2,\ldots, z_{c}) \in (1,C]^{c-1}$ and every $z \in (1,C]$, we have 
\[
\left\{
\begin{array}{l}
z \in   \mathcal{K}^-(\delta_d) / \mathfrak{K}^+ \\
(z,z_2,\cdots, z_{c}) \in Q
\end{array}
\right.
\quad \Longrightarrow \quad \mbox{Clos} \left(  \mathfrak{K} \cap \mathcal{K}^>(\delta_d)/z\right) \cap  \left( \cap_{i=2}^{c} \mathcal{K}(\delta_d)/z_i \right) \neq \emptyset.
\]
If $c=1$, then for every $z \in (1,C]$, we have 
\[
z \in   \mathcal{K}^-(\delta_d) / \mathfrak{K}^+ \quad \Longrightarrow \quad \mbox{Clos} \left(  \mathfrak{K} \cap \mathcal{K}^>(\delta_d)/z\right) \neq \emptyset \quad \Longleftrightarrow \quad   \mathfrak{K} \cap \mathcal{K}^>(\delta_d)/z \neq \emptyset.
\]
\end{itemize}
\end{lemma}

\begin{proof}[Proof of Lemma \ref{LEM24oct}]
To prove that (i) $\Rightarrow$ (ii), assume that (i) is satisfied and fix a tuple $(z_1,\ldots ,z_{c-1},z)\in (1,C]^{c}$ (where we only take $z\in (1,C]$ if $c=1$). Since $(1,\infty)^{c}$ is open, there is $\epsilon\in (0,1)$ such that for every $t\in (1-\epsilon,1)$, $(z_1,\ldots,z_{c-1},tz)\in (1,\infty)^{c}$, so by (i) we have
 \[
E\cap \mathcal{K}(\delta_d)/(tz) \neq \emptyset \quad \forall t\in (1-\epsilon,1) \quad \mbox{with} \quad E:=\mathfrak{K} \cap \left( \cap_{i=1}^{c-1} \mathcal{K}(\delta_d)/z_i \right)
 \]
 (where we take $E :=\mathfrak{K}$ if $c=1$).   Thus, for every $t\in (1-\epsilon,1)$, there are $\lambda^t\in E$ and $k^t\in \N$ such that 
\begin{eqnarray}\label{12nov1}
\mathcal{K}^-_{k^t}(\delta_d)/z < \mathcal{K}^-_{k^t}(\delta_d)/(tz) \leq \lambda^t \leq \mathcal{K}^+_{k^t}(\delta_d)/(tz).
\end{eqnarray}
As a consequence, if 
\begin{eqnarray}\label{EQ13nov}
\mbox{Clos} \left(  \mathfrak{K} \cap \mathcal{K}^>(\delta_d)/z\right) \cap  \left( \cap_{i=1}^{c-1} \mathcal{K}(\delta_d)/z_i \right) = \emptyset,
\end{eqnarray}
then we have $E\cap  \mathcal{K}_{k^t}^>(\delta_d)/z=\emptyset$ (in the case $c=1$, we only assume that $E\cap  \mathcal{K}^>(\delta_d)/z=\emptyset$ which implies $E\cap  \mathcal{K}_{k^t}^>(\delta_d)/z=\emptyset$), so that  $\lambda^t > \mathcal{K}^+_{k^t}(\delta_d)/z$ or $\lambda^t \leq \mathcal{K}^-(\delta_d)_{k^t}/z$ which, by (\ref{12nov1}), is prohibited. Then, if (\ref{EQ13nov}) is satisfied (if $E\cap  \mathcal{K}^>(\delta_d)/z=\emptyset$ for $c=1$), we have
\[
\mathcal{K}^+_{k^t}(\delta_d)/z < \lambda^t \leq \mathcal{K}^+_{k^t}(\delta_d)/(tz) \qquad \forall t\in (1-\epsilon,1)
\]
and by letting $t$ tend to $1$ (note that $E$ is compact and the set of $k$ such that $E\cap \mathcal{K}(\delta_d)/(tz) \neq \emptyset$ is finite), we obtain $\lambda \in E$ and $k\in \N$ such that $\lambda \in\mathcal{K}^+_k(\delta_d)/z$, which contradicts (\ref{EQ13nov}) ($E\cap  \mathcal{K}^>(\delta_d)/z=\emptyset$ if $c=1$).

Let us now prove that (ii)$\Rightarrow$ (i). The set $\mathcal{O}$ of $(z_1,\ldots, z_{c})\in [1,\infty)^{c}$ such that 
\[
 \mathfrak{K} \cap \left( \cap_{i=1}^{c} \mathcal{K}(\delta_d)/z_i \right)=\emptyset
 \]
 is open in $[1,\infty)^{c}$, symmetric and contained in $Q$ (by (\ref{EQ13nov1})). Therefore, since the set 
 \[
\Gamma := \left\{ \frac{k+\delta_d}{l+\delta_d} \cdot \frac{e^+}{f^+} \, \vert \, k,l \in \N, \, e^+, f^+ \in \mathfrak{K}^+\right\} \subset (0,\infty)
\]
is countable, there is $\bar{z}=(\bar{z}_1,\ldots, \bar{z}_{c})\in \mathcal{O}$ such that 
 \begin{eqnarray}\label{8nov2}
\frac{\bar{z}_i}{\bar{z}_j} \notin \Gamma \qquad \forall i,j \in \{1,\ldots, c\} \mbox{ such that }i \neq j.
\end{eqnarray}
Set
\[
\hat{t} := \sup \Bigl\{ t \geq 1\, \vert \, tz\in \mathcal{O}\Bigr\}  \quad \mbox{and} \quad \hat{z}:=\hat{t}\bar{z}.
\]
By (\ref{EQ13nov1}), $\hat{t}$ is well-defined, $\hat{z}\in Q$ and the compact set 
\[
A:= \mathfrak{K} \cap \left( \cap_{i=1}^{c} \mathcal{K}(\delta_d)/\hat{z}_i \right)
 \]
 is a non-empty disjoint union of closed intervals. We claim that there is $i\in \{1,\ldots, c\}$ such that 
 \[
 A = A^+ \subset \mathcal{K}^-(\delta_d)/\hat{z}_i.
 \]
 As a matter of fact, if $\lambda$ in $A$ does not belong to $ \mathcal{K}^-/\hat{z}_i$ for all $i=1,\ldots, c$, then for all $i\in \{1,\ldots, c\}$ and all $t<1$ sufficiently close to $1$ we have $t\lambda \in \mathcal{K}/\hat{z}_i$, which is equivalent to 
 \[
 \lambda \in \cap_{i=1}^{c} \mathcal{K}(\delta_d)/(t\hat{z}_i) \qquad \forall t<1 \mbox{ sufficiently close to } 1.
 \]
 But this property contradicts the fact that $t\hat{z}\in \mathcal{O}$ for  all $t<1$ sufficiently close to $1$. Therefore, we have shown that for every $\lambda \in A$, there is $i=i(\lambda) \in \{1,\ldots,c\}$ such that $\lambda \in  \mathcal{K}^-(\delta_d)/\hat{z}_i$. This shows that the compact set $A$ is a finite set satisfying 
 \[
A= A^+\subset \cup_{i=1}^c \mathcal{K}^-(\delta_d)/\hat{z}_i.
 \]
 But if there are $e^+,f^+\in A=A^+$ and $i\neq j$ in $\{1,\ldots, c\}$ such that $e^+ \in \mathcal{K}^-(\delta_d)/\hat{z}_i$ and $f^+ \in \mathcal{K}^-(\delta_d)/\hat{z}_j$, then there are $k_i,k_j\in \N$ such that
 \[
 \hat{z}_i = \frac{k_i+\delta_d}{e^+} \quad \mbox{and} \quad  \hat{z}_j = \frac{k_j+\delta_d}{f^+},
 \]
 which implies that $\bar{z}_i/\bar{z}_j=\hat{z}_i/\hat{z}_j$ and contradicts (\ref{8nov2}). In conclusion, the claim is proved and we note that by symmetry ($\mathcal{O}$ and $Q$ are symmetric) we can assume without loss of generality that $i=1$. We claim now that $A=A^+\subset \mathfrak{K}^+$. As a matter of fact, suppose for contradiction that there is $\lambda \in A=A^+$ such that $\lambda \notin \mathfrak{K}^+$. Then, we have
 \[
\frac{\lambda}{t} \in \mathfrak{K} \quad \mbox{and} \quad  \frac{\lambda}{t} \in \cap_{i=1}^{c} \mathcal{K}(\delta_d)/(t\hat{z}_i) \qquad \forall t<1 \mbox{ sufficiently close to } 1.
 \]
 But this property contradicts the fact that $t\hat{z}\in \mathcal{O}$ for  all $t<1$ sufficiently close to $1$. In conclusion, we have 
 \[
  A= A^+ \subset \mathcal{K}^-(\delta_d)/\hat{z}_{1} \cap \mathfrak{K}^+,
  \]
  which implies that 
  \[
   z=\hat{z}_{1} \in \left( \mathcal{K}^-(\delta_d)/ \mathfrak{K}^+\right) \cap (1,C] \quad \mbox{and} \quad \hat{z} =(z, \hat{z}_2,\ldots, \hat{z}_{c}) \in Q.
  \]
  Let us now distinguish the case $c\geq 2$ and $c=1$. If $c=1$, then $\mathfrak{K} \cap \mathcal{K}^>(\delta_d)/z=\emptyset$, because that set is contained in $A$ which is contained in  $ \mathcal{K}^-(\delta_d)/z \cap \mathfrak{K}^+$. If $c\geq 2$, then we have 
  \[
   \mbox{Clos} \left( \mathfrak{K} \cap \mathcal{K}^>/z\right) \cap \left( \cap_{i=2,\ldots,c} \mathcal{K}(\delta_d)/\hat{z}_i \right) \subset A \subset  \mathcal{K}^-(\delta_d)/z \cap \mathfrak{K}^+,
  \]
  which shows that $\mbox{Clos} \left( \mathfrak{K} \cap \mathcal{K}^>/z\right) \cap \left( \cap_{i=2,\ldots,c} \mathcal{K}(\delta_d)/\hat{z}_i \right)=\emptyset$ because  a point of $\mbox{Clos} \left( \mathfrak{K} \cap \mathcal{K}^>/z\right) $ cannot belong to $ \mathcal{K}^-(\delta_d)/z \cap \mathfrak{K}^+$. In the two cases, we obtain a contradiction to assumption (ii).
 \end{proof}

The equivalence of (i) and (ii) follows from Proposition \ref{PROPcaracx}. Let us prove (ii)  $\Leftrightarrow$ (iii) by induction over $d\in \N^*$. In fact, we need to prove more than this, we show by induction over $d\in \N^*$, that if we have $C>1$, $\mathfrak{K}_0$ a given element of $\I$, $P\subset \mathcal{S}^d$  a compact set such that  
\[
P\subset \mathcal{C}:= [1,C]^d \quad \mbox{and} \quad \mathfrak{K}_0 \cap \left( \cap_{i=1}^d \mathcal{K}(\delta_d)/z_i \right)\neq \emptyset \quad \forall (z_1,\ldots, z_d) \in \mathcal{S}^d \setminus P,
\]
and $\mathfrak{T}$  the rooted tree with nodes in $\R\times \I$, a root given by $\mathfrak{T}_0 := \left( 1, \mathfrak{K}_0\right)$
and such that the set of children of a $i$-child $(z,\mathfrak{K})$ with $i\in \{0,\ldots,d-1\}$ is given by the set of all $(z',\mathfrak{K}')\in \R\times \I$ such that
\[
z' \in \left( \mathcal{K}^-(\delta_d) / \mathfrak{K}^+\right) \cap (1,C] \quad \mbox{and} \quad \mathfrak{K}' := \mbox{Clos} \left( \mathfrak{K} \cap \mathcal{K}^>(\delta_d)/z'\right),
\]
then the two following properties are equivalent: 
\begin{itemize}
\item[(iv)] For every $(z_1,\ldots, z_d) \in [1,\infty)^d$, $\mathfrak{K}_0 \cap \left( \cap_{i=1}^d \mathcal{K}(\delta_d)/z_i \right)\neq \emptyset$.
\item[(v)]  For every $d$-child $(z,\mathfrak{K})$ of $\mathfrak{T}$ such that $\mbox{Anc}_{\R}(z,\mathfrak{K})\in \mathfrak{S}_d(P)$, we have $\mathfrak{K}\neq \emptyset$.
\end{itemize}
The result for $d=1$ follows from Lemma \ref{LEM24oct} and the formula $\mathfrak{K}_1 := \mbox{Clos} \left( \mathfrak{K} \cap \mathcal{K}^>(\delta_d)/z_1\right)$ for any $1$-child. We now assume that the equivalence (iv)  $\Leftrightarrow$ (v) has been proven for some $d\in \N^*$ and prove it for $d+1$. So, we assume that there are $C>1$, $\mathfrak{K}_0 \in \I$, $P\subset \mathcal{S}^{d+1}\cap [1,C]^{d+1}$  a compact set and $\mathfrak{T}$ a rooted tree satisfying the required assumption and we set $Q:=\mathfrak{S}_d(P)\in [1,C]^{c}$ with $c=d+1$. By construction, $Q$ is symmetric and the assumption of Lemma \ref{LEM24oct} is satisfied with $ \mathfrak{K}= \mathfrak{K}_0$. As a consequence, the above property (iv) is equivalent to the following one: For every $(z_2,\ldots, z_{c}) \in (1,C]^{d}$ and every $z \in (1,C]$, 
\[
\left\{
\begin{array}{l}
z \in   \mathcal{K}^-(\delta_d) / \mathfrak{K}_0^+ \\
(z,z_2,\cdots, z_{c}) \in Q
\end{array}
\right.
\quad \Longrightarrow \quad \mbox{Clos} \left(  \mathfrak{K}_0 \cap \mathcal{K}^>(\delta_d)/z\right) \cap  \left( \cap_{i=2}^{c} \mathcal{K}(\delta_d)/z_i \right) \neq \emptyset.
\]
But, if for every $z_1\in   \mathcal{K}^-(\delta_d) / \mathfrak{K}_0^+\cap (1,C]$, we set $ \mathfrak{K}_1:= \mbox{Clos} \left(  \mathfrak{K}_0 \cap \mathcal{K}^>(\delta_d)/z_1\right)$, then 
by setting 
\[
Q_1 =Q(z_1) := \Bigl\{ (z_2,\cdots, z_c) \in [1,\infty)^d \, \vert \, (z_1,z_2,\cdots, z_c) \in Q\Bigr\},
\]
we have (as a matter of fact, the proof of (i) $\Rightarrow$ (ii) in Lemma \ref{LEM24oct} shows that the required intersection is non-empty for any point outside of the compact set)
\[
Q_1 \subset [1,C]^d \quad \mbox{and} \quad \mathfrak{K}_1 \cap \left( \cap_{i=2}^{d+1} \mathcal{K}(\delta_d)/z_i \right) \neq \emptyset \qquad \forall (z_2,\ldots, z_{d+1}) \in [1,\infty)^{d} \setminus Q_1.
\]
By the induction hypothesis, for every $z_1 \in   \mathcal{K}^-(\delta_d) / \mathfrak{K}_0^+\cap (1,C]$, the latter property  is equivalent to a property of the form (v) for a tree with root  $\mathfrak{K}_1$, which shows exactly that our equivalence holds for $d+1$. 
\end{proof}

Applying the above method to the compact set $P\subset \mathcal{S}^4$ consisting of the points $(z_1,z_2,z_3,z_4)\in \mathcal{S}^4$ satisfying 
\begin{eqnarray}\label{16fev}
5z_1 \leq 47, \quad 2z_2\leq 5z_1, \quad z_2z_3+z_1z_3 \leq 8z_1z_2, \quad z_3z_4+z_2z_4\leq 10z_2z_3,
\end{eqnarray}
 we could establish by running a computer program in C++  that the two following covering properties hold:
 \[
 [1,\infty)^4 \setminus \left( P \setminus (1,31/5)^4 \right) \subset \mathcal{K}^{4}(\delta_4) / \mathcal{K}_{0}(\delta_4)
 \]
 and
 \[
 [1,31/5]^4 \subset  \mathcal{K}^{4}(\delta_4) / \mathcal{K}_{\llbracket 0,1 \rrbracket}(\delta_4).
 \]
Using a MacBook Pro with a processor of 2,6GHz, the time required to generate and check all the points provided by Proposition \ref{3August1} took about 9 hours. This method to check the covering property over compact set is rather naive, another method, in the spirit of what we did in Section \ref{SEC4Dinfinity}, could be to check all admissible ways to cover  $\mathcal{K}_{\llbracket 0,1 \rrbracket}(\delta_4)$ under the constraints (\ref{16fev}) by checking the emptyness or non-emptyness of convex polytopes. 

\subsection{Non-covering property of $\mathcal{K}^4(\delta_4)/\mathcal{K}_0(\delta_4)$}\label{SEC4Dnot}

To show that $[1,\infty)^4$ is not covered by $\mathcal{K}^4(\delta_4)/\mathcal{K}_0(\delta_4)$, we are going to prove that for every $\epsilon=(\epsilon_1, \epsilon_2,\epsilon_3,\epsilon_4)\in \R^4$ satisfying 
\[
0< \epsilon_1 < \epsilon_3 < \epsilon_2<\epsilon_4 < \frac{32}{3211}
\]
the point
\[
z=z_{\epsilon}:= \left(\frac{13(1-\epsilon_1) }{5}, \frac{3211(1-\epsilon_2)}{935}, \frac{247(1-\epsilon_3) }{55}, \frac{61009 (1-\epsilon_4)}{10285} \right)
\]
does not belong to $\mathcal{K}^4(\delta_4)/\mathcal{K}_0(\delta_4)$. If $z=z_{\epsilon} \in \mathcal{K}^4(\delta_4)/\mathcal{K}_0(\delta_4)$, then there is $k=(k_1,k_2,k_3,k_4)\in \N^4$ such that $z \in F_{k,0}^4(\delta_4)$. Then, by Proposition \ref{PROPFeather} (i), this means that we have 
\begin{eqnarray}\label{25sept2}
\frac{6k_i +1}{5}  \leq z_i \leq  6k_i+5 \qquad \forall i=1, 2, 3, 4
\end{eqnarray}
and
\begin{eqnarray}\label{25sept3}
\frac{6k_j +1}{6k_i+5}  \leq \frac{z_j}{z_i} \leq  \frac{6k_j+5}{6k_i+1} \qquad \forall i,j=1, 2, 3, 4.
\end{eqnarray}
We note that 
\[
 \frac{3211}{935} =  \frac{13^2 \cdot 19}{5\cdot 11 \cdot 17}, \quad \frac{247}{55} =  \frac{13\cdot 19}{5\cdot 11}, \quad  \frac{61009}{10285} =  \frac{13^2 \cdot 19^2}{5\cdot 11^2 \cdot 17},
\]
which gives
\[
 \frac{z_2}{z_1} = \frac{13\cdot 19}{11 \cdot 17}\left( \frac{1-\epsilon_2}{1-\epsilon_1} \right), \quad  \frac{z_3}{z_1} = \frac{19}{11} \left( \frac{1-\epsilon_3}{1-\epsilon_1} \right), \quad  \frac{z_4}{z_1} = \frac{13\cdot 19^2}{11^2 \cdot 17} \left( \frac{1-\epsilon_4}{1-\epsilon_1} \right)
 \]
  \[
 \frac{z_3}{z_2} = \frac{17}{13} \left( \frac{1-\epsilon_3}{1-\epsilon_2} \right), \quad  \frac{z_4}{z_2} = \frac{19}{11} \left( \frac{1-\epsilon_4}{1-\epsilon_2} \right), \quad  \frac{z_4}{z_3} = \frac{13\cdot 19}{11 \cdot 17}\left( \frac{1-\epsilon_4}{1-\epsilon_3} \right).
 \]
 And we also note that (\ref{25sept2}) yields
\[
 k_1 \leq  2 - \frac{13}{6}\, \epsilon_1, \quad  k_2  \leq  \frac{504}{187} -\frac{3211}{1122} \, \epsilon_2, \quad  k_3  \leq  \frac{118}{33} -\frac{247}{66} \, \epsilon_3
 \]
 and
 \[
  \frac{4792}{30855} -  \frac{61009}{61710} \, \epsilon_4  \leq k_4  \leq  \frac{29476}{6171} -\frac{61009}{12342} \, \epsilon_4, 
 \]
which, since $\epsilon_1,\epsilon_2,\epsilon_3,\epsilon_4>0$ and $\epsilon_4\in (0,32/3211)\subset (9584/61009)$, implies that
\begin{eqnarray}\label{16OCT}
 k_1 \in \{0,1\}, \quad  k_2 \in \{0,1,2\}, \quad  k_3 \in \{0,1,2,3\}   \quad \mbox{and} \quad k_4 \in \{1,2,3,4\}. 
 \end{eqnarray}
 Let us distinghish the cases $k_1=0$ and $k_1=1$. \\
 
\noindent Case 1: $k_1=0$.\\
The lower bound in (\ref{25sept3}) with $(i,j)=(1,2)$ gives, since $(1-\epsilon_2)/(1-\epsilon_1)\in (0,1)$, $6k_2+1\leq 5z_2/z_1 < 1235/187 < 7$, so we have $k_2=0$. Then, since $(1-\epsilon_3)/(1-\epsilon_2)\in (0,91/85)$ (it is equivalent to $6-91\epsilon_2+85\epsilon_3>0$ which is verified because $\epsilon_2<32/3211<6/91$), the lower bound in (\ref{25sept3}) with $(i,j)=(2,3)$ gives $6k_3+1\leq 5z_3/z_2 < 85(1-\epsilon_3)/(13(1-\epsilon_2) < 7$, so that $k_3=0$. Thus, since $(1-\epsilon_4)/(1-\epsilon_3)\in (0,1)$, the lower bound in (\ref{25sept3}) with $(i,j)=(3,4)$ gives $6k_4+1\leq 5z_4/z_3 < 1235/187 < 7$, so that $k_4=0$, which contradicts $ k_4 \in \{1,2,3,4\}$.\\ 
 
\noindent Case 2: $k_1=1$.\\
On the one hand, since $(1-\epsilon_3)/(1-\epsilon_1)\in (0,1)$, the lower bound in (\ref{25sept3}) with $(i,j)=(1,3)$ gives $6k_3+1 \leq 11 z_3/z_1<19$, so that $k_3\leq 2$. On the other hand, the upper bound in (\ref{25sept3}) with $(i,j)=(1,3)$ gives $6k_3+5 \geq 7z_3/z_1>11$ because $133 (1-\epsilon_3)/(11(1-\epsilon_1))>11$ is equivalent to $12+121\epsilon_1 -133 \epsilon_3>0$ which is satisfied because $\epsilon_3<32/3211<12/133$. Therefore, we have  $k_3 =2$. Recall that $k_2 \in \{0,1,2\}$. The upper bound in (\ref{25sept3}) with $(i,j)=(1,2)$ gives $6k_2+5 \geq 7z_2/z_1>5$ (because $1729 (1-\epsilon_2)/(187(1-\epsilon_1))>5$ is equivalent to $794+935 \epsilon_1-1729\epsilon_2 >0$ which is satisfied because $\epsilon_2<32/3211<794/1729$), so that $k_2\geq 1$. If $k_2=1$, then, since $(1-\epsilon_4)/(1-\epsilon_2)\in (0,1)$, the lower bound in (\ref{25sept3}) with $(i,j)=(2,4)$ gives $6k_4+1 \leq 11z_4/z_2<19$, so that $k_4\leq 2$. As a consequence, if $k_2=1$, then the upper bound in (\ref{25sept3}) with $(i,j)=(3,4)$ gives $6k_3+1 \leq (6k_4+5) z_3/z_4 <13$ (because $3179 (1-\epsilon_3)/(247(1-\epsilon_4))<13$ is equivalent to $32+3179 \epsilon_3-3211\epsilon_4 >0$ which is satisfied because $\epsilon_4<32/3211$),  that is, $k_3\leq 1$, a contradiction. If $k_2=2$, then, since $(1-\epsilon_3)/(1-\epsilon_2)>1$,  the upper bound in (\ref{25sept3}) with $(i,j)=(2,3)$ gives $17/13<z_3/z_2  \leq (6k_3+5)/(6k_2+1)=17/13$, a contradiction.  

\section{Proof of Theorem \ref{THM3}}\label{SECTHM3}

The case $n=2$ being obvious, we only prove in the next two sections the cases $n=3$ and $n=4$. The proof of the $n=4$ case refers implicitely to the notion of weak chains introduced in the proof of Theorem \ref{THM2}, but for sake of clarity we provide here a proof which is independent of the rest of the paper. 
  
\subsection{A preparatory result}

Recall that for every $\delta\in (0,1/2]$, $\mathcal{K}(\delta)$ and $\mathcal{B}(\delta)$ are defined by
\[
\mathcal{K}(\delta) = \bigcup_{k\in \N} \mathcal{K}_k(\delta) \quad \mbox{and} \quad \mathcal{B}(\delta) = \bigcup_{k\in \N} \mathcal{B}_k(\delta)
\]
where for each $k\in \N$, $\mathcal{K}_k(\delta)$ and $\mathcal{B}_k(\delta)$ are given by
\[
\mathcal{K}_k(\delta):= \left[ \mathcal{K}_k^-(\delta), \mathcal{K}_k^+(\delta)\right] =  \left[k+ \delta,k+1-\delta\right] \quad \mbox{and} \quad
\mathcal{B}_k(\delta) := \left( \mathcal{B}_k^-(\delta),  \mathcal{B}_k^+(\delta)  \right) = \left(k-\delta,k+\delta\right),
\]
and recall that we have set  for every integer $d$, $\delta_d:=1/(d+2)$. 

\begin{proposition}\label{PROPTHM3}
Let $m$ be an integer $\geq 2$ and $d:=m-1$. If the Lonely Runner Conjecture with Free Starting Positions in Obstruction Form does not hold for $m$, then for every $(a,b) \in \{0,\ldots, d\}$, there are $P_0, P_1, \ldots, P_d \in [0,1)^{d+1}$ and $(z_1,\ldots, z_d) \in (1,\infty)^d$ such that 
\begin{eqnarray}\label{23jan1}
P_{a}= P_{b}, \quad z_1 < z_2 < \cdots < z_d
\end{eqnarray}
and
\begin{eqnarray}\label{23jan3}
\mathcal{K}(\delta_d) -P_0 \subset \bigcup_{i=1}^d \left(\mathcal{B}(\delta_d)-P_i\right)/z_i.
\end{eqnarray}
\end{proposition}

\begin{proof}[Proof of Proposition \ref{PROPTHM3}]
Let $m$ an integer $\geq 2$ be such that the Lonely Runner Conjecture with Free Starting Positions in Obstruction Form does not hold for $m$, let $d:=m-1$, and let $(a,b) \in \{0,\ldots, d\}$ be fixed. By assumption (see Proposition \ref{PROPLRCFSPequivLRCFSPOF}), there are $\bar{P}=(\bar{P}_1,\ldots,\bar{P}_m)\in [0,1)^m$ and distinct real numbers $\bar{w}_1, \ldots, \bar{w}_{m} >0$ such that
\begin{eqnarray}\label{EQ5fev}
\bar{P} +t  \left(\bar{w}_1, \ldots, \bar{w}_m\right) \notin \mathcal{K}^m(1/(m+1)) \qquad \forall t > 0.
\end{eqnarray}
Without loss of generality, we may assume that $\bar{w}_1=1$ and $\bar{w}_2, \ldots, \bar{w}_m > 1$. Since $\bar{w}_{a+1}\neq \bar{w}_{b+1}$, there is $\bar{t}>0$ such that 
\[
\left\{ \bar{P}_{a+1} + \bar{t} w_{a+1} \right\} = \left\{ \bar{P}_{b+1} + \bar{t} w_{b+1} \right\}.
\]
Define $P_0, P_1, \ldots, P_d \in [0,1)^{d+1}$ and $(z_1,\ldots, z_d) \in (1,\infty)^d$ by
\[
P_i := \bar{P}_{i+1} + \bar{t} \bar{w}_{i+1} \quad \mbox{and} \quad z_i := \bar{w}_{i+1} \qquad \forall i=0, \ldots,d.
\]
Then (\ref{23jan1}) is satisfied and moreover (\ref{EQ5fev}) means that for every $s\geq 0$, there is $i\in \{0,\ldots,d\}$ such that
\[
P_i + s z_i = \bar{P}_i + (\bar{t}+s) \bar{w}_{i+1} \notin \mathcal{K}(1/(m+1)),
\]
which is equivalent to saying that for every $s\geq 0$, there is $i\in \{0,\ldots,d\}$ such that
\[
s \notin \left( \mathcal{K}(\delta_d) - P_i \right)/z_i.
\]
This property implies that if some $s\geq 0$ does belong to  $\mathcal{K}(\delta_d) - P_0=(\mathcal{K}(\delta_d) - P_0)/z_0$, then it has to be in some of the sets $(\mathcal{B}(\delta_d)-P_1)/z_1, \ldots, (\mathcal{B}(\delta_d)-P_d)/z_d$, which gives (\ref{23jan3}).
\end{proof}
\subsection{Case $n=3$} 

We argue by contradiction and apply Proposition \ref{PROPTHM3} with $d=1$, $a=0$, $b=1$ and $N=1$. If the Lonely Runner Conjecture with Free Starting Positions in Obstruction Form does not hold for $m=2$, then there are $P \in (0,1)$ and $z>1$ such that
\[
\mathcal{K}_{0}(1/3) -P \subset \mathcal{K}(1/3) -P \subset \left(\mathcal{B}(1/3)-P\right)/z.
\] 
Hence, there is $k \in \N$ such that $\mathcal{K}_0(1/3)-P\subset (\mathcal{B}_{k}(1/3)-P)/z$, that is,
\[
\frac{3k-1-3P}{z} <1-3P < 2-3P < \frac{3k+1-3P}{z}.
\]
We infer that $z<2$ and $k<((1-3P)z+(1+3P))/3<2/3$, so that $k=0$. Then the above inequalities give  
\[
  \left(\frac{z-1}{z}\right) -1   <  1-3P \left(\frac{z-1}{z}\right) < 2-3P \left(\frac{z-1}{z}\right)   < 1  -  \left(\frac{z-1}{z}\right),
\]
or by setting $\rho=(z-1)/z\in(0,1)$ and $Q=\rho P\in (0,\rho)$,
\[
0< \rho < 1, \quad 0<Q<\rho, \quad 0  < 2 -\rho - 3 Q, \quad 0  < -1-\rho +3Q.
\]
The above system has no solution because $Q<\rho$ and the fourth inequality inequality yield $\rho>1/2$ while the sum of the third and fourth inequalities gives $\rho<1/2$.

\subsection{Case $n=4$} 

We have $\delta_2=1/4$ and for every $k\in \N$, $\mathcal{K}_k=[\mathcal{K}_k^-,\mathcal{K}_k^+] :=\mathcal{K}_k (\delta_2)=[(4k+1)/4,(4k+3)/4]$,  $\mathcal{B}_k=(\mathcal{B}_k^-,\mathcal{B}_k^+) :=\mathcal{B}_k (\delta_2)=[(4k-1)/4,(4k+1)/4]$. We argue by contradiction and apply Proposition \ref{PROPTHM3} with $d=2$, $a=0$, $b=1$. If the Lonely Runner Conjecture with Free Starting Positions in Obstruction Form does not hold for $m=3$, then there are $P_0,P_1, P_2 \in [0,1)$ and $z_2 > z_1 >1$ such that
\[
\mathcal{K} -P_0 \subset \bigcup_{i=1}^2 \left(\mathcal{B}-P_i\right)/z_i \quad \mbox{with} \quad P_0=P_1.
\] 
Note that the equality $P_0=P_1$ won't be used in the sequel. For $i=1,2$, the set $(\mathcal{B}-P_i)/z_i$ is the union of disjoint open intervals $\tilde{\mathcal{B}}^i_k=(\tilde{\mathcal{B}}^{i,-}_k,\tilde{\mathcal{B}}^{i,+}_k):=(\mathcal{B}_k-P_i)/z_i$ of length $\b_i:=1/(2z_i)$ and its complement is the union of closed disjoint intervals of length $\k_i:=1/(2z_i)$ and moreover the length of the sets $[\tilde{\mathcal{B}}^{1,+}_k,\tilde{\mathcal{B}}^{1,-}_{k+1}]$ is strictly less than the length of the sets $\tilde{\mathcal{B}}^2_l$. Therefore, each set $\tilde{\mathcal{K}}_l=[\tilde{\mathcal{K}}_l^-,\tilde{\mathcal{K}}_l^+]:=\mathcal{K}_l-P_0$, with $l\in \N$, can be covered in three ways; there are $k_1,k_2,k_2' \in \N$ such that we have either
\[
\tilde{\mathcal{K}}_l \subset \tilde{\mathcal{B}}^1_{k_1}\cup \tilde{\mathcal{B}}^2_{k_2} \quad \mbox{with} \quad   \tilde{\mathcal{B}}^{1,-}_{k_1}  < l+\frac{1}{4} <  \tilde{\mathcal{B}}^{2,-}_{k_2} <  \tilde{\mathcal{B}}^{1,+}_{k_1} < l+\frac{3}{4} <  \tilde{\mathcal{B}}^{2,+}_{k_2},
\]
or
\[
\tilde{\mathcal{K}}_l \subset \tilde{\mathcal{B}}^2_{k_2}\cup \tilde{\mathcal{B}}^1_{k_1} \quad \mbox{with} \quad   \tilde{\mathcal{B}}^{2,-}_{k_2}  < l+\frac{1}{4} <  \tilde{\mathcal{B}}^{1,-}_{k_1} <  \tilde{\mathcal{B}}^{2,+}_{k_2} < l+\frac{3}{4} <  \tilde{\mathcal{B}}^{1,+}_{k_1}, 
\]
or 
\[
\tilde{\mathcal{K}}_l  \subset \tilde{\mathcal{B}}^2_{k_2}\cup \tilde{\mathcal{B}}^1_{k_1}  \cup \tilde{\mathcal{B}}^2_{k_2'}\quad \mbox{with} \quad   \tilde{\mathcal{B}}^{2,-}_{k_2}  < l+\frac{1}{4} < \tilde{\mathcal{B}}^{1,-}_{k_1} <  \tilde{\mathcal{B}}^{2,+}_{k_2} < \tilde{\mathcal{B}}^{2,-}_{k_2'} <  \tilde{\mathcal{B}}^{1,+}_{k_1} < l+\frac{3}{4} < \tilde{\mathcal{B}}^{2,+}_{k_2'}.
\]
In the first case, we say that we have a chain $\langle 12\vert k_1k_2\rangle_l$, in the second case a chain $\langle 21\vert k_2k_1\rangle_l$ and in the third case a chain $\langle 212\vert k_2k_1k_2'\rangle_l$,
In any case, for each $l\in \N$ there is a unique $k_1^l\in \N$ such that $\tilde{\mathcal{K}}_l\cap \tilde{\mathcal{B}}^1_{k_1^l}\neq \emptyset$ and in addition, the following result holds:

\begin{lemma}\label{LEM4fev}
For every $l\in \N$, we have $k_1^{l+1}=k_1^l+1$.
\end{lemma}

\begin{proof}[Proof of Lemma \ref{LEM4fev}]
Let $l\in \N$ be fixed, we need to check all possible coverings of $\tilde{\mathcal{K}}_l$ and $\tilde{\mathcal{K}}_{l+1}$.\\

\noindent Case 1: $\langle 12\vert k_1^lk_2^l\rangle_l$ for some $k_2^l\in \N$.\\
We have 
\begin{eqnarray}\label{EQ5fev1}
 \tilde{\mathcal{B}}^{1,-}_{k_1^l}  < l+\frac{1}{4} <  \tilde{\mathcal{B}}^{2,-}_{k_2^l} <  \tilde{\mathcal{B}}^{1,+}_{k_1^l} < l+\frac{3}{4} <  \tilde{\mathcal{B}}^{2,+}_{k_2^l}.
\end{eqnarray}
\noindent Subcase 1.1:  $\langle 12\vert k_1^{l+1}k_2^{l+1}\rangle_{l+1}$ for some $k_2^{l+1}\in \N$.\\
We have 
\[
 \tilde{\mathcal{B}}^{1,-}_{k_1^{l+1}}  < l+\frac{5}{4} <  \tilde{\mathcal{B}}^{2,-}_{k_2^{l+1}} <  \tilde{\mathcal{B}}^{1,+}_{k_1^{l+1}} < l+\frac{7}{4} <  \tilde{\mathcal{B}}^{2,+}_{k_2^{l+1}}.
\]
Since $\tilde{\mathcal{K}}_l\cap \tilde{\mathcal{B}}^1_{k_1^{l}+1}=\emptyset$ and $\tilde{\mathcal{B}}^{1,-}_{k_1^{l}+1}\leq  \tilde{\mathcal{B}}^{1,-}_{k_1^{l+1}} $ , both $ \tilde{\mathcal{B}}^{1,-}_{k_1^{l+1}}$ and $\tilde{\mathcal{B}}^{1,-}_{k_1^{l+1}} $ belong to the interval $[l+3/4,l+5/4]$. So, if $k_1^{l+1}\geq k_1^l+2$, then we have $1/z_1 \leq \tilde{\mathcal{B}}^{1,-}_{k_1^{l+1}} -\tilde{\mathcal{B}}^{1,-}_{k_1^{l+1}} \leq 1/2$, which implies that $z_1\geq 2$ and $\b_1\leq 1/4$, $\b_2<\b_1\leq 1/4$. This contradicts the fact that the union of $\tilde{\mathcal{B}}^1_{k_1^l}$ and $\tilde{\mathcal{B}}^2_{k_2^l}$ covers $\tilde{\mathcal{K}}_l$ of length $1/2$. In conclusion, we have $k_1^{l+1}= k_1^l+1$.\\

\noindent Subcase 1.2: $\langle 21\vert k_2^{l+1}k_1^{l+1}\rangle_{l+1}$ for some $k_2^{l+1}\in \N$.\\
In fact this case cannot occur. We have 
 \begin{eqnarray}\label{EQ5fev2}
 \tilde{\mathcal{B}}^{2,-}_{k_2^{l+1}}  < l+\frac{5}{4} <  \tilde{\mathcal{B}}^{1,-}_{k_1^{l+1}} <  \tilde{\mathcal{B}}^{2,+}_{k_2^{l+1}} < l+\frac{7}{4} <  \tilde{\mathcal{B}}^{1,+}_{k_1^{l+1}}.
 \end{eqnarray}
 We claim that $k_1^{l+1}=k_i^l+2$. As a matter of fact, on the one hand (\ref{EQ5fev1})-(\ref{EQ5fev2}) give $ \tilde{\mathcal{B}}^{1,-}_{k_1^{l+1}} - \tilde{\mathcal{B}}^{1,+}_{k_1^l}>1/2$ which implies $k_1^{l+1}-k_1^l\geq 2$, and on the other hand if $k_1^{l+1}\geq k_i^l+3$, then both $\tilde{\mathcal{B}}^{1}_{k_1^{l}+1}, \tilde{\mathcal{B}}^{1}_{k_1^{l}+2}$ must be contained in   the interval $[l+3/4,l+5/4]$ so that $z_1>2$ which contradicts, as in Subcase 1.1, the fact that the union of $\tilde{\mathcal{B}}^1_{k_1^l}$ and $\tilde{\mathcal{B}}^2_{k_2^l}$ covers $\tilde{\mathcal{K}}_l$. Now, set $s:=k_{2}^{l+1}-k_2^l$ and define $\delta_1, \delta_2, \nu_1, \nu_2, d_1, d_2>0$ by
\[
\delta_1 := \tilde{\mathcal{B}}_{k_1^l}^{1,-} - \tilde{\mathcal{K}}_{l-1}^+= \tilde{\mathcal{B}}_{k_1^l}^{1,-} - \left( l-\frac{1}{4}\right), \quad \delta_2 := \tilde{\mathcal{K}}_{l+2}^- - \tilde{\mathcal{B}}_{k_1^{l+1}}^{1,+} =  \left( l+\frac{9}{4}\right) - \tilde{\mathcal{B}}_{k_1^{l+1}}^{1,+},
\]
\[
\nu_1 := \tilde{\mathcal{B}}_{k_2^{l}}^{2,+} - \mathcal{K}_{l}^+ = \tilde{\mathcal{B}}_{k_2^{l}}^{2,+} - \left( l+\frac{3}{4}\right), \quad \nu_2 := \mathcal{K}_{l+1}^- - \tilde{\mathcal{B}}_{k_2^{l+1}}^{2,-} =  \left( l+\frac{5}{4}\right) - \tilde{\mathcal{B}}_{k_2^{l+1}}^{2,-} 
\]
and
\[
d_1 := \mathcal{K}_{l}^+ -  \tilde{\mathcal{B}}_{k_1^l}^{1,+} = \left( l+\frac{3}{4}\right) -  \tilde{\mathcal{B}}_{k_1^l}^{1,+}, \quad d_2 := \tilde{\mathcal{B}}_{k_1^{l+1}}^{1,-} - \mathcal{K}_{l+1}^- = \tilde{\mathcal{B}}_{k_1^{l+1}}^{1,-} - \left( l+\frac{5}{4}\right).
\]
Since $k_{1}^{l+1}=k_l+2$, we have 
\begin{eqnarray}\label{EQ16fev1}
\delta_1 + \delta_2 = \frac{5}{2} +  \tilde{\mathcal{B}}_{k_1^l}^{1,-}  - \tilde{\mathcal{B}}_{k_1^{l+1}}^{1,+}  = \frac{5}{2} - \frac{5}{2z_1} 
\end{eqnarray}
and by construction of $s, \nu_1,\nu_2, d_1,d_2$ we also have
\begin{eqnarray}\label{EQ16fev0}
 \frac{s-\frac{1}{2} }{z_2} + \nu_1 + \nu_2 = \left( \tilde{\mathcal{B}}_{k_2^{l+1}}^{2,-}  - \tilde{\mathcal{B}}_{k_2^{l}}^{2,+} \right) + \left( \tilde{\mathcal{B}}_{k_2^{l}}^{2,+} - \left( l+\frac{3}{4}\right)\right) + \left( \left( l+\frac{5}{4}\right) - \tilde{\mathcal{B}}_{k_2^{l+1}}^{2,-} \right) = \frac{1}{2},
\end{eqnarray}
\begin{eqnarray}\label{EQ16fev2}
d_1+d_2 + \frac{s-\frac{1}{2} }{z_2} + \nu_1+\nu_2 = \tilde{\mathcal{B}}_{k_1^{l+1}}^{1,-} -  \tilde{\mathcal{B}}_{k_1^l}^{1,+} < \tilde{\mathcal{B}}^{2,+}_{k_2^{l+1}} -   \tilde{\mathcal{B}}^{2,-}_{k_2^l} = \frac{s+\frac{1}{2} }{z_2} 
\end{eqnarray}
and
\[
d_i + \frac{1}{2z_1} + \delta_i = 1 \qquad \forall i=1,2.
\]
Hence by summing up the two last inequalities, we obtain
\[
d_1 + d_2 + \frac{1}{z_1}  +\delta_1+\delta_2= 2,
\]
which by (\ref{EQ16fev1}) gives
\[
d_1 + d_2 = \frac{3}{2z_1} - \frac{1}{2}
\]
and by (\ref{EQ16fev2}) implies
\[
\frac{1}{z_2} -\nu_1 - \nu_2 > \frac{3}{2z_1} - \frac{1}{2}.
\]
Therefore by (\ref{EQ16fev0}), we infer that
\[
 \frac{s+\frac{1}{2} }{z_2}  > \frac{3}{2z_1},
\]
which implies that $s\geq 2$. But we have (because $k_1^{l+1}=k_i^l+2$)
\[
\frac{s+\frac{1}{2} }{z_2} = \tilde{\mathcal{B}}_{k_2^{l+1}}^{2,+} - \tilde{\mathcal{B}}_{k_2^{l}}^{2,-} > \tilde{\mathcal{B}}_{k_1^{l+1}}^{1,-} - \tilde{\mathcal{B}}_{k_1^{l}}^{1,+} = \frac{3}{2z_1} 
\]
and
\[
\frac{s-\frac{1}{2} }{z_2} = \tilde{\mathcal{B}}_{k_2^{l+1}}^{2,-} - \tilde{\mathcal{B}}_{k_2^{l}}^{2,+} < \frac{1}{2}
\]
which give respectively $z_2<(2s+1)z_1/3$ and  $z_2>2s-1$ and $z_2<(2s+1)z_1/3$, so that since $s\geq 2$, we have $z_1 >3(2s-1)/(2s+1)\geq 9/5$, which is impossible because $5/(2z_1)>3/2$, that is $z_1<5/3<9/5$.\\

\noindent Subcase 1.3: $\langle 212\vert k_2^{l+1}k_1^{l+1}k_2'\rangle_{l+1}$ for some $k_2^{l+1},k_2'\in \N$.\\
In fact this case cannot occur and the proof is moreorless the same. We have 
\begin{eqnarray}\label{EQ16fev20}
 \tilde{\mathcal{B}}^{2,-}_{k_2^{l+1}}  < l+\frac{5}{4} <  \tilde{\mathcal{B}}^{1,-}_{k_1^{l+1}} <  \tilde{\mathcal{B}}^{2,+}_{k_2^{l+1}} <  \tilde{\mathcal{B}}^{2,-}_{k_2'} <  \tilde{\mathcal{B}}^{1,+}_{k_1^{l+1}}< l+\frac{7}{4} < \tilde{\mathcal{B}}^{2,+}_{k_2'} .
 \end{eqnarray}
 By the same reasoning as Subcase 1.2, we show that $k_1^{l+1}=k_i^l+2$, then we define $s, \delta_1, \delta_2, \nu_1,\nu_2$ exactly as in Subcase 1.2 and follows the same proof.\\
 
\noindent Case 2: $\langle 21\vert k_2^lk_1^l\rangle_l$ for some $k_2^l\in \N$.\\
We have 
\[ \tilde{\mathcal{B}}^{2,-}_{k_2^l}  < l+\frac{1}{4} <  \tilde{\mathcal{B}}^{1,-}_{k_1^l} <  \tilde{\mathcal{B}}^{2,+}_{k_2^l} < l+\frac{3}{4} <  \tilde{\mathcal{B}}^{1,+}_{k_1^l}.
\]
\noindent Subcase 2.1:  $\langle 12\vert k_1^{l+1}k_2^{l+1}\rangle_{l+1}$ for some $k_2^{l+1}\in \N$.\\
We have 
\[
\tilde{\mathcal{B}}^{1,-}_{k_1^{l+1}}  < l+\frac{5}{4} <  \tilde{\mathcal{B}}^{2,-}_{k_2^{l+1}} <  \tilde{\mathcal{B}}^{1,+}_{k_1^{l+1}} < l+\frac{7}{4} <  \tilde{\mathcal{B}}^{2,+}_{k_2^{l+1}}.
\]
If $k_1^{l+1}\geq k_1^l+2$, then we have $ \tilde{\mathcal{B}}^{1,-}_{k_1^{l+1}} -\tilde{\mathcal{B}}^{1,+}_{k_1^{l}} < 1/2$  which implies that $z_1> 3$ and $\b_2<\b_1\leq 1/6$. This contradicts the fact that the union of $\tilde{\mathcal{B}}^1_{k_1^l}$ and $\tilde{\mathcal{B}}^2_{k_2^l}$ covers $\tilde{\mathcal{K}}_l$ of length $1/2$. In conclusion, we have $k_1^{l+1}= k_1^l+1$.\\

\noindent Subcase 2.2: $\langle 21\vert k_2^{l+1}k_1^{l+1}\rangle_{l+1}$ for some $k_2^{l+1}\in \N$.\\
By symmetry, this case follows from the proof of Subcase 1.1.\\

\noindent Subcase 2.3: $\langle 212\vert k_2^{l+1}k_1^{l+1}k_2'\rangle_{l+1}$ for some $k_2^{l+1},k_2'\in \N$.\\
 We have 
 \[
 \tilde{\mathcal{B}}^{2,-}_{k_2^{l+1}}  < l+\frac{5}{4} < \tilde{\mathcal{B}}^{1,-}_{k_1^{l+1}} <  \tilde{\mathcal{B}}^{2,+}_{k_2^{l+1}} < \tilde{\mathcal{B}}^{2,-}_{k_2'} <  \tilde{\mathcal{B}}^{1,+}_{k_1^{l+1}} < l+\frac{7}{4} < \tilde{\mathcal{B}}^{2,+}_{k_2'} .
 \]
 If $k_1^{l+1}\geq k_1^l+2$, then both $ \tilde{\mathcal{B}}^{1,+}_{k_1^{l}}$ and $\tilde{\mathcal{B}}^{1,+}_{k_1^{l}+1} $ belong to the interval $[l+3/4,l+5/4]$. So, we have $z_2>2$ and we obtain a contradiction as before. In conclusion, we have $k_1^{l+1}= k_1^l+1$.\\
 
\noindent Case 3: $\langle 212\vert k_2^lk_1^lk\rangle_l$ for some $k_2^l,k\in \N$.\\

\noindent Subcase 3.1:  $\langle 12\vert k_1^{l+1}k_2^{l+1}\rangle_{l+1}$ for some $k_2^{l+1}\in \N$.\\
By symmetry, this case follows from the proof of Subcase 2.3.\\

\noindent Subcase 3.2: $\langle 21\vert k_2^{l+1}k_1^{l+1}\rangle_{l+1}$ for some $k_2^{l+1}\in \N$.\\
By symmetry, this case follows from the proof of Subcase 1.3.\\

\noindent Subcase 3.3: $\langle 212\vert k_2^{l+1}k_1^{l+1}k_2'\rangle_{l+1}$ for some $k_2^{l+1},k_2'\in \N$.\\
In fact this case cannot occur and the proof follows the same line as the proof of Subcases 1.2 and 1.3.

\end{proof}

By Lemma \ref{LEM4fev}, we have 
\begin{eqnarray}\label{EQ6fev}
l+ \frac{1}{4} < \tilde{\mathcal{B}}_{k_1^l}^{1,+} < l + \frac{3}{4} \qquad \forall l \in \N.
\end{eqnarray}
 Define the sequence $\{u_k\}_{k\in \N}$ by 
$$
u_k := \tilde{\mathcal{B}}_{k_1^0+k}^{1,+} - \left(k+\frac{1}{4} \right) =  \tilde{\mathcal{B}}_{k_1^0}^{1,+} - \frac{1}{4} + k \left(\frac{1}{z_1}- 1 \right) \qquad \forall k\in \N.
$$
By construction, we have $u_0>0$ and 
\[
u_{k+1}-u_k = \frac{1}{z_1}- 1 <0 \qquad \forall k\in \N. 
\]
Therefore, $\{u_k\}_{k\in \N}$ is decreasing and there is $\bar{k}\in \N^*$ such that $u_{\bar{k}}<0$, that is,
\[
\tilde{\mathcal{B}}_{k_1^k}^{1,+} = \tilde{\mathcal{B}}_{k_1^0+\bar{k}}^{1,+} < k+\frac{1}{4}.
\]
This contradicts (\ref{EQ6fev}) and concludes the proof of the case $n=4$.

\section{Upper bound on the $\mathcal{K}_0(\delta)$-measures}\label{SECK0Measure}
 
 The purpose of this section is to give a proof of the upper bound that we gave in the introduction on the $\mathcal{K}_0(\delta)$-measures defined by 
 \[
\mathcal{M}_0^{\delta}(z) := \mathcal{L}^1 \left( \mathcal{B}(\delta)/z \cap \mathcal{K}_0(\delta) \right) \qquad \forall z\geq 1, \, \forall \delta\in (0,1/2].
\]
 
 \begin{figure}[H]
\begin{center}
\includegraphics[width=8cm]{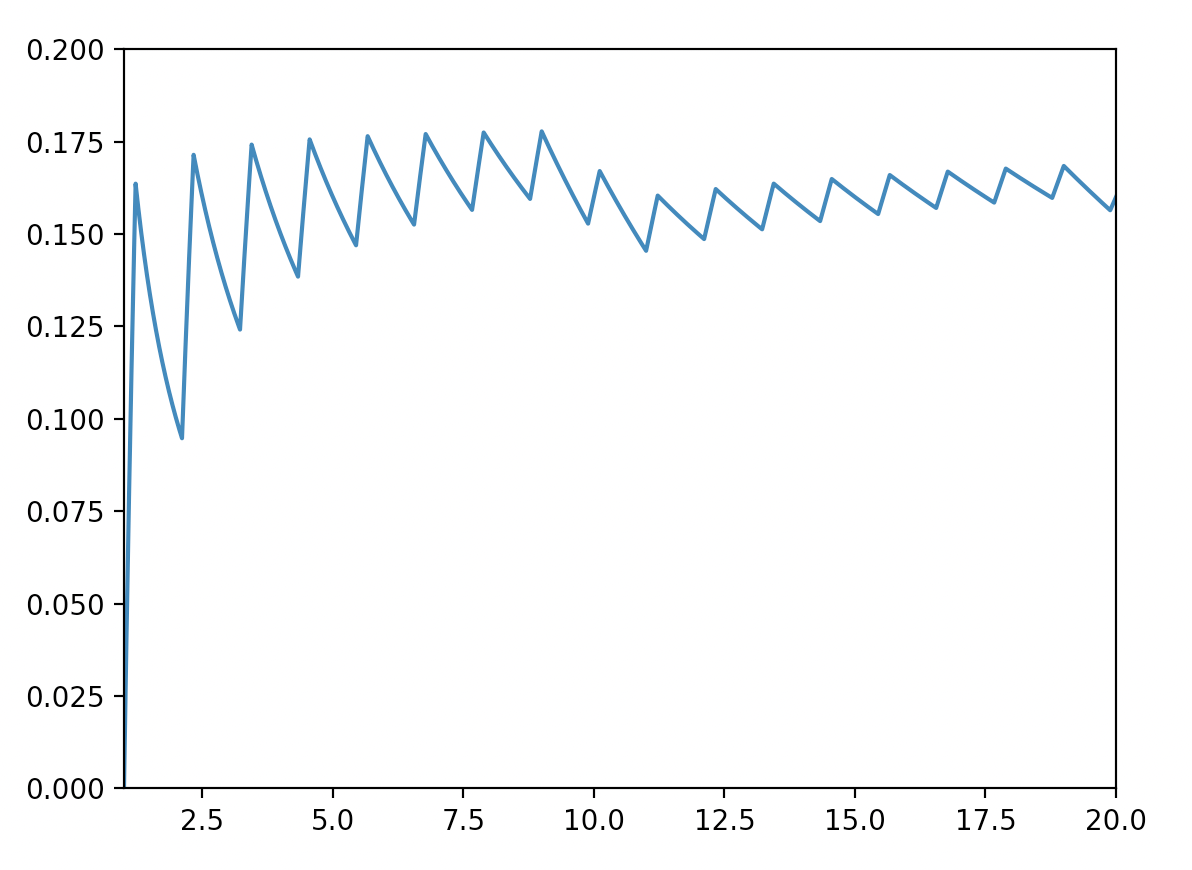}
\caption{The function $z\geq 1 \mapsto \mathcal{M}_0^{\delta}(z)$ is bounded from above by $2\delta(1-2\delta)/(1-\delta)$, here with $\delta =1/10$ \label{figMeasure}}
\end{center}
\end{figure}

We prove the following result (see Figure \ref{figMeasure}):

 \begin{proposition}\label{PROPMeasure}
 For every $\delta \in (0,1/3]$, we have
 \[
 \mathcal{M}_0^{\delta}(z) \leq  \frac{2\delta(1-2\delta)}{1-\delta} \qquad \forall z \geq 1.
 \]
 Moreover, for every $\delta=1/k$ with $k\in \N \cap [3,\infty)$, $z\geq 1$ defined by $z:=k-1$ satisfies $\mathcal{M}_0^{\delta}(z)=2\delta(1-2\delta)/(1-\delta)$.
 \end{proposition}
 
 \begin{proof}[Proof of Proposition \ref{PROPMeasure}]
Let $\delta\in (0,1/2)$ be fixed, for every $z\geq 1$, we define $a^{\delta}(z), b^{\delta}(z)\in \N$ by
\[
a^{\delta}(z) := \min \Bigl\{k\in \N \, \vert \, \mathcal{B}_k^+(\delta)/z \geq \delta\Bigr\}=  \left\lceil \delta (z-1) \right\rceil.
\]
and
\[
 b^{\delta}(z) := \max \Bigl\{k\in \N \, \vert \, \mathcal{B}_k^-(\delta)/z \leq 1-\delta\Bigr\}=  \left\lfloor (1-\delta)z+\delta \right\rfloor.
\]
We fix $\bar{z}> 1$, we set $\bar{a}:=a^{\delta}(\bar{z})\geq 1, \bar{b}:=b^{\delta}(\bar{z})\geq 1$ and we distinguish several cases:\\

\noindent Case 1: $\mathcal{B}_1^-(\delta)/\bar{z} \geq \delta$.\\
Then $\bar{z}\leq (1-\delta)/\delta$ and we have
\[
\mathcal{M}^{\delta}_0 (\bar{z}) = \left\{ \begin{array}{ccl}
 \frac{2\bar{b}\delta}{\bar{z}}  & \mbox{ if } & \mathcal{B}_{\bar{b}}^+(\delta)/\bar{z} \leq 1-\delta\\
 (\bar{b} -1) \frac{2\delta}{\bar{z}} + (1-\delta) - \left(\frac{ \bar{b}-\delta}{\bar{z}}\right) & \mbox{ if } & \mathcal{B}_{\bar{b}}^+(\delta)/\bar{z} > 1-\delta.
\end{array}
\right.
\]
We distinguish two subcases:\\

\noindent Case 1.1: $\mathcal{B}_{\bar{b}}^+(\delta)/\bar{z} \leq 1-\delta$. \\
Then we have $(\bar{b}+\delta)/(1-\delta) \leq \bar{z} \leq (1-\delta)/\delta$ which implies $\bar{b}\leq (1-\delta)^2/\delta - \delta = (1-2\delta)/\delta$. Thus, since the function $r\geq 0 \mapsto r/(r+\delta)$ is increasing, we obtain 
\[
\mathcal{M}_0^{\delta} (\bar{z})  =  \frac{2\bar{b}\delta}{\bar{z}}  \leq \frac{2\bar{b}\delta(1-\delta)}{\bar{b}+\delta} \leq 2\delta(1-\delta) \left( \frac{\frac{1-2\delta}{\delta}}{\frac{1-2\delta}{\delta}+\delta} \right) = \frac{2\delta(1-2\delta)}{1-\delta}.
\]

\noindent Case 1.2: $\mathcal{B}_{\bar{b}}^+(\delta)/\bar{z} > 1-\delta$.\\
 Then we have $\bar{z} < (\bar{b}+\delta)/(1-\delta)$. The function 
\[
f:z\mapsto  (\bar{b} -1) \frac{2\delta}{z} + 1-\delta - \left(\frac{ \bar{b}-\delta}{z}\right)  = 1-\delta + \frac{(2\delta-1)\bar{b}-\delta}{z}
\]
 is increasing on $[1,\infty)$ and satisfies $f((\bar{b}+\delta)/(1-\delta))=2\bar{b}\delta(1-\delta)/(\bar{b}+\delta)$. Thus, if $z:=(\bar{b}+\delta)/(1-\delta)$ satisfies $z\leq (1-\delta)/\delta \Leftrightarrow \mathcal{B}_1^-(\delta)/z \geq \delta$, then we have $b^{\delta}(z)=\bar{b}$ and $\mathcal{M}_0^{\delta}(z)=2\bar{b}\delta/z$ (because $ \mathcal{B}_{\bar{b}}^-(\delta)/z \leq  \mathcal{B}_{\bar{b}}^+(\delta)/z=(1-\delta)$), so Case 1.1 applies to $z$ and gives
 \[
 \mathcal{M}_0^{\delta} (\bar{z}) \leq f(z) = \frac{2b\delta}{z} =\mathcal{M}_0^{\delta}(z) \leq  \frac{2\delta(1-2\delta)}{1-\delta}.
 \]
If however, we have $z=(\bar{b}+\delta)/(1-\delta)> (1-\delta)/\delta$ (which implies $\bar{b}> (1-\delta)^2/\delta -\delta$), then $\hat{z}:=(1-\delta)/\delta \in [\bar{z},z]$ satisfies $\mathcal{M}_0^{\delta} (\bar{z}) \leq \mathcal{M}_0^{\delta} (\hat{z})$, $b^{\delta}(\hat{z})=\bar{b}$ (because $\mathcal{B}_{\bar{b}}^-(\delta)/\hat{z} \leq  \mathcal{B}_{\bar{b}}^-(\delta)/\bar{z}\leq (1-\delta)$ and $\mathcal{B}_{\bar{b}+1}^-(\delta)/\hat{z} > \mathcal{B}_{\bar{b}}^+(\delta)/\hat{z}\geq \mathcal{B}_{\bar{b}}^+(\delta)/z=1-\delta$) and (note that $2\delta-1<0$)
\begin{eqnarray*}
\mathcal{M}_0^{\delta} (\hat{z}) = 1-\delta + \frac{(2\delta-1)\bar{b}-\delta}{\hat{z}} & = &  1-\delta - \frac{\delta^2}{1-\delta} + \frac{\delta(2\delta-1)}{1-\delta} \bar{b} \\
& \leq & \frac{1-2\delta}{1-\delta} +  \frac{\delta(2\delta-1)}{1-\delta}  \left(\frac{(1-\delta)^2}{\delta} - \delta \right) = \frac{2\delta(1-2\delta)}{1-\delta}.
 \end{eqnarray*}
 
\noindent Case 2: $\delta \in \mathcal{B}_1(\delta)/\bar{z}$ \\
Then we have 
\[
 \frac{1-\delta}{\delta} < \bar{z} <  \frac{1+\delta}{\delta}.
 \]
 Again we distinguish two cases:\\
 
\noindent Case 2.1: $1-\delta \in \mathcal{B}_{\bar{b}}(\delta)/\bar{z}$.\\
In this case, we have
  \[
 \frac{\bar{b}-\delta}{1-\delta} < \bar{z} <  \frac{\bar{b}+\delta}{1-\delta}
 \]
 and
 \[
 \mathcal{M}_0^{\delta}(\bar{z}) = (\bar{b} -2) \frac{2\delta}{\bar{z}} + (1-\delta) - \left(\frac{ \bar{b}-\delta}{\bar{z}}\right) +  \left(\frac{ 1+\delta}{\bar{z}}\right) - \delta = \frac{(1-2\delta) (1+\bar{z}-\bar{b})}{\bar{z}}.
 \]
 But $(\bar{b}+\delta)/(1-\delta)>\bar{z}$ gives $\bar{b} >(1-\delta)\bar{z}-\delta$ which implies $(1+\bar{z}-\bar{b})/\bar{z}<(1+\delta)/\bar{z}+\delta$ which is $<2\delta/(1-\delta)$ because $\bar{z}>(1-\delta)/\delta$.\\

\noindent  Case 2.2: $\mathcal{B}_{\bar{b}}(\delta)^+/\bar{z}\leq 1-\delta$.\\
 Then we have $\bar{z}>(1-\delta)/\delta$, $\bar{z}\geq (\bar{b}+\delta)/(1-\delta)$ and
  \[
 \mathcal{M}_0^{\delta}(\bar{z}) = (\bar{b} -1) \frac{2\delta}{\bar{z}} + \left(\frac{ 1+\delta}{\bar{z}}\right) - \delta  = \frac{\delta(2\bar{b} -1)}{\bar{z}} - \delta.
 \]
 The function 
\[
g:z\mapsto \frac{\delta(2\bar{b} -1)}{z} - \delta
\]
is decreasing on $[1,\infty)$ so we can consider the least $t\geq 0$ such that $z$ defined by $z=\bar{z}-t$ satisfies $z=(1-\delta)/\delta$ or $z=(\bar{b}+\delta)/(1-\delta)$. If $z=(1-\delta)/\delta$ we get the result by Case 1.1 (note that $g(\bar{z})\leq g(z)$). Otherwise, we have $z=(\bar{b}+\delta)/(1-\delta)$ and $z>(1-\delta)/\delta$. We get the result by Case 2.1 (we have $1-\delta =\mathcal{B}_{\bar{b}}^+/\bar{z}$).\\
  
 \noindent Case 3: $\mathcal{B}_1^+(\delta)/\bar{z}<\delta$. \\
 Then we have $\bar{z}>(1+\delta)/\delta \geq  2(1-\delta)/(1-2\delta)$ (because $\delta \leq 1/3$) and 
  \[
 \mathcal{M}_0^{\delta}(\bar{z}) \leq (\bar{b}-\bar{a}+1) \frac{2\delta}{\bar{z}}, 
 \]
 with 
 \[
 \bar{b} \leq (1-\delta)\bar{z}+\delta \quad \mbox{and} \quad  \bar{a} \geq \delta (\bar{z}-1) +1.
 \]
 So, we obtain
 \begin{eqnarray*}
  \mathcal{M}_0^{\delta}(\bar{z}) & \leq & \frac{2\delta}{\bar{z}} \left( (1-\delta)\bar{z}+\delta - \delta (\bar{z}-1) -1 +1 \right) \\
  & = &  \frac{2\delta}{z} \left( (1-2\delta)\bar{z}+2\delta  \right)\\
  & = & 2\delta (1-2\delta) + \frac{4\delta^2}{\bar{z}}\\
  & \leq  & 2\delta (1-2\delta) + \frac{2\delta^2(1-2\delta)}{1-\delta} = \frac{2\delta(1-2\delta)}{1-\delta}.
 \end{eqnarray*}
 The proof is complete.
 \end{proof}
 
 \appendix
 
\section{Proofs of Auxiliary results}

This section is devoted to several results stated without proofs in the introduction. Our first result is concerned with the equivalence of the Lonely Runner Conjecture and the Lonely Runner Conjecture in Obstruction Form.

\begin{proposition}\label{PROPLRCequivLRCOF}
The Lonely Runner Conjecture  holds true for $n\geq 2$ if and only if for every $w_1, \ldots, w_{n-1} >0$ there is $t>0$ such that $\{tw_i\}$ belongs to $[1/n,1-1/n]$ for all $i=1, \ldots, n-1$.
\end{proposition} 

As we said after the statement of the Lonely Runner Conjecture, the above proposition follows from the following result:

\begin{lemma}\label{LEMLRCequivLRCOF}
For every $v, v' \in \R$, every $t>0$ and every $\delta \in (0,1/2]$, one has that $|\{tv\}-\{tv'\}|$ belongs to $[\delta,1-\delta]$ if and only if   $\{t|v'-v|\}$ belongs to   $[\delta,1-\delta]$.
\end{lemma}

\begin{proof}[Proof of Lemma \ref{LEMLRCequivLRCOF}]
Let $v, v' \in \R$, $t>0$ and $\delta \in (0,1/2]$ be fixed. There are $n,n' \in \Z$ and $\theta, \theta' \in [0,1)$ such that
\[
tv = n + \theta \quad \mbox{and} \quad tv' = n' + \theta'.
\]
Then, on the one hand we have 
\[
|\{tv\}-\{tv'\}|=|\theta-\theta'|
\]
and on the other hand we have 
\[
t (v'-v) = n'-n + \theta' - \theta  \quad \mbox{and} \quad  t (v-v') = n-n' + \theta - \theta'
\]
which shows that $\{t|v'-v|\}=\{\pm (\theta - \theta')\}$. We conclude easily. 
\end{proof}

\begin{proof}[Proof of Proposition \ref{PROPLRCequivLRCOF}]
Assume that the Lonely Runner Conjecture  holds true for some $n\geq 2$ and consider $w_1, \ldots, w_{n-1} >0$. Since some of the $w_i$ could be equal, we have to consider the set 
\[
\mathcal{V} := \bigl\{0\bigr\} \cup \bigl\{w_1,\ldots, w_{n-1} \bigr\}.
\]
If $\sharp (\mathcal{V})=n$, then we write $\mathcal{V} = \{v_1, \ldots, v_n\}$ with $v_1=0$ and apply the Lonely Runner Conjecture with $i=1$ to obtain $t>0$ such that $|\{tv_j\}| \in [1/n,1-1/n]$ for all $j=2, \ldots,n$, which by Lemma \ref{LEMLRCequivLRCOF} implies that $\{tw_j\}$ belongs to $[1/n,1-1/n]$ for all $j=1, \ldots, n-1$. If $\sharp (\mathcal{V})<n$, we write $\mathcal{V} = \{v_1, \ldots, v_r\}$, with $v_1=0$ and $\sharp (\mathcal{V})\geq 2$, and pick $n-r$ distinct real numbers $v_{r+1},\ldots, v_n$ not in $\mathcal{V}$. Then, we apply the Lonely Runner Conjecture with $i=1$ to obtain $t>0$ such that $|\{tv_j\}| \in [1/n,1-1/n]$ for all $j=2, \ldots,n$, which by Lemma \ref{LEMLRCequivLRCOF} implies that $\{tw_j\}$ belongs to $[1/n,1-1/n]$ for all $i=j, \ldots, n-1$. Now, we fix an integer $n\geq 2$ and assume that for every $w_1, \ldots, w_{n-1} >0$ there is $t>0$ such that $\{tw_i\}$ belongs to $[1/n,1-1/n]$ for all $i=1, \ldots, n-1$. Then, given $n$ runners with different constant speeds  $v_1, \ldots, v_n$ and $i\in \{1,\ldots,n\}$, we define $w_1, \ldots, w_{n-1} >0$ such that 
\[
\bigl\{w_1,\ldots, w_{n-1}\bigr\} = \bigl\{ |v_j-v_i| \, \vert \, j \in \{1,\ldots,n\} \setminus \{i\}\bigr\}.
\]
By applying the assumption, we get $t>0$ such that $\{tw_i\}$ belongs to $[1/n,1-1/n]$ for all $j=1, \ldots, n-1$, which thanks to Lemma \ref{LEMLRCequivLRCOF}, means that $|\{tv_j\}-\{tv_i\}| \in [1/n,1-1/n]$ for all $j$ in $\{1, \ldots, n\} \setminus \{i\}$. The proof is complete.
\end{proof}

Our second result refers to the equivalence of the Lonely Runner Conjecture with Free Starting Positions with the Lonely Runner Conjecture with Free Starting Positions in Obstruction Form. 

\begin{proposition}\label{PROPLRCFSPequivLRCFSPOF}
The Lonely Runner Conjecture with Free Starting Positions holds true for $n\geq 2$ if and only if for every $P_1, \ldots, P_{n-1}\in [0,1)$ and every distinct positive real numbers $w_1, \ldots, w_{n-1}$, there is $t>0$ such that $\{P_i+tw_i\}$ belongs to $[1/n,1-1/n]$ for all $i=1, \ldots, n-1$. As a consequence, the Lonely Runner Conjecture with Free Starting Positions is equivalent to the Lonely Runner Conjecture with Free Starting Positions in Obstruction Form.
\end{proposition} 

\begin{proof}[Proof of Proposition \ref{PROPLRCFSPequivLRCFSPOF}]
As in the proof of Proposition \ref{PROPLRCequivLRCOF}, we start with a preliminary lemma.

\begin{lemma}\label{LEMLRCFSPequivLRCFSPOF}
For every $a, a', v, v' \in \R$, every $t>0$ and every $\delta \in (0,1/2]$, one has that $|\{a+tv\}-\{a'+tv'\}|$ belongs to $[\delta,1-\delta]$ if and only if   $\{ b + t|v-v'|\}$ belongs to $[\delta,1-\delta]$, where $b:=a-a'$ if $|v-v'|=v-v'$ and $b:=a'-a$ if $|v-v'|=v'-v$. 
\end{lemma}
\begin{proof}[Proof of Lemma \ref{LEMLRCFSPequivLRCFSPOF}]
Let $a,a',v, v' \in \R$, $t>0$ and $\delta \in (0,1/2]$ be fixed. There are $n,n' \in \Z$ and $\theta, \theta' \in [0,1)$ such that
\[
a+tv = n + \theta \quad \mbox{and} \quad a'+tv' = n' + \theta'.
\]
Then, on the one hand we have 
\[
|\{a+tv\}-\{a'+tv'\}|=|\theta-\theta'|
\]
and on the other hand we have 
\[
(a-a') + t (v-v') = n-n' + \theta - \theta'  \quad \mbox{and} \quad  (a-a') + t (v'-v) = n-n' + \theta - \theta'
\]
which shows that $\{t|v'-v|\}=\{\pm (\theta - \theta')\}$. We conclude easily. 
\end{proof}

Assume that the Lonely Runner Conjecture with Free Starting Positions holds true for some $n\geq 2$ and consider $P_1,\ldots,P_{n-1} \in [0,1)$ and distinct real numbers $w_1, \ldots, w_{n-1} >0$. The constant speeds of the $n$ runners $x_1,\ldots, x_n:[0,\infty) \rightarrow \R/\Z$ given by 
\[
x_1(t) =0 \quad \mbox{and} \quad x_i(t) = \left\{P_{i-1}+tw_{i-1}\right\} \quad \forall i=2,\ldots, n,
\]
satisfy (\ref{CardSpeeds}) with $i=1$, so there is $t>0$ such that $x_i(t)\in [1/n,1-1/n]$ for all $i=2,\ldots,n$, which proves the result. Conversely, suppose that for every $P_1, \ldots, P_{n-1}\in [0,1)$ and every distinct positive real numbers $w_1, \ldots, w_{n-1}$, there is $t>0$ such that $\{P_i+tw_i\}$ belongs to $[1/n,1-1/n]$ for all $i=1, \ldots, n-1$ and consider $n$ runners $x_1,\ldots, x_n: [0,\infty) \rightarrow \R/\Z$ with constant speeds $v_1, \ldots, v_n$, starting positions $\bar{x}_1, \ldots, \bar{x}_n$ and $i\in \{1,\ldots,n\}$ such that (\ref{CardSpeeds}) is satisfied. Let $w_1, \ldots, w_{n-1} \in \R$ be such that 
\[
\bigl\{w_1,\ldots, w_{n-1}\bigr\} = \bigl\{ |v_j-v_i| \, \vert \, j \in \{1,\ldots,n\} \setminus \{i\}\bigr\}.
\]
By (\ref{CardSpeeds}), all $w_k$ with $k=1,\ldots ,n-1$ are positive and distincts and there is a one-to-one mapping $\sigma :\{1,\ldots, n-1\} \rightarrow \{1,\ldots,n\} \setminus \{i\}\bigr\}$ such that for every $k\in \{1,\ldots, n-1\}$, $w_k=|v_{\sigma(k)}-v_i|$. Then let  $P_1, \ldots, P_{n-1} \in [0,1)$ given by
\[
P_k :=  \left\{ 
\begin{array}{ccl}
\left\{\bar{x}_{\sigma(k)}-\bar{x}_i\right\} & \mbox{ if } & w_k=v_{\sigma(k)}-v_i\\
\left\{  \bar{x}_i- \bar{x}_{\sigma(k)}\right\} & \mbox{ if } & w_k=v_i - v_{\sigma(k)}.
\end{array}
\right.
\]
By applying the assumption, there is $t>0$ such that $\{P_k+tw_k\}$ belongs to $[1/n,1-1/n]$ for all $k=1, \ldots, n-1$, which thanks to Lemma \ref{LEMLRCequivLRCOF}, means that $|\{\bar{x}_j+tv_j\}-\{\bar{x}_i+tv_i\}| \in [1/n,1-1/n]$ for all $j$ in $\{1, \ldots, n\} \setminus \{i\}$. The proof of the first part is complete. The second part follows easily by periodicity of the set $\cup_{k\in \Z^m} \mathcal{K}_k^m(1/(m+1))$ and the fact that any affine line having a direction vector with positive coordinates enters the set $(0,\infty)^m$.
\end{proof}

\section{Proof of Proposition \ref{PROPFeather}}\label{APPPROPFeather}

(i) Let $k=(k_1,\ldots, k_d) \in \N^d$ and $l\in \N$ be fixed.  Let us show that $F_{k,l}^d (\delta)= R_{k,l}^d (\delta)\cap C_k^d (\delta)$. If $(z_1, \ldots, z_d) \in \R^d$ belongs to $F_{k,l}^d$, then there is $\lambda \in \mathcal{K}_l(\delta)$ such that 
\[
(z_1, \ldots, z_d) \in \frac{1}{\lambda} \mathcal{K}_k^d(\delta),
\]
which means that $(z_1, \ldots, z_d) \in C_k^d(\delta)$, and, since $0<l+\delta \leq \lambda \leq l+1-\delta$, that
\[
 \frac{k_i+\delta}{l+1-\delta} \leq z_i \leq \frac{k_i+1-\delta}{l+\delta}\, \forall i=1, \ldots, d,
\]
that is to say $(z_1, \ldots, z_d) \in R_{k,l}^d(\delta)\cap C_k^d(\delta)$. To prove the reverse inclusion, we consider the convex set 
\[
E= \bigcap_{i=1}^d \Bigl\{ (z_1,\ldots, z_d) \in \R^d \, \vert \, z_i \geq k_i+\delta \Bigr\},
\]
whose boundary, denoted by $\partial E$, is given by
\[
\partial E = \bigcup_{i=1}^d \Bigl\{  (z_1,\ldots, z_d) \in (0,\infty)^d \, \vert \, z_i=k_i+\delta, \, z_j\geq k_j+\delta \, \, \forall j=\{1,\ldots,d\} \Bigr\},
\]
and we define the sets $L_k^d(\delta)$ and $U_{k}^d(\delta)$, respectively  the lower and upper faces of $\mathcal{K}_{k}^d(\delta)$, by
\[
L_k^d(\delta):= \bigcup_{i=1}^d L_{k,i}^d(\delta) \quad \mbox{with} \quad L_{k,i}^d(\delta) := \mathcal{K}_k^d(\delta) \cap \Bigl\{ z\in \R^d \, \vert \, z_i=k_i+\delta \Bigr\} \, \, \forall i=1, \ldots, d
\]
and
\[
U_k^d(\delta):= \bigcup_{i=1}^d U_{k,i}^d(\delta) \quad \mbox{with} \quad U_{k,i}^d(\delta) := \mathcal{K}_k^d(\delta) \cap \Bigl\{ z\in \R^d \, \vert \, z_i=k_i+1-\delta \Bigr\} \, \, \forall i=1, \ldots, d.
\]
We have the following result:
 
\begin{lemma}\label{27sept1}
For every $z\in (0,\infty)^d$, there is a unique $\alpha = \alpha(z)>0$ such that $\alpha z\in \partial E$. Moreover the function $R : (0,\infty)^d \rightarrow \partial E$, defined by $R(z):=\alpha(z)z$ for all $z\in (0,\infty)^d$, is continuous and satisfies
\begin{eqnarray}\label{28sept1}
R \left(\mathcal{K}_{k}^d(\delta) \right) = R \left( L_k^d(\delta) \right)= R \left( U_k^d(\delta)\right) = L_k^d(\delta).
\end{eqnarray}
\end{lemma}

\begin{proof}[Proof of Lemma \ref{27sept1}]
For every  $z\in (0,\infty)^d$, the existence of $\alpha>0$ such that $\alpha z\in \partial E$ follows from the fact that $0\notin E$ and $\alpha z \in \mbox{Int}(E)$ for $\alpha>0$ sufficiently large.  If for some $z= (z_1,\ldots, z_d) \in (0,\infty)^d$, there are $\alpha, \alpha'>0$ such that $\alpha z, \alpha' z\in \partial E$, then there are $i,i'\in \{1,\ldots,d\}$ such that $\alpha z_i=k_i+\delta$ and  $\alpha' z_{i'} = k_{i'}+\delta$. But we have $\alpha' z_i\geq k_i+\delta=\alpha z_i$ (because $\alpha' z \in \partial E$) and $\alpha z_{i'}\geq k_{i'}+\delta=\alpha z_i$ (because $\alpha z \in \partial E$). We infer that $\alpha=\alpha'$, which proves the uniqueness. The function $R$ is locally bounded and its graph is closed (in $(0,\infty)^d\times \partial E$), so an easy argument of compactness proves that $R$ is continuous. Since $L_k^d(\delta)\subset \partial E$, we have $R (L_k^d(\delta))=L_k^d(\delta)$. Moreover, for every $z\in \mathcal{K}_l^d(\delta)$, $z_i\geq k_i+\delta$ for all $i=1,\ldots, d$ implies that $\alpha(z)\leq 1$, so that $\alpha(z)z_i \leq k_i+1-\delta$ for any $i=1,\ldots, d$, which shows that $R(z)\in  L_k^d(\delta)$. Thus we have $R(\mathcal{K}_{k}^d(\delta)) = L_k^d(\delta)$ and $R(U_{k}^d(\delta)) \subset R(\mathcal{K}_{k}^d(\delta)) = L_k^d(\delta)$. Define the set $F$ of $(z_1,\ldots,z_d)\in (0,\infty)^d$ such that $z_i\leq k_i+1-\delta$ for all $i=1, \ldots,d$ and denote its boundary by $\partial F$. Then by repeating the above proof, we can prove that for every $z\in (0,\infty)^d$ there is a unique $\beta>0$ such that $\beta z\in \partial F$ and show that the function $S : (0,\infty)^d \rightarrow \partial F$, defined by $S(z):=\beta(z)z$ for all $z\in (0,\infty)^d$, satisfies $S(L_{k}^d(\delta)) \subset S(\mathcal{K}_{k}^d(\delta)) = U_k^d(\delta)$ and $R(S(z))=z$ for all $z\in \partial E$. We infer that $L_k^d(\delta) \subset (R\circ S) (L_k^d(\delta)) \subset R(U_k^d(\delta)) \subset L_k^d(\delta)$. Therefore, we have $R(U_k^d(\delta)) = L_k^d(\delta)$.
\end{proof}

The property (\ref{28sept1}) means that 
\[
C_k^d (\delta) =(0,\infty)  \cdot L_{k}^d(\delta) = (0,\infty)  \cdot U_{k}^d(\delta)
\]
and thus shows that for every $z=(z_1, \ldots, z_d) \in R_{k,l}^d(\delta)\cap C_k^d(\delta)$ there are $\bar{\alpha},\hat{\alpha}>0$ and $\bar{y}=(\bar{y}_1, \ldots, \bar{y}_d), \hat{y}=(\hat{y}_1, \ldots, \hat{y}_d) \in (0,\infty)^d$ such that 
\[
z = \frac{\bar{y}}{\bar{\alpha}} =\frac{\hat{y}}{\hat{\alpha}}, \quad \bar{y} \in L_k^d(\delta) \quad \mbox{and} \quad  \hat{y} \in U_k^d(\delta).
\]
Thus, by definition $L_k^d(\delta) $ and $U_k^d(\delta)$, there are $\bar{i},\hat{i} \in \{1, \ldots, d\}$ such that $\bar{y}_{\bar{i}} =k_{\bar{i}} + \delta, \hat{y}_{\hat{i}} =k_{\hat{i}} + 1- \delta$ which implies, because 
\[
 \frac{k_i+\delta}{l+1-\delta} \leq z_i \leq \frac{k_i+1-\delta}{l +\delta} \qquad \forall i=1, \ldots,d,
\]
 that
\[
\bar{\alpha} = \frac{\bar{y}_{\bar{i}}}{z_{\bar{i}}} \leq \left( k_{\bar{i}} + \delta\right) \left( \frac{l+1-\delta}{ k_{\bar{i}}+\delta}\right) = l +1 - \delta
\]
and
\[
\hat{\alpha} = \frac{\hat{y}_{\hat{i}}}{z_{\hat{i}}} \geq \left( k_{\hat{i}} + 1- \delta\right) \left( \frac{l+\delta}{ k_{\hat{i}}+1-\delta}\right)=l+\delta.
\]
If $\hat{\alpha} \leq \bar{\alpha}$, we infer that $\bar{\alpha} \in [l+\delta,l+1-\delta]=\mathcal{K}_l(\delta)$ and $z = \bar{y}/\bar{\alpha}$ with $\bar{y} \in \mathcal{K}_k(\delta)$, which shows that $z\in F_{k,l}^d(\delta)$. If $\hat{\alpha} > \bar{\alpha}$, then the set $[\bar{\alpha},\hat{\alpha}] \cap [l+\delta,l+1-\delta]$ is non-empty so there are $\lambda \in \mathcal{K}_l(\delta)$ and $s\in [0,1]$ such that $\lambda = s\bar{\alpha}+(1-s)\hat{\alpha}$. Then, by convexity of $\mathcal{K}_k(\delta)$, we have 
\[
\lambda z = \left(  s\bar{\alpha}+(1-s)\hat{\alpha}\right) z = s \bar{\alpha}z + (1-s) \hat{\alpha}z = s \bar{y} + (1-s)\hat{y} \in  \mathcal{K}_k(\delta),
\] 
which shows that $z$ belongs to $F_{k,l}^d(\delta)$ too. In conclusion, we have $F_{k,l}^d (\delta)= R_{k,l}^d (\delta)\cap C_k^d (\delta)$ and as a consequence $F_{k,l}^d (\delta)$ is a convex compact set. Let us now show that $C_k^d (\delta)$ coincides with the set $S$ of points $z=(z_1,\ldots,z_d)\in (0,\infty)^d$ satisfying
\begin{eqnarray}\label{30sept1}
\frac{z_i}{k_i+\delta}  - \frac{z_j}{k_j+1-\delta} \geq 0 \qquad \forall i,j =1, \ldots,d \mbox{ with } i\neq j,
\end{eqnarray}
whose boundary satisfies
\begin{eqnarray}\label{29sept1}
\partial S \cap (0,\infty)^d = S \cap \left( \bigcup_{i,j=1,\ldots,d, i\neq j} (0,\infty) \cdot \left\{z\in  (0,\infty)^d \, \vert \, z_i=k_i+\delta, z_j=k_j+1-\delta \right\}\right).
\end{eqnarray}
First, we observe that $C_k^d (\delta)\subset S$. As a matter of fact, for every $z \in C_k^d (\delta)$ there are,  by (\ref{28sept1}), $\alpha>0$ and $y\in \mathcal{K}_k^d(\delta)$ such that $\alpha z =y$. Hence we have 
\[
k_i+\delta \leq y_i \leq k_i+1-\delta \qquad \forall j =i,\ldots,d,
\] 
which yields  for every $i,j\in \{1,\ldots,d\}$ with $i\neq j$,
\[
\frac{y_j}{\alpha(k_j+1-\delta)} \leq \frac{1}{\alpha} \leq \frac{y_i}{\alpha(k_i+\delta)},
\]
which means that $z \in S$. 

\begin{lemma}\label{LEM29sept}
Let $z\in C_k^d(\delta)$ be fixed, then $z\in \partial C_k^d(\delta)$ if and only if  there is $j\in \{1,\ldots,d\}$ such that $R(z)_j=k_j+1-\delta$.
\end{lemma}

\begin{proof}[Proof of Lemma \ref{LEM29sept}]
Let $z\in C_{k}^d(\delta)$ be fixed. If $R(z)_j\neq k_j+1-\delta$ for all $j =1,\ldots,d$, then we have $R(z)_j< k_j+1-\delta$ for all $j=1,\ldots,d$. As a consequence, by continuity of $R$ (see Lemma \ref{27sept1}) and since $R(R(z))=z$ (because $R(z)\in \partial E$), the image of a sufficiently small open ball $B$ centered at $R(z)$ satisfies
\[
R(B) \subset \partial E \cap \left(\bigcap_{j=1,\ldots,d} \{z_j<k_j+1-\delta\}\right) \subset L_k(\delta),
\]
so that $(0,+\infty)\cdot R(B)$ is a neighborhood of $z$ which is contained in $C_k^d(\delta)$. This shows that if $z\in \partial C_k^d(\delta)$, then there is  $j\in \{1,\ldots,d\}$ such that $z_j=k_j+1-\delta$. If  $R(z)_j=k_j+1-\delta$ for some $j\in \{1,\ldots,d\}$, then the segment $I:=[R(z),R(z)+se_j]$ with $s\in [0,1]$ (where $e_j$ denotes the $j$-th vector of the canonical basis of $\R^d$) satisfies $I\subset \partial E$ and $I\cap L_{k}^d(\delta)=\{R(z)\}$. This shows that $R(z)\in \partial C_k^d(\delta)$ and in turn that $z\in \partial C_k^d(\delta)$.
\end{proof}

Let us show that the boundary of $C_k^d(\delta)$ (in $(0,\infty)^d$) is contained in $\partial S \cap (0,\infty)^d$). If $z\in (0,\infty)^d \cap \partial C_k^d(\delta)$, then by Lemma \ref{LEM29sept}, there is $j\in \{1,\ldots \}$ such that $R(z)_j=k_j+1-\delta$. Thus, since $R(z)\in L_k^d(\delta)$, there  is also $i\in \{1,\ldots,d\}$ such that $R(z)_i=k_i+\delta$. We infer easily that $R(z)$ is in $\partial S$ and so $z$ too.  In conclusion, $C_k^d (\delta)\subset S$ and $(0,\infty)^d \cap \partial C_k^d(\delta) \subset (0,\infty)^d \cap \partial S$, so the two sets are equal. 

To prove (ii), we suppose for contradiction that $\mathcal{K}_0\subset \mathcal{B}_k/z$ for some $z\geq 1$ and $k\in\N$. Then we have $\mathcal{B}_k^-/z= (k-\delta)/z<\delta<1-\delta<(k+\delta)/z=\mathcal{B}_k^+/z$, which implies that $2\delta/z>1-2\delta$, a contradiction because $z\geq 1$ and $\delta\in (0,1/2)$. 

To prove (iii), we assume that $\delta \leq 1/4$ and consider some $z\in [1,+\infty)$. We set 
\[
 b^{\delta}(z) := \max \Bigl\{k\in \N \, \vert \, \mathcal{B}_k^-(\delta)/z \leq 1-\delta\Bigr\}=  \left\lfloor (1-\delta)z+\delta \right\rfloor.
\] 
We claim that $\mathcal{K}_0(\delta)\setminus (\mathcal{B}(\delta)/z)$ contains at least one $z$-bridge if and only if $\mathcal{B}_2^-/z\leq 1-\delta \Leftrightarrow b^{\delta}(z)\geq 2 \Leftrightarrow z\geq (2-\delta)/(1-\delta)$.  As a matter of fact, if $\mathcal{K}_0(\delta)\setminus (\mathcal{B}(\delta)/z)$ contains at least one $z$-bridge, then we have necessarily $\mathcal{B}_2^-/z\leq 1-\delta$. Moreover, if $b:=b^{\delta}(z)\geq 2$, then we have by definition of $b^{\delta}(z)$,
\[
z\geq \frac{b-\delta}{1-\delta} \quad \mbox{and} \quad z < \frac{b+1-\delta}{1-\delta}.
\]
So we have 
\[
\mathcal{B}^+_{b-1}/z = \frac{b-1+\delta}{z} > \frac{(b-1+\delta)(1-\delta)}{b+1-\delta} > \delta \quad (\mbox{because } \delta \leq 1/4),
\]
which means that the $z$-kwai $\mathcal{K}_{b-1}(\delta)$ is contained in $\mathcal{K}_0$. In conclusion, if $\mathcal{B}_2^-/z\leq 1-\delta \Leftrightarrow z\geq (2-\delta)/(1-\delta)$, then $\mathcal{K}_0(\delta)\setminus (\mathcal{B}(\delta)/z)$ contains a $z$-kwai so it contains a closed interval of length $(1-2\delta)/z$. Let us treat the case $z\in [1,(2-\delta)/(1-\delta)$. If $z\in [1,(2-\delta)/(1-\delta))$, then we have $\mathcal{B}_1^-(\delta)/z>\delta$ (because $(2-\delta)/(1-\delta)<(1-\delta)/\delta$) which implies
\[
\left[ \delta, \mathcal{B}_1^-(\delta)/z\right] \subset \mathcal{K}_0\setminus (\mathcal{B}/z),
\]
so that $\mathcal{K}_0(\delta)\setminus (\mathcal{B}(\delta)/z)$ contains a closed interval of length at least 
$$
 \mathcal{B}^-_1/z -\delta = \frac{1-\delta}{z}-\delta.
 $$
 To prove (iv), we note that for every $(z_1, \ldots, z_d) \in [1,\infty)^d$ and every $l \in \llbracket 0,N-1\rrbracket$, the set $P^d_{\delta,l} (z_1, \ldots, z_d)$ is compact and given by a finite set of affine inequalities, so it is a a convex polytope.

\end{document}